\documentclass[A4,11pt]{article}
\usepackage{amsmath}
\usepackage{amssymb}
\usepackage{amscd}

\newtheorem{theoreme}{Th\'eor\`eme}[subsection]
\newtheorem{proposition}[theoreme]{Proposition}
\newtheorem{corollaire}[theoreme]{Corollaire}
\newtheorem{lemme}[theoreme]{Lemme}
\newtheorem{definition}[theoreme]{D\'efinition}
\newtheorem{exemple}[theoreme]{Exemple}
\newtheorem{remarque}[theoreme]{Remarque}
\newtheorem{remarques}[theoreme]{Remarques}
\newenvironment{preuve}{\begin{trivlist} \item[]{\it Preuve---}}
{\par\hfill $\square$\end{trivlist}}
\renewcommand{\P}{\mathbb{P}}
\newcommand{\C}{\mathbb{C}}
\newcommand{\B}{{\cal B}}
\newcommand{\R}{\mathbb{R}}
\newcommand{\N}{\mathbb{N}}
\newcommand{\Z}{\mathbb{Z}}

\newcommand{\h}{{\rm h}}
\newcommand{\M}{{\cal M}}

\renewcommand{\H}{{\cal H}}

\newcommand{\K}{{\cal K}}
\newcommand{\HD}{{\rm HD}}
\newcommand{\PSH}{{\rm PSH}}

\newcommand{\Cr}{{\rm C}}
\newcommand{\PC}{{\rm PC}}
\newcommand{\PR}{{\cal P}}

\newcommand{\ddc}{{\rm dd^c}}
\renewcommand{\d}{{\rm d}}
\newcommand{\Dr}{{\rm D}}
\newcommand{\cad}{{\it c.-\`a-d. }}
\newcommand{\ie}{{\it i.e. }}
\newcommand{\lov}{{\rm lov}}
\newcommand{\vol}{{\rm vol}}
\newcommand{\loc}{{\rm loc}}
\newcommand{\Lone}{{{\rm L}^1}}
\newcommand{\Ltwo}{{{\rm L}^2}}
\newcommand{\Lp}{{{\rm L}^p}}
\newcommand{\Lploc}{{{\rm L}^p_{\rm loc}}}
\newcommand{\Linfty}{{{\rm L}^\infty}}

\newcommand{\Loneloc}{{{\rm L}^1_{\rm loc}}}
\newcommand{\Ltwoloc}{{{\rm L}^2_{\rm loc}}}
\newcommand{\dist}{{\rm dist}}
\newcommand{\diam}{{\rm diam}}

\renewcommand{\S}{{\cal S}}
\newcommand{\F}{{\cal F}}
\newcommand{\A}{{\cal A}}
\newcommand{\supp}{{\rm supp}}
\newcommand{\E}{{\cal E}}

\newcommand{\voir}{{\it voir }}

\renewcommand{\O}{{\rm O}}
\renewcommand{\o}{{\rm o}}
\newcommand{\tran}{{\, ^\tau\!}}
\newcommand{\Id}{{\rm id}}

\title{Dynamique des applications d'allure polynomiale}
\author{Tien-Cuong Dinh et Nessim Sibony}
\date{}
\begin{document}
\maketitle
\begin{abstract}
We study the dynamics of polynomial-like mappings in several variables.
A special case of our results is the following theorem.
\\
\it
Let $f:U\longrightarrow V$ be a proper holomorphic map
from an open set $U\subset \subset V$ onto a Stein manifold $V$. 
Assume $f$ is
of topological degree $d_t\geq 2$. Then there
is a probability measure $\mu$ supported on $\K:=\bigcap_{n\geq 0}f^{-n}(V)$
satisfying the following properties.
\begin{enumerate}
\item[{\rm 1.}] The measure $\mu$ is invariant, K-mixing, of maximal entropy
  $\log d_t$. 
\item[{\rm 2.}] If $J$ is the Jacobian of $f$ with respect to a volume form
  $\Omega$ then $\int \log J \d \mu \geq
  \log d_t$. 
\item[{\rm 3.}] For every probability measure $\nu$ on $V$ with no mass on
  pluripolar sets $d_t^{-n} (f^n)^* \nu \rightharpoonup \mu$.
\item[{\rm 4.}]
If the p.s.h. functions on $V$ are $\mu$-integrables ($\mu$ is
  PLB) then
\begin{enumerate}
\item[{\rm (a)}] The Lyapounov exponents for $\mu$ are strictly positive.
\item[{\rm (b)}] $\mu$ is exponentially mixing.
\item[{\rm (c)}] There is a proper analytic subset $\E$ of $V$ such that for
  $z\not\in\E$, $\mu^z_n:=d_t^{-n} (f^n)^*\delta_z \rightharpoonup
  \mu$. 
\item[{\rm (d)}] The measure $\mu$ is a limit of Dirac masses on the repelling
  periodic points.
\end{enumerate}
\end{enumerate}
The condition $\mu$ is PLB is stable under small pertubation of
$f$. This gives large families where it is satisfied.
\rm
\end{abstract}
{\bf Mots cl\'es:} application d'allure polynomiale, mesure d'equilibre,
  K-m\'elange, vitesse de m\'elange, exposant de Lyapounov.
\\
{\bf Classification math\'ematique:} 37F, 32H50, 32Q, 32U.
\section{Introduction}
Depuis une vingtaine d'ann\'ees, l'\'etude de la dynamique des
applications holomorphes a connu une grande activit\'e. Pour la
th\'eorie des applications rationnelles dans $\P^1$, l'utilisation du
th\'eor\`eme de Riemann mesurable, des th\'eor\`emes de distorsion et
la technique des modules d'anneaux sont les outils fondamentaux qui
ont permis la d\'emonstration des th\'eor\`emes de non errance de
Sullivan et des
progr\`es dans les probl\`emes de renormalisation. Pour les aspects les
plus \'el\'ementaires, on pourra consulter les ouvrages de
Carleson-Gamelin \cite{CarlesonGamelin} et Milnor \cite{Milnor}.
\par
Les outils de la th\'eorie de Fatou-Julia \`a une variable complexe,
particuli\`erement le th\'eor\`eme de Montel,
n'admettent pas une extension imm\'ediate pour traiter les probl\`emes
analogues en plusieurs variables. L'utilisation de la th\'eorie des
courants positifs ferm\'es et de la th\'eorie du potentiel s'est
revel\'ee utile dans nombre de questions concernant les
automorphismes de $\C^2$, les endomorphismes de $\P^k$ ou plus
g\'en\'eralement la dynamique des applications m\'eromorphes. Les
trois articles Bedford-Smillie \cite{BedfordSmillie1}, 
Forn\ae ss\cite{Fornaess} et \cite{Sibony2}
contiennent un panorama des questions trait\'ees ainsi qu'une
importante bibliographie. On trouve dans \cite{FornaessSibony4}
l'\'etude dynamique d'exemples non triviaux de $\P^2$.
\par
Dans le pr\'esent article, on \'etudie le probl\`eme suivant. Soit $V$
une vari\'et\'e complexe de dimension $k$. 
Pour simplifier, supposons que $V$ est de Stein.
On consid\`ere une
application holomorphe d\'efinie dans un ouvert $U\subset\subset V$ et
telle que $f:U\longrightarrow V$ soit un rev\^etement (ramifi\'e ou
non) au dessus de $V$. On note $d_t$ le degr\'e topologique de $f$. Il
s'agit d'\'etudier la dynamique de $f$ \`a l'aide de mesures
ou de courants invariants associ\'es \`a $f$.
\par
Cette situation (applications d'allure polynomiale) est
stable par pertubation et est tr\`es riche en exemples. En dimension
1, elle a fait l'objet d'un travail de Douady et Hubbard
\cite{DouadyHubbard}. L'outil essentiel en dimension 1 est le
th\'eor\`eme de Riemann mesurable qui permet d'en ramener l'\'etude \`a
celle des polyn\^omes \`a une variable. Il est improbable qu'en
dimension strictement sup\'erieure \`a un, les applications que nous
\'etudions soient conjugu\'ees \`a des applications polynomiales comme
c'est le cas \`a une variable.
\par
Notre but est de construire une mesure de probabilit\'e invariante
maximisant l'entropie et dont les exposants de Lyapounov soient
strictement positifs, ensuite d'\'etudier les propri\'et\'es de cette
mesure. En somme il s'agit de construire une mesure hyperbolique. Pour
les automorphismes de H\'enon dans $\C^2$, ce programme a \'et\'e
r\'ealis\'e par Bedford, Lyubich, Smillie
\cite{BedfordLyubichSmillie,BedfordSmillie2}
pour une mesure introduite par le second auteur du pr\'esent
article. On sait que la construction d'un tel objet  est une question
centrale en dynamique. Elle est difficile dans le cadre r\'eel (\voir
les travaux de Benedicks-Carleson-Young
\cite{BenedicksCarleson, BenedicksYoung}). Le cadre holomorphe facilite
grandement les choses et fournit de vastes classes d'exemples. 
\par
Pour plus de clart\'e quant aux techniques et afin de mettre en
\'evidence les propri\'et\'es du cas holomorphe, nous introduisons les
mesures d'\'equilibre qui nous int\'eressent dans le cadre riemannien.
\par
Soit $V$ une vari\'et\'e riemannienne munie d'une forme
volume $\Omega$.
Soit $f$ 
une application r\'eelle 
de classe ${\cal C}^1$
d\'efinissant un rev\^etement ramifi\'e de degr\'e $d_t\geq 2$
au voisinage d'un compact
$X$ de $V$.
On suppose que $X$ est de mesure positive et que 
$f^{-1}(X)\subset X$. 
On se propose d'introduire
une ``mesure d'\'equilibre'' sur $X$. 
\par
Plus pr\'ecis\'ement, soit $J$ le jacobien r\'eel
de $f$ d\'efini par
$f^*\Omega=J\Omega$; on suppose $J$ non n\'egatif, \cad que $f$
pr\'eserve l'orientation en dehors de l'ensemble critique. 
Notons $\Omega_{|X}$ la restriction
normalis\'ee  de la
forme de volume $\Omega$ \`a $X$.
Lorsque $f$
d\'efinit un rev\^etement ramifi\'e de $f^{-1}(X)$ au dessus de $X$,
toute mesure $\sigma$, valeur
d'adh\'erence de la suite
$$\sigma_N:=\frac{1}{N}\sum_{n=1}^N\frac{(f^n)^*\Omega_{|X}}{d_t^n},$$
est invariante et v\'erifie $f^*\sigma=d_t\sigma$.
De plus, la mesure $\sigma$ ne
charge pas l'ensemble critique $\Cr$ de $f$, on a m\^eme
$\int \log J \d\sigma\geq \log d_t$. On en d\'eduit que la mesure
$\sigma$ est d'entropie au moins $\log d_t$.
\par
Notons $\mu$ une mesure de probabilit\'e obtenue par la construction
ci-dessus. On se pose le probl\`eme des propri\'et\'es dynamiques 
de $\mu$:
est-elle l'unique mesure d'entropie maximale, quels sont ses exposants
de Lyapounov, les points p\'eriodiques sont-ils \'equidistribu\'es par
rapport \`a $\mu$? Dans le cadre r\'eel, il est facile de construire
des contres exemples en prenant des produits de vari\'et\'es. 
Dans le cadre complexe, la situation est
totalement diff\'erente.
\par
Lorsque $V=\P^k$ et $f$ est un endomorphisme holomorphe de degr\'e
$d_t>1$, Briend-Duval \cite{BriendDuval1}, \cite{BriendDuval2}, ont
r\'ecemment montr\'e que les mesures $\mu_n^z$ \'equidistribu\'ees aux
pr\'eimages de $z$
$$\mu^z_n:=\frac{1}{d_t^n} \sum_{f^n(w)=z}\delta_w$$
convergent vers une mesure $\mu$ pour tout $z$ n'appartenant pas \`a un
ensemble analytique $\E$. Ils ont \'egalement montr\'e que la mesure
$\mu$ est l'unique mesure d'entropie maximale $\log d_t$, que les
exposants de Lyapounov sont strictement positifs et que les points
p\'eriodiques r\'epulsifs sont \'equidistribu\'es par rapport \`a
$\mu$. Ant\'erieurement
\cite{FornaessSibony2}, Forn\ae ss et le second auteur avaient
montr\'e que la mesure $\mu$ est m\'elangeante dans
$\P^k$ et que l'ensemble exceptionnel $\E$ est pluripolaire.
Le cas de dimension 1 avait d\'ej\`a \'et\'e r\'esolu par Lyubich
\cite{Lyubich} et Freire-Lopes-Ma\~ne \cite{FreireLopesMane}. Les
techniques utilis\'es sont d'ailleurs assez proches. 
\par
Ce sont ces r\'esultats qui ont servi de point de d\'epart \`a notre
\'etude. Observons cependant que le cadre des applications d'allure
polynomiale est moins rigide que celui des endomorphismes holomorphes
de $\P^k$. Il le contient d'ailleurs car il suffit de consid\'erer le
relev\'e d'un endomorphisme de $\P^k$ \`a $\C^{k+1}$ qui est alors une
application d'allure polynomiale pour $U$, $V$ convenables.
\par
Le paragraphe 3 est consacr\'e aux applications d'allure
polynomiale. Soit $f:U\longrightarrow V$ un rev\^etement ramifi\'e
holomorphe de degr\'e $d_t>1$, $U\subset\subset V$.
On d\'efinit
$\K:=\bigcap_{n\geq 0} f^{-n}(V)$. On montre qu'il existe une mesure
de probabilit\'e $\mu$ port\'ee par $\partial \K$ d'entropie $\log
d_t$, $K$-m\'elangeante donc en particulier m\'elangeante de tout
ordre. En particulier, si $B$ est un bor\'elien tel que $\mu(B)>0$
alors $\lim_{n\rightarrow\infty} \mu(f^n(B))=1$.
La d\'emonstration est bas\'ee sur des propri\'et\'es de
convergence des fonctions plurisousharmoniques.
\par
Suivant l'id\'ee de
Gromov-Yomdin \cite{Gromov1,Gromov3,Yomdin}
(\voir \'egalement \cite{Friedland}), nous consid\'erons 
dans le cas non compact, les 
degr\'es dynamiques de l'application $f$ qui d\'ecrivent la croissance
des volumes des sous-vari\'et\'es par it\'eration. Pour $1\leq l\leq
k$, on pose pour une forme de K\"ahler $\omega$ sur $V$
$$d_l:=\limsup_{n\rightarrow \infty} \left(\int_U (f^n)_* \omega^{k-l}
  \wedge \omega^l \right)^{1/n}.$$
Au paragraphe 3.3, nous estimons les degr\'es dynamiques de
  l'application $f$. 
On montre que $d_l\leq d_t$ pour tout $1\leq l\leq k$ et 
que $f$ est d'entropie $\log d_t$.
Donc la mesure $\mu$ maximise l'entropie. 
\par
Soit $X$ un sous-ensemble analytique de $V$. 
Le paragraphe 3.4 donne une caract\'erisation du plus grand sous ensemble
analytique $\E_X$ de $X$ qui est 
totalement invariant, \ie $f^{-1}(\E_X) = \E_X\cap U$. Un point $z$ est dans
$\E_X$ si et seulement si la ``proportion'' d'orbites de points de
$f^{-n}(z)$ restant dans $X$ est strictement positive. Ce r\'esultat
est vrai pour toute vari\'et\'e complexe $V$.
\par
De m\^eme que la croissance du volume des it\'er\'es de sous-vari\'et\'es est
tr\`es li\'ee \`a l'entropie, la croissance du volume des it\'er\'es 
de l'ensemble
critique $\Cr$ de $f$ est li\'ee \`a l'\'etude fine de la
mesure $\mu$. Pour d\'ecrire cette croissance, consid\'erons
$\delta_n$ le volume normalis\'e des images $f^n(\Cr\cap
U_{-n})$. Plus pr\'ecis\'ement, 
$$\delta_n:=\int_{\Cr\cap f^{-n}(U)}\frac{(f^n)^* \omega^{k-1}}{d_t^n} 
\ \ \mbox{ et } \ \ \delta:=\limsup_{n\rightarrow \infty}
\sqrt[n]{\delta_n}.$$
On a alors l'analogue du th\'eor\`eme de Briend-Duval pour l'espace
projectif. Si la s\'erie 
$\sum \delta_n$ converge (en particulier si $\delta<1$), 
alors pour tout $z$ hors d'un
ensemble analytique $\E$, la suite $\mu^z_n$ converge vers $\mu$. Les
exposants de Lyapounov sont non n\'egatifs et minor\'es par
$\frac{1}{2}\log(d_t/d_{k-1})$. La mesure $\mu$ est limite de masses
de Dirac en des points p\'eriodiques r\'epulsifs.
\par
Afin d'\'etudier la croissance des $\delta_n$, nous sommes amen\'es
\`a introduire une propri\'et\'e analytique des mesures de
probabilit\'e. Nous dirons que $\nu$ est PLB si toutes les fonctions
p.s.h. sont $\nu$-int\'egrables. En une variable, cela \'equivaut \`a
dire que $\nu$ est \`a Potentiel Localement Born\'e. Lorsque la mesure
d'\'equilibre $\mu$ de $f$ est PLB, on a un contr\^ole de $\delta_n$
et les exposants de Lyapounov sont strictement positifs. De plus, la
mesure $\mu$ est m\'elangeante d'ordre exponentiel. Nous donnons
une estimation de l'ordre de m\'elange qui est nouvelle m\^eme dans le
cas des applications holomorphes de $\P^k$. 
Nous montrons que
lorsqu'une application a une mesure associ\'ee $\mu$ qui est PLB, il
en est de m\^eme pour les applications \`a allure polynomiale
voisines. Cela permet de construire de vastes classes d'exemples
satisfaisant nos hypoth\`eses.   
\par
Cet article reprend une version pr\'ec\'edente de juin 2001 et une
partie d'une
pr\'epublication d'Orsay de mars 2002. Dans un prochain travail, nous
donnerons une construction de la mesure d'\'equilibre comme produit
g\'en\'eralis\'e de courants positifs ferm\'es pour les applications
polynomiales de $\C^k$ et pour les endomorphismes d'une vari\'et\'e
complexe compacte. Les id\'ees de construction de la mesure
d'\'equilibre peuvent \^etre \'etendues au cadre de l'it\'eration
al\'eatoire.  
\par
C'est un plaisir de remercier A. Ancona qui a r\'epondu \`a 
plusieurs questions de th\'eorie du potentiel.
\section{Applications r\'eelles}
Dans ce paragraphe, nous donnons quelques propri\'et\'es abstraites
sur le m\'elange pour les mesures
invariantes associ\'ees \`a des rev\^etements ramifi\'es. 
Dans le cas riemannien, les mesures d'\'equilibre que nous construisons, ne
chargent pas l'ensemble critique. 
Une version quantitative de cette propri\'et\'e
permet de montrer qu'un exposant de Lyapounov au moins
est strictement positif.
\subsection{Rev\^etements ramifi\'es}
Soient $X$ et $Y$ deux espaces m\'etriques localement compacts. 
Soit $f:Y\longrightarrow X$
une application continue. Pour toute fonction continue $\varphi$
sur $X$, on peut d\'efinir une fonction continue
$f^*\varphi:=\varphi\circ f$ sur $Y$. Par dualit\'e, pour toute
mesure $\nu$ \`a support compact dans $Y$, on peut d\'efinir une
mesure $f_*\nu$ \index{$f^*\nu$}
\`a support compact dans $X$ par la relation
$$\int_X \varphi \d (f_*\nu) := \int_Y f^*\varphi \d \nu \ \ \ \mbox{ pour
toute } \varphi \mbox{ continue sur } X.$$
Cet op\'erateur est continu sur l'ensemble des mesures positives
\`a support compact dans $Y$.
\par
En g\'en\'eral, on ne peut pas d\'efinir 
l'op\'erateur $f^*$ sur l'ensemble des mesures. Nous allons  
donner un cadre o\`u cet op\'erateur est bien d\'efini.
\begin{definition} \rm
Nous dirons que $(Y,f,X)$ est {\it un espace
\'etal\'e ramifi\'e au dessus de $X$} 
s'il existe une fonction $n:
Y\longrightarrow \N^+$ v\'erifiant les propri\'et\'es suivantes:
\begin{enumerate}
\item Pour tout $x\in X$, l'ensemble $f^{-1}(x)$ est discret.
\item $n(y)=1$ dans un ouvert dense de $Y$.
\item Pour tout $y_0\in Y$ et $Y_0$ 
un voisinage suffisamment petit de $y_0$, on a
$$\sum_{y\in Y_0,\ f(y)=x} n(y)= n(y_0)$$
lorsque $x\in X$ est suffisamment proche de $f(y_0)$.
\end{enumerate} 
On dira que $n(y)$ est {\it la multiplicit\'e}
de $f$ en
$y$.
\end{definition}
\par
Pour toute fonction $\varphi$ \`a
support compact dans $Y$, on d\'efinit sur $X$ la fonction
$f_*\varphi$ par la formule 
$$f_*\varphi(x):=\sum_{f(y)=x} n(y)\varphi(y).$$
Nous laissons au lecteur la preuve de la proposition suivante 
qui donne les premi\`eres propri\'et\'es de cette
notion. 
\begin{proposition} Soit $f:Y\longrightarrow X$ une application
continue d\'efinissant un espace \'etal\'e comme
pr\'ec\'edemment. Alors
\begin{enumerate}
\item L'application $f$ est ouverte. La fonction $n$ est
semi-continue sup\'erieurement.
\item La fonction $n$ est l'unique fonction v\'erifiant les
conditions donn\'ees dans la d\'efinition 2.1.1.
\item La condition 3 de la d\'efinition 2.1.1 \'equivaut \`a la
condition suivante:
pour toute $\varphi$ continue, \`a support compact
dans $Y$, $f_*\varphi$ est continue.
\end{enumerate}
\end{proposition}
\par
La propri\'et\'e 3 de la proposition 2.1.2 permet de d\'efinir
pour toute mesure $\nu$ \`a support dans $X$, 
une mesure $f^*\nu$ 
de $Y$ par la relation
$$\int\varphi \d f^*\nu=\int f_*\varphi \d\nu$$
o\`u $\varphi\in {\cal C}_c(X)$.
L'op\'erateur $f^*$ est continu sur les mesures.
\begin{definition} \rm
Soit $(Y,f,X)$ un espace \'etal\'e ramifi\'e.
On dira que $f$ d\'efinit {\it un rev\^etement ramifi\'e} 
de {\it degr\'e} 
$d_t$ si pour tout
$x\in X$ on a 
$$\sum_{f(y)=x} n(y)=d_t$$
\end{definition}
\par
Observons que si $f$ d\'efinit un rev\^etement ramifi\'e de
degr\'e $d_t$, pour
toute fonction $\varphi$ \`a support dans $Y$ (compact ou non),
la fonction $f_*\varphi$ est bien d\'efinie. De plus, si $\nu$ est
une mesure positive de masse $m$ sur $X$, la mesure $f^*\nu$ est
de masse $d_tm$. On a aussi
$f_*(f^*)\nu= d_t\nu$. 
\subsection{K-m\'elange pour les rev\^etements ramifi\'es}
Soit $f:Y\longrightarrow X$ une application continue
d\'efinissant un rev\^etement ramifi\'e de degr\'e $d_t$.
Consid\'erons le cas o\`u $Y$ est un ouvert relativement compact de
$X$. L'op\'erateur $d_t^{-1} f^*$ est continu sur le 
convexe des mesures de probabilit\'e \`a support dans $\overline
Y$. 
Le th\'eor\`eme du point fixe entra\^{\i}ne l'existence d'une
mesure de probabilit\'e $\mu$ \`a support compact 
dans $X$ v\'erifiant 
$$f^*\mu=d_t\mu.$$
On a
$$f_*\mu= f_* \left(\frac{1}{d_t} f^*\mu\right)
=\frac{1}{d_t}f_*(f^*)\mu= \mu.$$
La mesure $\mu$ est donc invariante par $f_*$.
\par
Pour toute fonction $\varphi$ continue sur $X$, on d\'efinit la
fonction $\Lambda\varphi$ 
continue sur $X$ par la relation 
$$\Lambda\varphi(x):=\frac{f_*\varphi(x)}{d_t}=
\frac{1}{d_t}\sum_{f(y)=x}
n(y) \varphi(y).$$
L'op\'erateur $\Lambda$ est parfois appel\'e
op\'erateur de Perron-Frobenius \cite{Lyubich}.
Il se prolonge par continuit\'e en un op\'erateur
de $\Ltwo(\mu)$ dans lui-m\^eme. D'apr\`es l'in\'egalit\'e de
Cauchy-Schwarz, on a $|f_*\varphi|^2\leq d_t f_* |\varphi|^2$. Par
cons\'equent, $\|\Lambda\|=1$. 
L'op\'erateur adjoint de $\Lambda$ est
d\'efini par $(\tran\Lambda)(\varphi)=\varphi\circ f$.
On a bien s\^ur
$\Lambda \tran\Lambda=\Id$. Mais $\Lambda$ n'est pas injectif en
g\'en\'eral. Pour
tout $n\geq 1$, posons 
$$V_n:=\big\{\varphi\in\Ltwo(\mu)|\ \Lambda^n\varphi=0\big\}.$$
On v\'erifie que la suite $V_n$ est croissante et que
$$V_n^\perp=\big\{\theta|\ \theta=\psi_n\circ f^n, \
\psi_n\in\Ltwo(\mu)\big\}.$$
Notons $H_0$ l'adh\'erence de $\bigcup_{n\geq 1} V_n$. Soit
$H_0^\perp$ son orthogonal. Il est clair que 
$$H_0^\perp=\big\{\theta| \mbox{ pour tout } n\geq 1 \mbox{ il existe }
\psi_n\in\Ltwo(\mu) \mbox{ v\'erifiant } \theta=\psi_n\circ f^n\big\}.$$
L'op\'erateur $\Lambda$ est injectif sur $H_0^\perp$. Il en r\'esulte
que la restriction de $\Lambda$ sur $H_0^\perp$ est un op\'erateur
unitaire. Autrement dit, pour $\varphi\in H_0^\perp$
$$(\Lambda\varphi)\circ f=\varphi.$$
Notons $\A$ la famille des bor\'eliens $A$ v\'erifiant
$$f^{-n}(f^n(A))=A \mbox{ pour tout } n\geq 1.$$
On v\'erifie que $\A$ est une tribu, en effet,  
$H_0^\perp$ est engendr\'e par les fonctions
indicatrices $1_A$ avec $A\in\A$. Rappelons les notions de K-m\'elange
et de $r$-m\'elange
\cite{CornfeldSinai}.
\begin{definition} \rm Soit $\nu$ une mesure de probabilit\'e
  invariante pour $f$. On dira
que $\nu$ est
{\it K-m\'elangeante} 
si
$$\sup_{\|\varphi\|_{\Ltwo(\nu)}\leq 1}
\left| \int \varphi(f^n)\psi\d \nu -
\left(\int \varphi\d \nu\right) \left(\int \psi\d
  \nu\right)\right|\longrightarrow 0$$
quand $n\rightarrow \infty$
et que $\nu$ est {\it $r$-m\'elangeante ou m\'elangeant d'ordre $r$} 
si
$$\lim_{n_1,...,n_r\rightarrow\infty} \int\psi_0\psi_1(f^{n_1})\ldots 
  \psi_r(f^{n_1+n_2+\cdots+n_r})\d \nu \longrightarrow
  \prod_{i=0}^r \left(\int \psi_i\d \nu \right)$$
o\`u $\varphi$, $\psi$ sont dans $\Ltwo(\nu)$ et les $\psi_i$ sont
  dans $\Linfty(\nu)$.
\end{definition}
Dans le cadre de la dynamique des rev\^etements ramifi\'es, 
on a la proposition suivante adapt\'ee \`a nos besoins.
\begin{proposition} Soient $X$ un espace m\'etrique localement
  compact, $Y$ un ouvert relativement compact de $X$.
Soit $(Y,f,X)$ un rev\^etement  ramifi\'e de degr\'e
$d_t\geq 2$. Si $\mu$ est une mesure de probabilit\'e \`a support
compact dans $X$ v\'erifiant
$f^*\mu=d_t\mu$, alors les propri\'et\'es suivantes
sont \'equivalentes:
\begin{enumerate}
\item $H_0^\perp=\C$.
\item Pour toute $\varphi\in\Ltwo(\mu)$,
  $\Lambda^n\varphi\longrightarrow c_\varphi:=\int\varphi\d \mu$ dans
  $\Ltwo(\mu)$.
\item Si $A$ est un bor\'elien v\'erifiant $f^{-n}(f^n(A))=A$ pour
  tout  $n\geq 1$, alors $\mu(A)=0$ ou $1$.
\item Si $B$ est un bor\'elien v\'erifiant $\mu(B)>0$ alors
  $\lim_{n\rightarrow \infty}
  \mu(f^n(B))=1$. 
\item  La suite d'op\'erateurs $(\Lambda^n)$ de $\Ltwo(\mu)$ dans
  $\Ltwo(\mu)$ est convergente au sens fort et la mesure $\mu$
est K-m\'elangeante.
\item La suite d'op\'erateurs $(\Lambda^n)$ de $\Ltwo(\mu)$ dans
  $\Ltwo(\mu)$ est 
convergente au sens fort et la mesure $\mu$ est
$r$-m\'elangeante pour tout $r\geq 1$.
\item La suite d'op\'erateurs $(\Lambda^n)$ de $\Ltwo(\mu)$ dans
  $\Ltwo(\mu)$ est convergente au sens fort
  et la mesure $\mu$ est ergodique.  
\end{enumerate}
\end{proposition}
\begin{preuve}
1.$\Longrightarrow$ 2. Si $\varphi_0\in H_0$ on a
$\Lambda^n\varphi_0\longrightarrow 0$
car $\|\Lambda\|= 1$. Sur $H_0^\perp$, on a $\Lambda^n 1=1$. Il
suffit donc de d\'ecomposer $\varphi=\varphi_0+c$, $\varphi_0\in H_0$,
$c\in \C$.
\par
2.$\Longrightarrow$ 3. Si $A=f^{-n}(f^n(A))$, posons $A_n:=f^n(A)$. 
On a $1_A=1_{A_n}\circ f^n$ et $\Lambda^n 1_A=1_{A_n}$. La suite
$1_{A_n}$ convergeant dans $\Ltwo(\mu)$ vers une constante, 
la limite ne peut \^etre que 0 ou 1.
D'o\`u $\mu(A)=\mu(A_n)=0$ ou $1$. 
\par
3.$\Longrightarrow$ 4. Soit $B$ tel que $\mu(B)>0$. Posons
$B_n:=f^{-n}(f^n(B))$, on a $B_{n+1}\supset B_n$. Si $\tilde
B:=\bigcup_{n\geq 0} B_n$ on a $f^{-n}(f^n(\tilde B))=\tilde B$. Donc 
$$\mu(f^n(B))=\mu(f^{-n}(f^n(B))=\mu(B_n)\longrightarrow 
\mu(\tilde B)>0.$$
D'apr\`es 3., on a $\mu(\tilde B)=1$.
\par
4.$\Longrightarrow$ 3. Puisque $f^{-n}(f^n(A))=A$ et $f^*\mu=d_t\mu$, 
on a
\begin{eqnarray*}
\mu(A) & = & \langle \mu,1_A \rangle 
= \langle\mu,1_{f^{-n}(f^n(A))}\rangle = \langle \mu,1_{f^n(A)}\circ
  f^n\rangle\\
& = & \langle \mu,1_{f^n(A)} \rangle = \mu(f^n(A)).
\end{eqnarray*}
Si $\mu(A)>0$, d'apr\`es 4.,
$\mu(f^n(A))$ tend vers 1 quand
  $n\rightarrow\infty$. Par cons\'equent, on a $\mu(A)=0$ ou 1.
\par 
3.$\Longrightarrow$ 1.  Puisque $H_0^\perp$ est engendr\'e par les
fonctions $1_A$ avec $A\in \A$, il est clair que $H_0^\perp=\C$. 
\par
1. et 2. $\Longrightarrow$ 5. Fixons $\varphi\in\Ltwo(\mu)$ et
$\epsilon>0$. Il existe $n_0>0$, $\tilde\varphi\in V_{n_0}$ et
$c\in H_0^\perp=\C$ tels que $\|\varphi-\tilde\varphi-c\|_{\Ltwo(\mu)}\leq
\epsilon$. Pour tout $n\geq n_0$ on a
$$\|\Lambda^n\varphi-c\|_{\Ltwo(\mu)} =
\|\Lambda^n(\varphi-\tilde\varphi-c)\|_{\Ltwo(\mu)} \leq \epsilon$$
car $\Lambda^n\tilde\varphi=0$ et 
$\|\Lambda\|=1$. D'o\`u la convergence de la
suite $(\Lambda^n)$.
\par
Pour le K-m\'elange, il
suffit de remarquer que la propri\'et\'e 2 implique 
\begin{eqnarray*}
\left| \int \varphi(f^n)\psi\d \mu -
\left(\int \varphi\d \mu\right) \left(\int \psi\d
  \mu\right)\right| & = & 
\left|\int\varphi(\Lambda^n\psi-c_\psi)\d\mu\right|\\
& \leq & 
\|\varphi\|_{\Ltwo(\mu)} \|\Lambda^n\psi -c_\psi\|_{\Ltwo(\mu)}.
\end{eqnarray*}
\par
5.$\Longrightarrow$ 6. 
On montre le $r$-m\'elange par r\'ecurrence sur
$r$. Le $1$-m\'elange est une cons\'equence du K-m\'elange. Pour
$r>1$, on a
$$\int c_{\psi_0}\psi_1(f^{n_1})\ldots 
  \psi_r(f^{n_1+n_2+\cdots+n_r})\d \mu = c_{\psi_0}\int \psi_1\ldots 
  \psi_r(f^{n_2+\cdots+n_r})\d \mu.$$ 
Par hypoth\`ese de r\'ecurrence, le dernier terme tend vers 
$$c_{\psi_0}\ldots c_{\psi_r}=\prod_{i=0}^r
\left(\int \psi_i\d \mu \right).$$
On a aussi que
\begin{eqnarray*}
\lefteqn{\int (\psi_0-c_{\psi_0})\psi_1(f^{n_1})\ldots 
  \psi_r(f^{n_1+n_2+\cdots+n_r})\d \mu} \hspace{3cm}\\
& = & \int
  (\Lambda^{n_1}\psi_0-c_{\psi_0}) \psi_1\ldots 
  \psi_r(f^{n_2+\cdots+n_r})\d \mu \\
& \leq &
  \|\Lambda^{n_1}\psi_0-c_{\psi_0}\|_{\Ltwo(\mu)}
\|\psi_1\|_{\Linfty(\mu)} 
  \ldots \|\psi_r\|_{\Linfty(\mu)}.
\end{eqnarray*}
Le dernier terme tend vers $0$ et on obtient le r\'esultat.
\par
6.$\Longrightarrow$ 7. est clair.
\par
7.$\Longrightarrow$ 3. Soit $A\in \A$. Posons $A_n:=f^n(A)$. 
Il existe une fonction $\theta\in H_0^\perp$
telle que $\Lambda^n 1_A=1_{A_n}\longrightarrow \theta$ avec
$\Lambda\theta=\theta$. 
Puisque $\Lambda$ est inversible dans $H_0^\perp$, on a
$\theta\circ f=\theta$. L'ergodicit\'e de $\mu$ implique que $\theta$
est constante. Cette constante ne peut \^etre que 0 ou 1. On a
$$\mu(A)=\langle \mu,1_A \rangle = \langle \mu,\Lambda^n 1_A \rangle
= \langle \mu,1_{A_n} \rangle \longrightarrow
\langle \mu,\theta \rangle = 0 \mbox{ ou } 1.$$
\end{preuve}
\par
Notons $\tilde\M$ 
l'ensemble des mesures
de probabilit\'e $\nu$ \`a support
compact dans $X$ telles que pour tout $n\geq 1$ il existe une mesure
$\nu_n$  \`a support compact
dans $X$ satisfaisant la relation
$\nu=d_t^{-n}(f^n)^*\nu_n$.
C'est un convexe compact.
Notons \'egalement $\M\subset \tilde{\M}$
le convexe compact des mesures de probabilit\'e $\nu$  \`a support
compact dans $X$ v\'erifiant $f^*\nu=d_t\nu$.
Si $\nu'$ est une mesure de
probabilit\'e \`a support compact
dans $X$ et si $\nu$ est une valeur adh\'erente \`a la suite 
$d_t^{-n}(f^n)^*\nu'$ alors $\nu$ appartient
\`a $\tilde{\M}$. 
\begin{corollaire} Soient $f$ et $\mu$
v\'erifiant les propri\'et\'es de la proposition 2.2.2. 
Alors pour toute mesure $\nu\in \tilde{\M}$,
il existe une constante $0\leq
  c\leq 1$ et une mesure $\nu^s$ singuli\`ere par rapport \`a $\mu$ 
telles que $\nu=c\mu +\nu^s$. En particulier, $\mu$ est
extr\'emale dans $\M$.
\end{corollaire}
\begin{preuve} Si $\nu=d_t^{-n}(f^n)^*\nu_n$,
on peut \'ecrire $\nu=c\mu+\nu^s$ et $\nu_n=c_n\mu+\nu^s_n$ 
o\`u $\nu^s$, $\nu^s_n$ sont
singuli\`eres par rapport \`a $\mu$ et o\`u $c$, 
$c_n$ sont des fonctions
positives dans $\Lone(\mu)$. On a
$$\nu=\frac{(f^n)^*\nu_n}{d_t^n}=c_n(f^n)\mu +
\frac{(f^n)^*\nu^s_n}{d_t^n}.$$ 
Puisque $f^*\mu=d_t\mu$, la mesure
$d_t^{-n}(f^n)^*\nu^s_n$ est singuli\`ere par
rapport \`a $\mu$ et $c_n(f^n)\in\Lone(\mu)$. On
en d\'eduit que $c=c_n\circ f^n$. D'apr\`es la proposition 2.2.2,
la fonction $\min\{c,\alpha\}$ de $H_0^\perp$ est constante
pour tout $\alpha\geq
0$. Par suite, $c$ est constante. Les mesures $\mu$ et $\nu$ \'etant de masse
1, on a $0\leq c\leq 1$.
\end{preuve}
\par
Nous verrons que l'hypoth\`ese faite dans la
proposition 2.2.2 est toujours v\'erifi\'ee pour les applications
holomorphes d'allure polynomiale. 
Donnons cependant un exemple dans le cadre riemannien.
\par
Soient 
$U\subset\subset V$ deux ouverts simplement connexes
d'une vari\'et\'e riemannienne $M$. On suppose que $V$ est connexe et
que $U$ contient $d_t\geq 2$ composantes connexes. 
Soit $f:U\longrightarrow V$ une application de
classe ${\cal C}^1$ d\'efinissant des bijections entre chaque
composante de $U$ et $V$.
On suppose que $\K:=\bigcap_{n\geq 0} f^{-n}(V)$ n'a qu'un nombre au
plus d\'enombrable
de composantes connexes non r\'eduites \`a un point. Ceci est vrai en
particulier si $V$ est de dimension 1. 
Montrons que la mesure $\mu$, obtenue comme point fixe de l'op\'erateur
$d_t^{-1}f^*$, 
est K-m\'elangeante et pour tout $x\in V$ on a
$\mu^x_n:=d_t^{-n}(f^n)^*\delta_x\rightharpoonup
\mu$ 
ce qui montre en particulier l'unicit\'e de $\mu$. 
\par
Soient $(f^{-n}_i)$ les branches inverses de $f^n$ avec $1\leq i\leq
d_t^n$. Posons $R^n_i:=f^{-n}_i(V)$. Pour tout $i$, il existe $j$
tel que $R^n_i\subset\subset R^{n-1}_j$. Soit $\B$ la
$\sigma$-alg\`ebre engendr\'ee par les $R^n_i$. D'apr\`es un
th\'eor\`eme de Blackwell \cite[p. 62]{Meyer},
l'alg\`ebre $\B$ contient
toute alg\`ebre $\B'$ dont les atomes (\cad les bor\'eliens minimaux) 
sont r\'eunion d'atomes de $\B$.
Il nous faut donc d\'eterminer les atomes de $\B$. Chaque $R^n_i$
\'etant connexe et tout atome \'etant intersection d\'ecroissante de
$R^n_i$, les atomes sont connexes.
\par
Soit $\nu:=\lim \mu^x_{n_i}$. Il est clair que $\nu(\varphi)=0$ si
$\Lambda^n\varphi=0$ pour un $n\geq 1$. 
On consid\`ere pour $i\not=j$ 
la fonction $\varphi$ d\'efinie par $\varphi:=1$ sur
$R^n_j$, $\varphi:=-1$
sur $R^n_i$ et $\varphi=0$ ailleurs. On a que 
$\nu(R^n_i)=\nu(R^n_j)$ et $\mu(R^n_i)=\mu(R^n_j)$. Donc 
$\nu(R^n_i)=d_t^{-n}=\mu(R^n_i)$. 
On en d\'eduit que $\nu$ et $\mu$ s'annulent sur les
atomes de $\B$. De plus, $\mu$ et $\nu$ co\"{\i}ncident sur $\B$. Si
les atomes de $\B$ sont triviaux sauf pour un nombre  d\'enombrable
au plus, il en r\'esulte que $\mu=\nu$. L'hypoth\`ese sur $\K$
assure pr\'ecis\'ement que l'ensemble des atomes non r\'eduits \`a un
point est au plus d\'enombrable.
\par
Pour montrer que $\mu$ est K-m\'elangeante, il suffit
d'observer que pour toute fonction continue $\varphi$, on a
$\Lambda^n\varphi\rightarrow c_\varphi$ dans $\Ltwo(\mu)$. En effet,
$\Lambda^n\varphi\rightarrow c_\varphi$ ponctuellement; on peut
appliquer le th\'eor\`eme de convergence domin\'ee.
\subsection{Mesure d'\'equilibre}
Nous allons maintenant construire une mesure d'\'equilibre dans le
cadre riemannien. Soit $V$ une vari\'et\'e
riemannienne munie d'une forme volume
$\Omega$. Soit $f$ une application de classe
${\cal C}^1$ d'un ouvert
$U$ de $V$
dans $V$. On notera $\Cr$ l'ensemble critique de $f$. 
On suppose que $\#f^{-1}(x)=d_t>1$ pour tout
$x\in V\setminus f(\Cr)$.
Soit $X$ un compact de $U$ de mesure positive tel que
$Y:=f^{-1}(X)\subset X$. Posons 
$$\sigma_N:=\frac{1}{N}\sum_{n=1}^N\frac{(f^n)^*\Omega_{|X}}{d_t^n}$$
o\`u on a not\'e $\Omega_{|X}$ la restriction de $\Omega$ \`a $X$. On
peut supposer $\int\Omega_{|X}=1$. La forme $\Omega$ \'etant de
degr\'e maximal, il existe une fonction
continue $J$ 
v\'erifiant
$f^*\Omega=J\Omega$. C'est le jacobien de $f$. On suppose
$J\geq 0$ sur $X$, \cad que $f$ pr\'eserve l'orientation sur
$X\setminus \Cr$. 
\par 
Observons que
$d_t^{-n}(f^n)^*\Omega_{|X}$ d\'efinit une mesure 
bien que l'op\'erateur $(f^n)^*$ ne soit pas
\`a priori d\'efini sur l'ensemble des mesures. 
En g\'en\'eral, si $\nu$ est une mesure de $V\setminus f(\Cr)$ on peut
d\'efinir la mesure $f^*\nu$ sur $U\setminus f^{-1}(f(\Cr))$. Cet
op\'erateur est continu. 
D'apr\`es le
th\'eor\`eme de Sard, l'ensemble des valeurs critiques de $f^n$ est de
mesure nulle. On en d\'eduit que $d_t^{-n}(f^n)^*\Omega_{|X}$
est une mesure de probabilit\'e. 
\par
Rappelons ici la notion d'entropie d'une mesure invariante $\sigma$ de
$f$. D\'esignons par $d$ la distance sur $V$. 
On pose
$$d_n^f(x,y):=\max_{0\leq i\leq n-1} d(f^i(x), f^i(y)).$$
Notons $B^f_n(x,r)$ la boule de centre $x$ et de rayon $r$
pour la distance $d^f_n$. D'apr\`es Brin-Katok \cite{BrinKatok}, la
mesure $\sigma$ \'etant invariante, pour $\sigma$-presque tout
$x$ la limite suivante existe
$$\h_\sigma(f,x):=\lim_{\delta\rightarrow 0} \limsup_{n\rightarrow
\infty}\frac{-\log \sigma(B^f_n(x,\delta))}{n}$$
et d\'efinit une fonction invariante par $f$ telle que
$$\int \h_\sigma(f,x) \d\sigma =\h_\sigma(f).$$
C'est {\it l'entropie} 
de $\sigma$ pour $f$. 
\par
Dans \cite{Mane}, Ma\~ne montre que pour les fractions rationnelles de
$\P^1$ la fonction $\log J$ est int\'egrable pour les mesures
d'entropie positive. C'est ce r\'esultat que nous a inspir\'e.
\begin{theoreme} Avec les notations ci-dessus, 
soit $\sigma$ une valeur d'adh\'erence de la suite
  $(\sigma_N)$. Alors
\begin{enumerate}
\item $\sigma$ est une mesure invariante: $f_*\sigma=\sigma$.
\item $\int \log J \d\sigma \geq \log d_t$.
\item $f^*\sigma=d_t\sigma$ sur $V\setminus f^{-1}(f(\Cr))$. On a
  $f^*\sigma=d_t\sigma$ lorsque $f$ d\'efinit un rev\^etement
  ramifi\'e de $f^{-1}(X)$ au dessus de $X$.
\item L'entropie $\h_\sigma$ de la mesure
$\sigma$ est sup\'erieure ou \'egale
  \`a $\log d_t$.
\end{enumerate}
\end{theoreme}
\begin{preuve} 
1. Observons que 
$$f_*\sigma_N-\sigma_N=
\frac{1}{N}\left(\Omega_{|X}-\frac{(f^N)^*\Omega_{|X}}{d_t^N}\right)$$
et que ces mesures tendent vers $0$ quand
$N\rightarrow\infty$. L'op\'erateur $f_*$ \'etant continu, on a
$f_*\sigma=\sigma$. 
\par
2. 
Soit $M>0$ tel que $J\leq M$ sur $X$. 
Fixons $m>0$, $\delta>0$ et posons
$$g_m(x):=\min\left(\log\frac{M}{J(x)}, m+\log M\right)=
\min\left(\log\frac{M}{J(x)}, m'\right)$$
o\`u $m':=m+\log M$. 
C'est une famille de fonctions continues, positives,
born\'ees sur $X$ 
qui tend vers $\log M/J(x)$ quand $m$ tend vers
l'infini.  
Posons 
$$\mu_n:=\frac{(f^n)^*\Omega_{|X}}{d_t^n}\ \  \mbox{ et }\ \ 
s_N(x):=\frac{1}{N}\sum_{q=0}^{N-1} g_m(f^q(x)).$$ 
On a
\begin{eqnarray}
\int s_N\d\mu_N & = & \frac{1}{N}\sum_{q=0}^{N-1} \int
g_m(f^q(x))\frac{(f^N)^*\Omega_{|X}}{d_t^N}\nonumber \\
& = & \frac{1}{N}\sum_{q=0}^{N-1} \int
g_m(x)\frac{(f^{N-q})^*\Omega_{|X}}{d_t^{N-q}}\nonumber \\
& = & \frac{1}{N}\sum_{q=0}^{N-1} \int
g_m(x) \d\mu_{N-q} =\int g_m \d\sigma_N.
\end{eqnarray}
Pour minorer $\int \log J \d\sigma$, il suffit donc de
majorer $\int s_N \d\mu_N$.
\par
Pour $\alpha>0$, posons $X^\alpha_N:=\{x, \
s_N(x)>\alpha\}$. Puisque $s_N(x)\leq m'$, on a 
\begin{eqnarray*}
\int g_m\d\sigma_N=\int s_N \d\mu_N & \leq &
m'\mu_N(X^\alpha_N)+\alpha(1-\mu_N(X^\alpha_N))\\
& = & \alpha+(m'-\alpha)\mu_N(X^\alpha_N).
\end{eqnarray*}
Si $\mu_N(X^\alpha_N)$ tend vers 0 quand $N\rightarrow\infty$, on a 
$$\limsup_{N\rightarrow\infty}\int g_m
\d\sigma_N\leq \alpha
\ \ \mbox{ et par suite, }\ \  \int\log\frac{M}{J}
\d\sigma \leq \alpha.$$
Il s'agit donc de d\'eterminer la borne inf\'erieure des $\alpha$ tels
que $\mu_N(X^\alpha_N)$ tende vers 0 pour $m$ fix\'e. Par
d\'efinition de $\mu_N$, on a
\begin{eqnarray*}
\mu_N(X^\alpha_N) & = & \int_{X^\alpha_N}
\frac{\prod_{q=0}^{N-1}J\circ f^q}{d_t^N} \Omega_{|X}. 
\end{eqnarray*}
Posons pour $\delta>0$ fix\'e et $j\in \Z$,
$$W_j:=\big\{\exp(-j\delta)<J\leq
\exp(-(j-1)\delta)\big\}$$ et
$$\tau_j(x):=\frac{1}{N} \# \big\{q, \ 
f^q(x)\in W_j \mbox{ et } 0\leq q\leq N-1\big\}.$$
On a $\sum\tau_j=1$ et
\begin{eqnarray}
\mu_N(X^\alpha_N) & \leq & 
\int_{X^\alpha_N}\left[\frac{1}{d_t}\exp\left(\sum -(j-1)\delta
\tau_j\right)\right]^N \Omega_{|X}.
\end{eqnarray}
Or sur $X^\alpha_N$ on a 
$$\alpha < s_N< \sum \tau_j (\log M + j\delta) =\sum
j\delta\tau_j +\log M.$$
Donc 
$$ -\sum(j-1)\delta\tau_j <-\alpha +(\log M +\delta)$$
et d'apr\`es (2)
$$\mu_N(X^\alpha_N)\leq \int_{X^\alpha_N}\left[\frac{\exp(-\alpha) M
\exp(\delta)}{d_t}\right]^N \Omega_{|X}.$$
Pour tout $\alpha>\log(M/d_t)+\delta$, on a
$\mu_N(X^\alpha_N)\rightarrow 0$. Or $\delta$ est arbitrairement
petit. On a $\int g_m\d \sigma_N\leq \log (M/d_t)$. 
Comme $g_m$ est continue, $\int g_m \d\sigma \leq \log (M/d_t)$
et par suite,
$\int\log J \d\sigma\geq \log d_t$.
\par
3. Sur $V\setminus f^{-1}(f(\Cr))$, on a
$$f^*\sigma_N-d_t\sigma_N=
\frac{1}{N}\left[\frac{(f^{N+1})^*\Omega_{|X}}{d_t^N}
-f^*\Omega_{|X}\right].$$
Donc ces mesures tendent vers $0$ quand
$N\rightarrow\infty$. Sur $V\setminus f^{-1}(f(\Cr))$, l'op\'erateur
$f^*$ \'etant 
continu, on a $f^*\sigma=d_t\sigma$.
Si $f$ d\'efinit un rev\^etement ramifi\'e
de $f^{-1}(X)$ au dessus de $X$, l'op\'erateur $f^*$ est continu sur
l'ensembles des mesures \`a support dans $X$. Dans ce cas, on a 
$f^*\sigma =d_t\sigma$ sur $V$. 
\par
4. La mesure $\sigma$ est de jacobien constant $d_t$, \ie pour tout
bor\'elien $B$ sur lequel $f$ est injective, on a $1_{f(B)}=f_*1_B$,
par suite,
\begin{eqnarray}
\sigma(f(B)) & = & d_t\sigma(B)
\end{eqnarray}
Notons que d'apr\`es la partie 2, l'ensemble
$\Cr$ est de mesure nulle pour $\sigma$, 
donc le jacobien
$J_\sigma$ de
$f$ par rapport \`a la mesure $\sigma$ est presque partout \'egal
\`a $d_t$.
Une formule de Parry \cite{Parry} assure que
$$\h_\sigma(f)\geq \int\log J_\sigma \d\sigma=\log d_t.$$ 
\end{preuve}
\begin{remarque} \rm 
Soit $f:Y\longrightarrow X$ un rev\^etement
  ramifi\'e de degr\'e $d_t\geq 2$ avec
$Y\subset\subset X$. Soit $\nu$ une
  mesure de probabilit\'e sur $X$  telle que $f^*\nu=J\nu$ avec une
  fonction $J\geq
  0$ semi-continue sup\'erieurement sur $X$. 
La d\'emonstration du th\'eor\`eme 2.3.1 montre que $\int \log
  J\d\sigma\geq \log d_t$ pour toute valeur
d'adh\'erence $\sigma$ de la suite de
  $\frac{1}{N}\sum_1^N d_t^{-n}(f^n)^*\nu$. 
\end{remarque}
\begin{remarque} \rm Supposons que la restriction de $f$ \`a
  $Y:=f^{-1}(X)$ soit un rev\^etement ramifi\'e de degr\'e $d_t\geq
  2$ au dessus de $X$.
Il est facile de v\'erifier que les points
extr\'emaux de $\M$ sont des mesures ergodiques. 
En utilisant la
d\'ecomposition de Choquet de $\sigma$ relativement aux
\'el\'ements extr\'emaux de $\M$, on a
$$\sigma=\int\nu_\alpha \d\tau(\alpha)$$
o\`u $\nu_\alpha$ est extr\'emale et $\tau$ est une mesure de
probabilit\'e sur l'ensemble des mesures extr\'emales.
Il en r\'esulte qu'il existe une mesure extr\'emale $\nu$, donc
ergodique,  pour laquelle
\begin{eqnarray}
\int \log J \d\nu & \geq &  \log d_t
\end{eqnarray}
Par suite, $\h_\nu(f)\geq \log d_t$.
Observons que l'ensemble $\M'$ des mesures de
probabilit\'e $\nu\in\M$ v\'erifiant (4) est un convexe
compact. En prenant des vari\'et\'es produit, il est facile de
construire des exemples montrant qu'une telle mesure n'est pas unique.  
\par
Rappelons ici la notions d'entropie topoplogique 
\cite{KatokHasselblatt}. 
On dit qu'un ensemble $F$ est
{\it $(n,\epsilon)$-s\'epar\'e}
si pour $x, y\in F$ distincts, on a $d^f_n(x,y)\geq \epsilon$.  
C'est-\`a-dire qu'\`a l'echelle $\epsilon$, les orbites de $x$ et $y$
peuvent \^etre distingu\'ees avant le temps $n$. 
On d\'efinit {\it l'entropie topologique} 
de $f$ par
$$\h(f):=\sup_{\epsilon>0}\limsup_{n\rightarrow\infty}
\frac{1}{n}\log \max \big\{ \# F \mbox{ pour } F\subset X\ (n,\epsilon)
\mbox{-s\'epar\'e}\big\}.$$
{\it Le principe variationnel} 
affirme que
$$\h(f)=\sup\big\{\h_\nu(f)\mbox{ pour }\nu \mbox{ invariante
  ergodique}\big\}.$$
\par
On retrouve dans le cadre du th\'eor\`eme 2.3.1, o\`u on a suppos\'e
$J\geq 0$, le r\'esultat de
Misiurewicz-Przytycki \cite{KatokHasselblatt}
qui dit que l'entropie topologique $\h_t(f)$ de $f$
est minor\'ee par $\log d_t$. 
\end{remarque}
\begin{proposition} Soit $f:U\longrightarrow V$ une application de
  classe ${\cal C}^1$. Soit $X$ un compact de $U$ tel que
  $f^{-1}(X)\subset X$. Supposons que la restriction de $f$ \`a
  $f^{-1}(X)$ d\'efinisse un rev\^etement ramifi\'e de degr\'e
  $d_t\geq 2$
au dessus de $X$. Si  $\mu\in\M'$ 
est une mesure ergodique, alors
\begin{enumerate}
\item $J_k:=\supp(\mu)$ est  parfait. 
\item Les exposants de Lyapounov $(\lambda_1,\ldots,
  \lambda_k)$ relatifs \`a $\mu$ sont constants et v\'erifient
  $\lambda_1+\cdots +\lambda_k\geq \log d_t$.
\end{enumerate}
En particulier, si $\dim_\R V=1$, la mesure $\mu$ est
  hyperbolique. 
\end{proposition}
\begin{preuve}
1. Si $x$ est un point isol\'e de $\supp(\mu)$ alors
$\mu\{x\}>0$. L'ergodicit\'e de $\mu$ entra\^{\i}ne que
$\mu=\delta_x$ et en particulier $f^{-1}(x)=\{x\}$. Ce qui
entra\^{\i}ne que $x$ est critique. Cela contredit l'int\'egrabilit\'e
de $\log J$ par rapport \`a $\mu$.
\par
2. Pour montrer l'existence des
exposants de Lyapounov pour $\mu$
il suffit \cite{KatokHasselblatt} de v\'erifier les in\'egalit\'es
$$\int \log^+\|Df\| \d\mu<+\infty \ \ \mbox{ et }\ \ 
\int \log^+\|(Df)^{-1}\|\d\mu<+\infty.$$
La premi\`ere in\'egalit\'e est \'evidente, 
l'int\'egrabilit\'e de $\log J$
par rapport \`a $\mu$ suffit pour v\'erifier la seconde.
Il existe donc 
$\lambda_{(1)}>\lambda_{(2)}>\cdots > \lambda_{(m)}$ et une
d\'ecomposition $\mu$-mesurable du fibr\'e tangent \`a $V$ en fibr\'es
invariants 
$E_1\oplus\cdots \oplus E_m$ tels que $\lim
\frac{1}{n}\log\|Df^n(x)v\|=\lambda_{(j)}(x)$ $\mu$-presque partout
si $|v|=1$ et $v\in E_j$. De
plus, les $\lambda_{(j)}$ sont constants car
$\mu$ est ergodique. Ce sont
{\it les exposants de Lyapounov}
de {\it multiplicit\'e}
respective $\dim E_j$. Un
changement de variable classique \cite[p.666]{KatokHasselblatt} 
permet de montrer que 
$$\lim\frac{1}{n}\log J_{f^n}=\sum(\dim E_j)
\lambda_{(j)}=\sum\lambda_j$$ 
o\`u $J_{f^n}$ est le jacobien
r\'eel de $f^n$ et $\lambda_1,\ldots, \lambda_k$
sont les $k$ exposants de
Lyapounov \'ecrits en tenant compte
les multiplicit\'es. En appliquant le
th\'eor\`eme ergodique \`a la fonction $\mu$-int\'egrable $\log J$,
on trouve
$$\sum \lambda_j=\lim \frac{1}{n}
\log J_{f^n}=\lim\frac{1}{n}\sum_{q=0}^{n-1} \log J\circ
f^q=\int\log J \d\mu\geq \log d_t. $$
\end{preuve}
\section{Applications d'allure polynomiale}
Dans ce paragraphe, nous \'etudions la dynamique d'une grande classe
d'applications holomorphes: les applications d'allure
polynomiale. Nous construisons une mesure d'\'equilibre $\mu$ 
d'entropie maximale. Nous montrons qu'elle est K-m\'elangeante. Nous
\'etudions les ensembles exceptionnels,
les points p\'eriodiques et les exposants de Lyapounov de cette
mesure. Nous montrons, en particulier, que si $\mu$ est PLB (\cad que
les fonctions p.s.h. sont $\mu$-int\'egrables), l'ensemble
exceptionnel est analytique, les points p\'eriodiques r\'epulsifs sont
denses dans $\supp(\mu)$, les exposants de Lyapounov sont
strictement positifs et la vitesse de m\'elange est d'ordre
exponentiel. La propri\'et\'e pour la mesure $\mu$ associ\'ee \`a $f$,
d'\^etre PLB, est stable par
pertubation. Cela nous permet d'exhiber une vaste famille 
d'applications avec une telle mesure d'\'equilibre. Nous avons
rassembl\'e au paragraphe 3.9 les propri\'et\'es non dynamiques des
mesures PLB. Le lecteur peut d'abord se familiariser avec ces
propri\'et\'es. 
\subsection{Quelques d\'efinitions}
Introduisons quelques notions qui d\'efinissent le cadre de notre
\'etude.
\begin{definition} \rm
Une vari\'et\'e complexe $V$ est dite
{\it S-convexe} 
si elle poss\`ede au moins une fonction
continue strictement p.s.h. Elle est dite {\it S-s\'epar\'ee} 
si pour tout point $a$ de $V$ et tout ensemble fini $A\subset
V\setminus\{a\}$, 
il existe une fonction
p.s.h. continue $\varphi$ v\'erifiant $\varphi(a)>\sup_A\varphi(z)$.
\end{definition}
\par
Il est clair qu'une vari\'et\'e S-convexe ne peut contenir des
ensembles analytiques compacts de dimension strictement positive.
Les ouverts relativement compacts d'une telle vari\'et\'e sont Kobayashi
hyperboliques \cite{Sibony1}. 
Tout ouvert d'une vari\'et\'e de Stein, par exemple $\C^k$, 
est S-convexe et S-s\'epar\'e. Si $V$ est
S-convexe, il r\'esulte du th\'eor\`eme de Richberg \cite{Richberg}, 
qu'elle poss\`ede une fonction
strictement p.s.h. $\Phi$ de classe
${\cal C}^\infty$ et donc une m\'etrique k\"ahl\'erienne associ\'ee
\`a la forme $\omega:=\ddc\Phi$.
\begin{definition} \rm
On appelle {\it application (holomorphe)
d'allure polynomiale}
toute application holomorphe propre $f$ de
$U$ sur $V$ o\`u $V$ est une vari\'et\'e complexe connexe,
S-convexe et $U$ est un
ouvert relativement compact de $V$ (\'eventuellement non-connexe).
On d\'efinit {\it l'ensemble de Julia rempli}
de $f$ par
$\K:=\bigcap_{n\geq 0} U_{-n}$
o\`u $U_{-n}:=f^{-n}(V)$. 
On appelle
{\it degr\'e topologique}
de $f$ le nombre $d_t$ de
pr\'eimages d'un point $z\in V$ compt\'ees avec mulitiplicit\'es.
Dans la plupart des r\'esultats, la connexit\'e de $V$ n'est pas
n\'ecessaire. On a seulement besoin que $V$ ait un nombre fini de
composantes et que chaque point de $V$ admet exactement $d_t$
pr\'eimages compt\'ees avec multiplicit\'es.
\end{definition}
\par
L'application $f$ est ouverte et d\'efinit un rev\^etement
ramifi\'e de $U$ au dessus de $V$. En particulier, le degr\'e
topologique $d_t$ de $f$ ne d\'epend pas du point $z\in V$.
L'ensemble
$\K$ est le plus grand compact {\it totalement invariant}
par $f$ au sens o\`u
$f^{-1}(\K)=\K$. 
\par
Les applications d'allure polynomiale en
dimension 1 lorsque les ouverts $U$ et $V$ sont simplement connexes, 
ont \'et\'e consid\'er\'ees par Douady-Hubbard \cite{DouadyHubbard}. 
Ils en ram\`enent
l'\'etude dynamique \`a celle des polyn\^omes. Cette r\'eduction est
obtenue gr\^ace au th\'eor\`eme de Riemann-mesurable. Ce th\'eor\`eme
n'ayant pas d'analogue \`a plusieurs variables, nous adoptons une
approche diff\'erente dont l'ingr\'edient essentiel est un
th\'eor\`eme de convergence pour les fonctions plurisousharmoniques.
\par
Notons ici que les endomorphismes polynomiaux de $\C^k$ dont l'infini
est ``attirant'' sont \`a allure polynomiale. Plus pr\'ecis\'ement,
la restriction d'une telle application \`a un ouvert convenable
est une application \`a allure polynomiale. En g\'en\'eral, un
endomorphisme  polynomial de $\C^k$ avec $k\geq 2$ n'est pas \`a
allure polynomiale m\^eme s'il est propre. 
La terminologie ``allure polynomiale'' est en fait
utilis\'ee par Douady-Hubbard dans le cas de dimension $1$, nous
l'avons conserv\'ee.   
%
%\par
%Dans tout ce paragraphe, on note $\Phi$ une fonction strictement
%p.s.h. de classe ${\cal C}^\infty$ de $V$ et $\omega:=\ddc \Phi$
%la forme k\"ahl\'erienne associ\'ee. Il existe alors des fonctions
%strictement p.s.h. ${\cal C}^\infty$ et positives sur $U$. Nous 
%supposons donc que $\Phi$ est positive sur $U$. 
\par
Rappelons que dans le cas d'une application holomorphe $g:
\P^k \longrightarrow \P^k$, on d\'efinit pour tout
$1\leq l\leq k$ {\it le degr\'e dynamique d'ordre $l$}
de $g$, not\'e
$d_l$, comme le degr\'e de la vari\'et\'e $g^{-1}(H)$ o\`u $H$ est
un sous-espace projectif g\'en\'erique de dimension $k-l$ de
$\P^k$. Le degr\'e topologique de $g$ est \'egal \`a $d_k$.
Si $\omega_0$ d\'esigne la forme de Fubini-Study de $\P^k$, on a
$$d_l=\int g^*(\omega_0^l)\wedge \omega_0^{k-l}=\int \omega_0^l
\wedge g_*(\omega_0^{k-l}).$$ 
Ces degr\'es jouent le r\^ole crucial pour calculer l'entropie de
l'application.
\par
Lorsque $f$ est une application d'allure polynomiale, nous
consid\'erons les notions locales analogues. Pour tout
$1\leq l\leq k:=\dim V$ et tout $n\geq 1$, on pose
$$d_{l,n}:=\int_{U_{-n-1}} (f^n)^*(\omega^l) \wedge \omega^{k-l} =
\int_U (f^n)_*(\omega^{k-l})\wedge \omega^l.$$
On v\'erifie que $d_{l,n}$ est r\'eel positif. 
\begin{definition} \rm On appelle 
{\it degr\'e dynamique d'ordre $l$}
de $f$ le nombre r\'eel
positif 
$$d_l:=\limsup_{n\rightarrow \infty} \sqrt[n]{d_{l,n}}.$$ 
\end{definition}
On a $d_{k,n}=d_t^n\int_U\omega^k$ et
donc $d_k=d_t$.
Notons que, si $f$ est la restriction d'une application polynomiale,
les degr\'es dynamiques d\'efinis ci-dessus sont en g\'en\'eral plus
petits que ceux d\'efinis globalement (\voir l'exemple 3.10.6).
La proposition
suivante montre que les $d_l$ sont bien 
d\'efinis du point de vue dynamique.
\begin{proposition} 1. Les degr\'es dynamiques 
$d_l$ ne d\'ependent pas de la forme
k\"ahl\'erienne $\omega$ choisie.
\par
2. Soit $U'\subset\subset V$ un voisinage de $\K$. Alors
$$d_l=\limsup_{n\rightarrow \infty} \left(\int_{U'}
(f^n)_*(\omega^{k-l}) \wedge \omega^l\right)^{1/n}.$$
\par
3. Les $d_l$ sont des nombres r\'eels strictement positifs, invariants
par conjugaison holomorphe. 
\end{proposition}
\begin{preuve} 1. Soit $\tilde\omega$ une autre forme k\"ahl\'erienne
sur $V$. Il existe une constante $c\geq 1$ telle que
$c^{-1}\omega\leq \tilde\omega \leq c\omega$ sur $U$. Posons
$$d'_{l,n}:=\int_{U}
(f^n)_*(\tilde\omega^{k-l}) \wedge \tilde\omega^l\ \ \mbox{ et }\ \ 
d_l':=\limsup \sqrt[n]{d'_{l,n}}.$$ 
On a 
$$c^{-k} d_{l,n}\leq d'_{l,n} \leq
c^k d_{l,n}.$$ 
En cons\'equence, on obtient $d'_l=d_l$. Ce qui montre aussi
l'invariance par conjugaison holomorphe.
\par
2. Fixons $m$ assez grand tel que $U_{-m}$ soit
contenu dans $U'$. Posons
$$d^*_{l,n}:=\int_{U_{-m}}
(f^n)_*(\omega^{k-l}) \wedge \omega^l\ \  \mbox{ et } \ \ d^*_l:=\limsup
\sqrt[n]{d^*_{l,n}}.$$
Il est clair que $d^*_{l,n}\leq d_{l,n}$ et donc $d_l^*\leq d_l$.
\par
Puisque la forme $(f^{m-1})^*(\omega^{k-l})$ est lisse sur
$U_{-m+1}$, il existe une constante $c>0$ telle que sur
$U_{-m}$ on ait  $(f^{m-1})^*(\omega^l) \leq c\omega^l$.
On en d\'eduit que pour tout $n\geq m$
\begin{eqnarray*}
d_{l,n} & = & \int_{U_{-m}} (f^{n-m+1})_*(\omega^{k-l}) \wedge
(f^{m-1})^*(\omega^l) \\
& \leq &  c \int_{U_{-m}} (f^{n-m+1})_*(\omega^{k-l}) \wedge\omega^l
=c d^*_{l,n-m+1}
\end{eqnarray*}
Par cons\'equent, $d_l\leq d_l^*$ et donc $d_l=d_l^*$.
\par
De la m\^eme mani\`ere, on montre que
$$\limsup_{n\rightarrow \infty} \left(\int_{U'}
(f^n)_*(\omega^{k-l}) \wedge \omega^l\right)^{1/n} = d_l^*.$$
\par
3. On a $f^*\omega\leq c\omega$ sur $U_{-1}$ pour un $c>0$, 
donc $(f^n)^*\omega \leq
c^n \omega$ sur $U_{-n-1}$. Par cons\'equent,
$$d_{l,n}  =  \int_{U_{-n-1}} (f^n)^*(\omega^l) \wedge
\omega^{k-l}\leq c^{ln}\int_{U_{-n-1}} \omega^k \leq
c^{ln}\int_U \omega^k$$
et
$$d_{l,n} \geq c^{-(k-l)n}\int_{U_{-n-1}} (f^n)^*(\omega^l)\wedge
(f^n)^* (\omega^{k-l})=c^{-(k-l)n}d_t^n\int_U \omega^k.$$
On en d\'eduit que $c^{-k+l}d_t\leq d_l\leq c^l$.
\end{preuve}
\par
L'ensemble critique $\Cr$ de $f$ a une grande influence \`a la
dynamique de $f$. Nous aurons besoin de mesurer 
la croissance des images de $\Cr$. 
Notons 
\begin{enumerate}
\item[] $\Cr_{-n}:=\bigcup_{j=0}^{n-1} f^{-j}(\Cr)$
l'ensemble critique de $f^n$; 
\item[] $\Cr_n:=f^n(\Cr\cap U_{-n})$;
\item[] $\PC_n:=\bigcup_{j=1}^n \Cr_j$ 
{\it l'ensemble postcritique
d'ordre $n$};
\item[] $\PC_\infty:=\bigcup_{j=1}^\infty \Cr_j$ 
{\it l'ensemble postcritique d'ordre infini}
de $f$. 
\end{enumerate}
Le volume de
$\Cr_n$ dans $U$ est \'egal \`a $d_t^n\delta_n$ o\`u 
$$\delta_n:=\frac{1}{d_t^n}\int_{\Cr_n\cap U}\omega^{k-1}=
\frac{1}{d_t^n}\int_{\Cr\cap
U_{-n-1}} (f^n)^*(\omega^{k-1}).$$
Soit $J$  le jacobien r\'eel de $f$.
Si $V$ est un ouvert de $\C^k$ et $\omega=\ddc |z|^2$, on a $\ddc\log
J=2[\Cr]$. Dans le cas g\'en\'eral, il existe une constante $c\geq
0$ telle que sur $U$ on ait 
$\ddc\log J\geq 2[\Cr] -2c\omega$. Par cons\'equent,
$$\delta_n\leq \frac{1}{2}\int_U \ddc \left[ \frac{(f^n)_* \log
    J}{d_t^n} \right] \wedge \omega^{k-1} +\frac{c d_{k-1,n}}{d_t^n}.$$
Rappelons que $\d^c=i(\overline\partial -\partial)$ et que $[\Cr]$
    d\'esigne le courant d'int\'egration associ\'e \`a $\Cr$ (\voir
    \cite{Lelong}, \cite{GriffithsHarris}). On peut introduire un
    exposant $\delta$ d\'ecrivant la croissance relative du volume de
    $f^n(\Cr\cap U_{-n})$ par rapport au degr\'e dynamique. Posons
    $$\delta:=\limsup_{n\rightarrow \infty}\sqrt[n]{\delta_n}.$$  
\subsection{Mesure d'\'equilibre}
Dans le cadre des applications d'allure polynomiale,
on peut am\'eliorer 
le th\'eor\`eme 2.3.1 de mani\`ere suivante:
\begin{theoreme} Soit $f:U\longrightarrow V$ une application
d'allure polynomiale de degr\'e topologique $d_t\geq 2$. 
Il existe une mesure de probabilit\'e $\mu$ port\'ee par
$\partial \K$ v\'erifiant la relation
$f^*\mu=d_t\mu$ et telle que pour toute forme volume $\Omega$, de masse
$1$ dans $\Ltwo(V)$, la suite de mesures
$d_t^{-n}(f^n)^*\Omega$ converge vers $\mu$. 
Pour toute $v\in \Ltwo(\mu)$ la suite $v_n:=d_t^{-n}(f^n)_* v$ 
converge dans
$\Ltwo(\mu)$ vers la constante $c_v:=\int v\d\mu$. L'application
$f$ est K-m\'elangeante pour la mesure $\mu$.
\end{theoreme}
\par
On dira que $\mu$ est {\it la mesure
    d'\'equilibre} de $f$, 
son support $J_k:=\supp(\mu)$ est {\it
l'ensemble de Julia (d'ordre maximal)}.
On a not\'e $k$ la dimension
de $V$.
\begin{lemme} Soit $\varphi$ une fonction 
p.s.h. born\'ee 
au voisinage de $\K$. Alors pour $n$ assez
grand $\varphi_n:=d_t^{-n}(f^n)_*\varphi$
est une fonction
p.s.h., born\'ee
sur $V$ et ces fonctions convergent dans $\Lploc(V)$
vers une constante $c_\varphi$ quand $n\rightarrow+\infty$ pour tout
$p\geq 1$.
De plus, l'ensemble 
$$\E'_\varphi:=\big\{z\in
V,\ \limsup_{n\rightarrow +\infty}
\varphi_n(z)<c_\varphi\big\}$$
est pluripolaire et invariant dans le sens o\`u 
$f(\E'_\varphi\cap U)\subset\E'_\varphi$. (On
verra dans la proposition 3.2.5 
que si $\varphi$ est born\'ee et strictement p.s.h. au voisinage
de $\K$, $\E'_\varphi$ ne d\'epend pas de $\varphi$).
\end{lemme}
\begin{preuve} Soit $\K_\delta$ un voisinage de $\K$ dans
lequel $\varphi$ est p.s.h., born\'ee par une constante $c>0$.
La fonction
$$\varphi_n(z):=\Lambda^n\varphi=\frac{(f^n)_*\varphi}{d_t^n}
=\frac{1}{d_t^n}\sum_{f^n(w)=z}\varphi(w)$$ 
est p.s.h. et born\'ee par $c$
dans $\K_n:=\{z\in V,\ 
f^{-n}(z)\subset \K_\delta\}$. On a
$\K_n=V$ pour $n$ assez grand.
\par
Montrons que $\varphi_n$ tend faiblement vers une
constante $c_\varphi$. 
Observons que si $\Psi$ est p.s.h. dans $V$ telle que  
$f_*\Psi\geq d_t\Psi$, alors $\Psi$ est constante. Ceci est
une cons\'equence du principe du maximum et de
l'in\'egalit\'e 
$$\sup_U\Psi\geq \sup_V\frac{f_*\Psi}{d_t}\geq \sup_V\Psi.$$
Posons $\Psi_0:=\limsup \varphi_n$ et $\Psi$ la r\'egularis\'ee
de $\Psi_0$. Alors $\Psi$ est une fonction p.s.h. On a pour
tout $z\in V$
\begin{eqnarray*}
d_t \Psi_0(z) & = & \limsup d_t \varphi_{n+1} (z) =\limsup f_*\varphi_n
(z)\\
& \leq & f_*\limsup \varphi_n(z) =f_*\Psi_0(z).
\end{eqnarray*}
On en d\'eduit que $f_*\Psi\geq d_t\Psi$ et donc $\Psi$ est
constante. Posons 
$$c_\varphi:=\Psi=(\limsup\varphi_n)^*.$$
\par
On peut extraire de
la suite $(\varphi_n)$ des sous-suites convergentes dans $\Lploc(V)$.
Supposons que $\varphi_{n_j}$ tend vers une fonction
p.s.h. $\psi$. Montrons que $\psi=c_\varphi$.
Sinon il existe
$\epsilon>0$ tel que $\psi\leq c_\varphi-2\epsilon$ sur $U$. D'apr\`es
le lemme de Hartogs \cite{Hormander}, $\varphi_{n_j}\leq c_\varphi-\epsilon$
sur $U_{-2}$ pour $j$ assez grand.
Or $\varphi_{n_j+l}\leq c_\varphi-\epsilon$ pour tout
$l\geq 1$ car
$\varphi_{n_j+l}=d_t^{-l}(f^l)_*\varphi_{n_j}$ et
$f^{-l}(U)\subset U_{-2}$. Cela contredit le fait que
$(\limsup\varphi_n)^*=c_\varphi$.
\par
Il existe un ensemble pluripolaire
$\E'_\varphi$ tel que $\limsup\varphi_n(z)=c_\varphi$ 
pour tout $z\in V\setminus \E'_\varphi$ \cite{Lelong, BedfordTaylor}. On a
$\E'_\varphi=\{z\in V,\ \limsup
\varphi_n(z)<c_\varphi\}$.
Montrons que $\E'_\varphi$ est invariant.
Fixons un point $z\in V$. On a
$$d_t\limsup \varphi_n(z)=
\limsup \sum_{f(w)=z}\varphi_{n-1}(w)\leq \sum_{f(w)=z}\limsup
\varphi_{n-1}(w).$$
D'autre part, $\limsup \varphi_{n-1}(w)\leq c_\varphi$ pour tout
$w$. Soit $z\not \in\E'_\varphi$. On a $\limsup \varphi_n(z)=c_\varphi$.
Par cons\'equent, $\limsup\varphi_{n-1}(w)=c_\varphi$ pour tout
$w\in f^{-1}(z)$. On en d\'eduit que $f^{-1}(z)\cap \E'_\varphi=\emptyset$
et donc $f(\E'_\varphi\cap U)\subset \E'_\varphi$. 
\end{preuve}
\par
Le raffinement suivant du lemme de Hartogs nous sera utile.
\begin{proposition} Soit $(v_j)$ une suite de fonctions
  p.s.h born\'ees dans $V$ convergeant dans $\Loneloc(V)$ vers une
  fonction continue $v$. Soit $\sigma$ une mesure de probabilit\'e \`a
  support compact dans $V$ telle que 
$\int v_n\d \sigma \rightarrow \int
  v\d\sigma$. Alors pour tout $1\leq p<\infty$, $(v_n)$ converge vers
  $v$ dans $\mbox{\rm L}^p(\sigma)$ et $\limsup v_n=v$ $\sigma$-presque
  partout. 
\end{proposition}
\begin{preuve} Pour $\epsilon>0$, 
d'apr\`es le lemme de Hartogs, on a pour $n$ assez
  grand $v_n\leq v+\epsilon$ sur $\supp (\sigma)$. 
Pour $\delta>0$ posons
$A_n^\delta:=\{v_n<v-\delta\}$. On a
$$\int_V v_n\d\sigma\leq \int_{A_n^\delta}(v-\delta) \d\sigma
+\int_{V\setminus A_n^\delta} (v+\epsilon)\d\sigma.$$
D'o\`u 
$$\sigma(A_n^\delta)\leq \frac{1}{\delta}\left(\epsilon +\int_V
v\d\sigma -\int_V v_n\d\sigma\right).$$
Donc pour tout $\delta$ fix\'e, on a
$\limsup_{n\rightarrow \infty}\sigma(A^\delta_n)\leq \epsilon/\delta$
pour tout $\epsilon>0$. Par cons\'equent, $\lim \sigma(A^\delta_n)=0$.
La
suite $(v_n)$ \'etant born\'ee, on en d\'eduit que $v_n\rightarrow v$
dans $\mbox{L}^p(\sigma)$. 
\par
Le lemme de Fatou entra\^{\i}ne que 
$$\int_V v\d\sigma =\limsup \int_V v_n \d\sigma \leq 
\int_V\limsup v_n \d\sigma.$$
On obtient donc $\int(\limsup v_n-v)\d\sigma\geq 0$. 
Or d'apr\`es le lemme de
Hartogs, le terme sous le signe est n\'egatif. Donc $\limsup v_n=v$
$\sigma$-presque partout. 
\end{preuve}
{\it Fin de la d\'emonstration du th\'eor\`eme 3.2.1}---
Soit $\Omega$ une forme volume dans $\Ltwo(V)$, de masse 1. On peut
supposer $\Omega$ \`a support compact. La suite
$d_t^{-n}(f^n)^*\Omega$ n'admet qu'une valeur d'adh\'erence. En
effet, pour toute $\varphi$ p.s.h. on a
$$\int\varphi \frac{(f^n)^*\Omega}{d_t^n}
=\int\varphi_n\Omega\longrightarrow c_\varphi$$
car $(\varphi_n)$ converge dans $\Ltwoloc(V)$ vers $c_\varphi$, de plus,
les fonctions p.s.h. engendrent un espace dense dans $\Ltwoloc(V)$.
\par
Soit $\Omega$ une forme volume de masse 1 
\`a support compact dans $V\setminus \K$. 
Alors $d_t^{-n}(f^n)^*\Omega$ tend
vers une mesure $\mu$ et son support est contenu dans $U_{-n}\setminus
  \K$. Par cons\'equent, $\mu$ est port\'ee par $\partial \K$.
\par
Lorsque $v$ est une fonction p.s.h. de classe ${\cal C}^2$ on a $\int
v_n\d\mu =\int v \d\mu =c_v$. On a vu que $v_n\rightarrow c_v$ dans
$\Ltwoloc(V)$, il r\'esulte de la proposition 3.2.3 que
$v_n\rightarrow c_v$ dans $\Ltwo(\mu)$. Les op\'erateurs
$d_t^{-n}(f^n)_*$ 
sont de norme 1. De plus, l'espace vectoriel engendr\'e par les
fonctions p.s.h. born\'ees \'etant dense dans $\Ltwo(\mu)$, la
convergence de $v_n$ vers $c_v$ dans $\Ltwo(\mu)$ pour toute $v$ en
d\'ecoule. D'apr\`es la proposition 2.2.2, 
$\mu$ est K-m\'elangeante.
\par
\hfill $\square$
\begin{theoreme} Soit $f:U\longrightarrow V$ une application d'allure
  polynomiale. Les propri\'et\'es suivantes sont alors v\'erifi\'ees:
\begin{enumerate}
\item Pour toute mesure de probabilit\'e $\nu_0$ ne chargeant pas les
  ensembles pluripolaires, 
  $\nu_n:=d_t^{-n}(f^n)^*\nu_0$ tend faiblement vers $\mu$. 
\item Soit $\nu$ une mesure de probabilit\'e v\'erifiant
  $f^*\nu=d_t\nu$. Si $\nu$ ne charge pas les ensembles pluripolaires
  alors $\nu=\mu$. Si $\nu$ est ergodique et n'est pas port\'ee par un
  ensemble pluripolaire alors $\nu=\mu$.
\item Soient $K$ un compact de $U$ et $c>0$.
La suite $\nu_n$ converge vers $\mu$ uniform\'ement sur les
  mesures $\nu$ contenues dans $\M^p_c(K,U)$ 
(\voir la d\'efinition au paragraphe 3.9).
\end{enumerate}
\end{theoreme}
\begin{preuve}
1. Soit $\varphi$ une fonction lisse strictement 
p.s.h. dans $V$. Posons
$c_\varphi:=\int \varphi\d\mu$. On sait que
$\varphi_n:=d_t^{-n}(f^n)^*\varphi$ converge dans
$\Ltwoloc(V)$ vers $c_\varphi$. D'apr\`es la proposition 3.9.4, 
pour toute suite
$(\varphi_{n_j})$, on peut extraire une sous-suite $(\varphi_{m_j})$
convergeant en dehors d'un ensemble pluripolaire vers $c_\varphi$. Si
$\nu_{n_j}\rightharpoonup \nu$, d'apr\`es le th\'eor\`eme de
convergence domin\'ee, on a 
$$\int\varphi \d\nu_{m_j}= \int \varphi_{m_j} \d \nu_0 \longrightarrow
c_\varphi= \int\varphi \d\mu.$$ 
Par cons\'equent, 
les mesures $\nu$ et $\mu$ co\"{\i}ncident sur l'ensemble
des fonctions lisses strictement p.s.h. qui engendrent un sous-espace
dense. On a donc $\mu=\nu$.
\par
2. Si $\nu$ ne charge pas les ensembles pluripolaires et si
$f^*\nu=d_t\nu$, on a d'apr\`es la proposition 3.9.4
$$\int \varphi \d\nu = \int \varphi_{m_i} \d \nu \longrightarrow
c_\varphi= \int \varphi\d \mu.$$
Donc $\nu=\mu$.
\par 
Si une mesure $\nu$ v\'erifiant $f^*\nu=d_t\nu$ n'est pas port\'ee
par un ensemble 
pluripolaire et est ergodique, elle ne charge pas les ensembles
pluripolaires comme on le v\'erifie ais\'ement. On a alors $\nu=\mu$. 
\par
3. 
Rappelons qu'une mesure $\nu$ port\'ee par $K$ appartient \`a 
$\M^p_c(K,U)$ si pour toute fonction p.s.h. $\psi$ dans $U$ on a 
$\|\psi\|_{\Lone(\nu)}\leq c\|\psi\|_{\Lp(U)}$.
Soient $\nu^j\in\M^p_c(K,U)$ et $(m_j)\rightarrow \infty$ tels que 
$d_t^{-m_j}(f^{m_j})^*\nu^j$ tend faiblement vers une mesure
$\nu$. Soit $\varphi$ une fonction p.s.h. lisse dans $V$. Alors la
suite $(\varphi_{m_j})$ tend vers $c_\varphi$ dans $\Lploc(V)$. 
Par d\'efinition de
$\M^p_c(K,U)$, on a 
$$\left|\int \varphi_{m_j} \d\nu^j 
-c_\varphi\right| = 
\left|\int(\varphi_{m_j}-c_\varphi)\d \nu^j\right|\leq
c\|\varphi_{m_j}-c_\varphi\|_{\Lp(U)}.$$
Donc $\int \varphi \d\nu =c_\varphi= \int\varphi \d\mu$. On en d\'eduit
comme pr\'ec\'edemment que $\nu=\mu$.
\end{preuve}
\par
Posons pour tout $z\in V$
$$\mu_n^z:=\frac{(f^n)^*\delta_z}{d_t^n}=\frac{1}{d_t^n}
\sum_{f^n(w)=z} \delta_w$$
o\`u $\delta_z$ d\'esigne la masse de Dirac en $z$. Posons \'egalement
$$\E':=\big\{z\in V, \ \mu \mbox{ \rm n'est pas adh\'erente \`a la
  suite }  (\mu^z_n) \big\}$$
$$\E:=\big\{z\in V, \ \mu^z_n \mbox{ ne tend pas vers } \mu
\big\}$$
et pour toute suite d'entiers $(n_j)$ tendant vers l'infini
$$\E^{(n_j)}:=\big\{z\in V, \ \mu \mbox{ \rm n'est pas adh\'erente \`a la
  suite }  (\mu^z_{n_j}) \big\}.$$
Observons que $\E$ est la r\'eunion des $\E^{(n_j)}$ et que $\E'$ est
  l'intersection des  $\E^{(n_j)}$.
Notons $\PSH(.)$ le c\^one des fonctions p.s.h.
On a la proposition suivante:
\begin{proposition} Soit $f$ comme au th\'eor\`eme 3.2.1 et soit $p\geq
  1$. Alors 
\begin{enumerate}
\item L'op\'erateur $\Lambda:=d_t^{-1}f_*:\PSH(W)\cap
  \Lp(W)\longrightarrow \PSH(U)\cap\Lp(U)$ est born\'e pour 
pour tout ouvert $W$ contenant $\overline U_{-2}$.
\item Il existe une constante $0<c_1<1$ telle que pour toute $\varphi$
  pluriharmonique dans $V$ on a 
$\sup_U|\Lambda\varphi - c_\varphi|\leq
  c_1\sup_U|\varphi-c_\varphi|$. En particulier, $\varphi_n$ tend vers
  $c_\varphi$ uniform\'ement et g\'eom\'etriquement.
\item Pour toute fonction $\varphi$, quasi-p.s.h. au voisinage de 
$\K$ et $\mu$-int\'egrable, on a
  $\varphi_n\rightarrow c_\varphi$ dans $\Lploc(V)$; si $\int\varphi\d
  \mu =-\infty$, $\varphi_n$ tend uniform\'ement vers $-\infty$.
\item Les ensembles $\E^{(n_j)}$ et $\E'$ sont pluripolaires. De plus, 
si $\varphi$ est une fonction strictement p.s.h. au
voisinage de $\K$ et  $\mu$-int\'egrable, on a 
\begin{eqnarray*}
\E^{(n_j)} & = & \big\{z\in V,\ 
\limsup \varphi_{n_j}(z)<c_\varphi\big\} \\
\E' & = & \big\{z\in V, \ \limsup\varphi_n(z)<c_\varphi\big\}\\
\E  & = & \big\{z\in V,\ \liminf\varphi_n(z)<c_\varphi\big\}.
\end{eqnarray*}
Si la s\'erie $\sum\|\varphi_n-c_\varphi\|_{\Lp(U)}$ converge, 
alors $\E$ est pluripolaire.  
\item Si $X$ est un ferm\'e non pluripolaire de $V$ v\'erifiant
  $f^{-1}(X)\subset X$, alors $X$ contient l'ensemble de Julia $J_k$. 
\end{enumerate}
\end{proposition}
\begin{preuve} 
1. Sinon, il existe $\varphi^j$ p.s.h. sur $W$ telles que 
$\|\Lambda\varphi^j\|_{\Lp(U)}=1$ et $\varphi^j$ tend vers 0 sur $W$. 
Il est clair que ceci contredit le fait que $W$ contient $\overline U_{-2}$.
\par
2. On peut supposer $c_\varphi=0$. La propri\'et\'e 2 est une
cons\'equence directe du principe du maximum et du fait que la famille
$\{\varphi \mbox{ pluriharmonique }, c_\varphi=0, |\varphi|\leq 1
\mbox{ sur } U\}$ 
est compacte. 
\par
3. Rappelons que $\varphi$ est {\it quasi-p.s.h.}
si $\ddc\varphi\geq -c\omega$,
$c>0$. Puisque $\varphi$ s'\'ecrit au voisinage de $\K$ 
comme diff\'erence de deux fonctions
p.s.h. qui sont $\mu$-int\'egrables, il suffit de
consid\'erer le cas o\`u $\varphi$ est p.s.h. au voisinage de $\K$.
Supposons que $\varphi$ est $\mu$-int\'egrable.
Pour raisonner comme au lemme 3.2.2, il suffit de montrer que 
$(\limsup \varphi_n)^*$ n'est pas identiquement \'egale \`a $-\infty$. 
Sinon, $\varphi_n$ converge  
uniform\'ement vers $-\infty$. Par cons\'equent, $\langle\mu,\varphi_n
\rangle$ tend vers $-\infty$. Ceci est impossible car
$\langle\mu,\varphi_n \rangle =\langle 
d_t^{-n}(f^n)^*\mu,\varphi\rangle=
\langle \mu,\varphi \rangle$.
La deuxi\`eme partie est aussi claire.
\par
4. Montrons que si $\varphi_{n_j}(z)\longrightarrow
  c_\varphi$ alors $\mu_{n_j}^z\rightharpoonup\mu$. Soit $\psi$
  une fonction de classe ${\cal C}^2$ \`a support compact. Pour
  $\epsilon>0$ suffisamment petit les fonctions
  $\varphi^+:=\varphi+\epsilon\psi$ et
  $\varphi^-:=\varphi-\epsilon\psi$ sont p.s.h. D'apr\`es le lemme
  3.2.2, la suite de fonctions $\varphi^+_n$ (resp. $\varphi^-_n$) tend
  dans $\Ltwoloc(V)$ vers une constante $c_{\varphi^+}$
  (resp. $c_{\varphi^-}$). Il en r\'esulte que $\psi_n$ converge dans
  $\Ltwoloc(V)$ vers la constante 
$$c_\psi:=\frac{1}{\epsilon}(c_{\varphi^+}-c_\varphi)
=\frac{1}{\epsilon}(c_\varphi-c_{\varphi^-}).$$
Si $\psi_{n_j}(z)\longrightarrow \alpha$, on a, gr\^ace au lemme de
Hartogs, 
$c_\varphi+\epsilon\alpha\leq c_{\varphi^+}$ et 
 $c_\varphi-\epsilon\alpha\leq c_{\varphi^-}$. Donc
$\alpha=c_\psi$. On a donc montr\'e que
$$\E^{(n_j)}=\big\{z\in V,\ \limsup \varphi_{n_j}(z)<c_\varphi\big\}.$$
Cet ensemble est ind\'ependant de la fonction strictement p.s.h 
et $\mu$-int\'egrable $\varphi$.
Or l'ensemble o\`u $\limsup
\varphi_{n_j}(z)<(\limsup\varphi_{n_j})^*=c_\varphi$ est
pluripolaire. L'ensemble $\E^{(n_j)}$ est pluripolaire.
Les assertions analogues sur $\E'$ et $\E$ sont imm\'ediates. 
\par
Montrons
que $\E$ est pluripolaire si la s\'erie $\sum
\|\varphi_n-c_\varphi\|_{\Lp(U)}$ converge. D'apr\`es
la proposition 3.9.2, la s\'erie
$\sum \|\varphi_n-c_\varphi\|_{\Lp(\nu)}$ est aussi convergente pour
toute mesure PLB 
$\nu$ \`a support compact dans $U$. Il en r\'esulte
que $\E=\{z\in V,\ \lim \varphi_n(z)<c_\varphi\}$ est de $\nu$ mesure 
nulle. 
Comme dans la proposition 3.9.4, on d\'eduit que
$\E$ est pluripolaire.
\par
Bien s\^ur, si la s\'erie  
$\sum\|\varphi_{n_j}-c_\varphi\|_{\Lp(U)}$ converge, 
on a aussi $\mu^z_{n_j}
\rightharpoonup\mu$ sauf pour un ensemble pluripolaire $\E^{(n_j)}$.
\par
5. Puisque $\E'$ est
pluripolaire, on a $X\not \subset \E'$. Il existe donc
$z\in X$ tel que $\mu$ est adh\'erente \`a la suite $(\mu_n^z)$. Comme
$X$ est un ferm\'e totalement invariant, $\mu_n^z$ est \`a support
dans $X$. Par cons\'equent, le support $J_k$ de $\mu$ est contenu dans $X$.
\end{preuve}
\par
Notons $\S(\K)$
le c\^one convexe des fonctions continues p.s.h. au
voisinage de $\K$.  
On d\'efinit {\it l'enveloppe par rapport \`a $\S(\K)$}
de l'ensemble de Julia $J_k=\supp(\mu)$
dans $\K$ par
$$\widehat J_k:=\big\{z\in \K,\ \varphi(z)\leq \sup_{J_k}
\varphi(w) \mbox{ pour toute } \varphi \in \S(\K)\big\}. $$
On a la proposition suivante:
\begin{proposition} Soient $\nu_0$ une
mesure de probabilit\'e de
$V$ et $\nu$ une valeur adh\'erente de la suite
$d_t^{-n}(f^n)^*\nu_0$. On a
$\int\varphi \d\nu\leq \int \varphi \d\mu$
pour toute fonction $\varphi$ p.s.h. au voisinage de $\K$ et 
$\mu$-int\'egrable. De plus,
$\nu$ est port\'ee par $\widehat J_k$. 
Si $\nu_0$ ne charge pas $\E$, on a $\nu=\mu$. Si $\nu\not=\mu$, on a
$\int \varphi \d\nu <\int \varphi \d\mu$ pour toute fonction $\varphi$
strictement p.s.h. au voisinage de $\K$ et $\mu$-int\'egrable. 
Si $u$ est une fonction
pluriharmonique au voisinage de $\K$, on a $\int u\d\nu =\int u \d\mu$. 
\end{proposition}
\begin{preuve} Soit $(n_i)$ une suite croissante d'entiers telle
que $d_t^{-n}(f^{n_i})^*\nu_0$ tende vers $\nu$. Puisque $\varphi$ est
s.c.s, on a 
\begin{eqnarray*}
\int\varphi \d\nu & \leq & \limsup \int \varphi \d\nu_{n_i}
=\limsup \int\varphi_{n_i} \d\nu_0 \leq \int \limsup
\varphi_{n_i} \d\nu_0\\
& \leq & \int c_\varphi \d\nu_0=c_\varphi=\int
\varphi \d\mu.
\end{eqnarray*}
\par
On sait que
la suite $\varphi_n$ converge ponctuellement vers $c_\varphi$ sur
$V\setminus\E$. Si $\nu_0$ ne charge pas $\E$ et si $\varphi$ est continue, les
in\'egalit\'es ci-dessus deviennent des \'egalit\'es.
D'o\`u $\nu=\mu$.
\par
Montrons que le support de $\nu$ est contenu dans
$\widehat{J_k}$. Sinon, soit $z_0\in\supp (\nu)
\setminus\widehat{J_k}$. Il existe une fonction p.s.h.
continue $\psi$ au voisinage de $\K$ telle que $\psi(z_0)>0$
et $\psi<0$ sur $J_k$. Soit $\psi^+(z):=\max(0,\psi(z))$.
Cette fonction p.s.h. ne peut v\'erifier $\int \psi^+\d\nu\leq
\int\psi^+ \d\mu$.
\par
Supposons qu'il existe une fonction $\varphi$ strictement p.s.h. au
voisinage de $\K$ et $\mu$-int\'egrable 
telle que $\int \varphi \d\nu = \int \varphi
\d\mu$. Montrons que $\nu=\mu$. Soit $\psi$ une fonction ${\cal
  C}^2$. Posons $c_\psi:=\int \psi\d\mu$. 
Soit $\epsilon>0$ telle que $\varphi^+:=\varphi+\epsilon \psi$ soit
p.s.h. La fonction $\varphi^+$ \'etant $\mu$-int\'egrable et s.c.s, on a 
$$\int\varphi^+ \d\nu\leq \limsup \int\varphi^+_n \d\nu_0 \leq \int
\limsup \varphi^+_n \d\nu_0 
\leq c_{\varphi^+}=c_\varphi+\epsilon c_\psi.$$
Donc $\int \psi \d\nu\leq c_\psi=\int \psi \d\mu$.
Pour la fonction $-\psi$,
on obtient $\int -\psi \d\nu\leq - \int \psi \d\mu$. D'o\`u $\int \psi
\d\nu = \int \psi \d\mu$ et $\nu=\mu$.
\par
L'assertion sur les fonctions pluriharmoniques est claire puisque les
fonctions $u$ et $-u$ sont p.s.h.
\end{preuve}
\begin{remarques} \rm 
Si $\mu$ charge un sous-ensemble analytique $X$ de $V$ alors
$\mu$ est port\'ee par $X$. En effet, on peut supposer $X$
irr\'eductible et de dimension
minimale. Le th\'eor\`eme de r\'ecurrence de Poincar\'e permet de
supposer alors que $f(X)=X$. L'ergodicit\'e permet de conclure. 
En particulier, si $X$ est irr\'eductible de dimension minimale, il est
totalement invariant. 
Si $\mu$ charge un ensemble
pluripolaire $A$, alors $\mu$ est port\'ee par l'ensemble
pluripolaire $\bigcup_{n\geq 0} f^{-n} A$. 
%\par
%2. Si $f$ est un endomorphisme holomorphe de $\P^k$, son \'etude se
%ram\`ene \`a ce cadre. Il suffit, en effet, de consid\'erer
%l'application holomorphe $F:\C^{k+1}\longrightarrow \C^{k+1}$ telle
%que $f\circ\pi=\pi\circ F$ o\`u $\pi$ d\'esigne l'application
%canonique de $\C^{k+1}\setminus\{0\}$ dans $\P^k$. L'application $F$
%est d'allure polynomiale au sens pr\'ec\'edent et les r\'esultats
%dynamiques que nous consid\'erons pour $f$ se d\'eduisent facilement
%des propri\'et\'es correspondantes pour $F$. 
\end{remarques}
\par
Notons $\partial_S \K$ 
la fronti\`ere de
Silov 
de $\K$ associ\'ee aux fonctions continues p.s.h. aux voisinage
de $\K$. Rappelons que {\it la fronti\`ere de Choquet}
$\partial_C \K$ de $\K$ 
est
l'ensemble des points $x$ de $\K$ pour lesquels la seule mesure de
probabilit\'e $\nu$ \`a support dans $\K$ satisfaisant
$$u(x)\leq \int u\d\nu \mbox{ pour toute } u\in \S(\K)$$
est la masse de Dirac $\delta_x$ en $x$. 
On a $\overline{\partial_C \K}=\partial_S
\K$ \cite{Meyer}. 
De plus, pour tout $x\in \K$, il existe une mesure de probabilit\'e
$\nu$ port\'ee par $\partial_C \K$ telle que 
$$u(x)\leq \int u\d\nu \mbox{ pour toute } u\in\S(\K).$$
On dit que $\nu$ est {\it une mesure repr\'esentante}
de $x$.
\begin{corollaire}
Soit $f:U\longrightarrow V$ comme au th\'eor\`eme
  3.2.1. Si $V$ est S-s\'epar\'e, $\mu$ est port\'ee par $\partial_S
  \K$. Si $V$ est un ouvert d'une vari\'et\'e de Stein, 
$f$ poss\`ede au plus un
point fixe totalement invariant. Si $V$ est une vari\'et\'e de Stein
  et si $X_1$ et
  $X_2$ sont deux sous-ensembles analytiques
totalement invariants alors
  $X_1\cap X_2\not=\emptyset$. 
\end{corollaire}
\begin{preuve} Supposons que $V$ est S-s\'epar\'ee. Montrons d'abord
  que  $\partial_C \K$ est totalement invariant par $f$.
Soit $z\in \partial_C \K$. 
Notons $w_1$, $\ldots$, $w_m$ les images
r\'eciproques de $z$. 
Soit $\nu$ une mesure \`a support dans $\K$ repr\'esentant $w_1$. La
mesure $f_*\nu$ repr\'esente $z$ et
donc $f_*\nu=\delta_z$ car $z\in \partial_C \K$. On en d\'eduit que
$\nu=\sum_{j=1}^m \alpha_j\delta_{w_j}$ avec $\sum\alpha_j=1$ et
$\alpha_j\geq 0$. Comme $V$ est S-s\'epar\'ee, on a
$\nu=\delta_{w_1}$. 
Par suite, $f^{-1}(\partial_C
\K)\subset \partial_C \K$ et donc  $f^{-1}(\partial_S
\K)\subset \partial_S \K$. 
\par
Puisque tout point de $\K$ est
repr\'esent\'e par une mesure port\'ee par
$\partial_C \K$, pour toute
mesure de probabilit\'e $\sigma$ port\'ee par $\K$, il existe une
mesure de probabilit\'e $\widehat\sigma$ port\'ee par $\partial_S
\K$ telle que 
$$\int\varphi \d\sigma \leq \int \varphi \d\widehat
\sigma \mbox{ pour toute }\varphi\in \S(\K).$$
Appliquant cela \`a $\mu$, on a $\int \varphi
\d\mu\leq \int\varphi \d\widehat \mu$.
\par
Soit $\tilde \mu$ la limite d'une suite
$d_t^{-n}(f^{n_i})^*\widehat\mu$.
On a, d'apr\`es l'observation
pr\'ec\'edente,
$$\int\varphi \d\mu=\int\varphi_{n_i} \d\mu \leq
\lim\int \varphi_{n_i} \d\widehat \mu
=\int\varphi d\tilde \mu.$$
D'apr\`es la proposition 3.2.5,
la derni\`ere int\'egrale est major\'ee par $\int \varphi \d\mu$. 
Donc $\mu=\tilde \mu$. Il est clair que $\tilde \mu$ est port\'ee
par $\partial_S \K$ car ce dernier ensemble est totalement
invariant.
\par
Soient $X_1$ et $X_2$ deux sous-ensembles analytiques totalement
invariants.
En utilisant des limites de sommes de C\'esaro, 
on peut construire des mesures de
probabilit\'e  
$\nu_1$, $\nu_2$  
port\'ees par $X_1$, $X_2$ et v\'erifiant
$f^*\nu_j=d_t\nu_j$.  D'apr\`es la proposition 3.2.5, 
on a $\int u\d\nu_1=\int
u\d\nu_2$ pour $u$ pluriharmonique. Or si $X_1$ et $X_2$
sont disjoints et si $V$ est une vari\'et\'e de Stein,
il existe une telle fonction $u$ avec
$u=0$ sur $X_1$ et $u>0$ sur $X_2$ ce qui est impossible. 
\par
Supposons maintenant que $V$ 
est un ouvert d'une vari\'et\'e de Stein. Alors les fonctions
pluriharmoniques de $V$ s\'eparent les points. De la m\^eme mani\`ere,
on montre qu'il existe au plus un point fixe totalement invariant.
\end{preuve}
\begin{proposition} On a $\mu(\E)=0$ ou $1$. On a aussi $\mu(\E')=0$.
L'ensemble $\E$ est vide 
dans les deux cas suivants: 
\begin{enumerate}
\item L'ensemble $\K$ est $B$-r\'egulier,
\cad que les fonctions continues sur $\K$ sont uniform\'ement
approximables par des fonctions p.s.h. au voisinage de $\K$.
\item L'ensemble critique $\Cr$ de $f$ ne rencontre pas
$\K$.
L'ensemble $\K$ est alors $B$-r\'egulier, 
$\mu_n^z$ tend vers $\mu$ uniform\'ement
sur $V$. De plus, $J_k=\K$ et $f$ est hyperbolique (\cad expansive) sur $\K$.
\end{enumerate}
\end{proposition}
\begin{preuve} 
Soit $\varphi$ une fonction continue strictement p.s.h. 
Posons $\widehat\varphi:=\liminf
\varphi_n$. On a 
$$\frac{f_*\widehat \varphi}{d_t}=\frac{f_*(\liminf\varphi_n)}{d_t} \leq
\frac{\liminf f_*\varphi_n}{d_t}=\liminf
\varphi_{n+1}=\widehat\varphi$$ 
et
$$\widehat\varphi\leq (\limsup \varphi_n)^*= c_\varphi.$$ 
D'apr\`es la proposition 3.2.5, on a $\E=\{\widehat\varphi
<c_\varphi\}$ et donc $f(\E\cap U)\subset \E$. 
Soit $z\in V$ un point tel que
$d_t^{-1}f_*\widehat\varphi(z) <c_\varphi$. On a
$f^{-1}(z)\cap \E\not=\emptyset$. 
On en d\'eduit que $z\in\E$ et donc
$\widehat \varphi(z)<c_\varphi$. 
Par cons\'equent, $\widehat\varphi(z) <c_\varphi$
si et seulement si $d_t^{-1}f_*\widehat\varphi(z) <c_\varphi$. Soit
$\tilde \varphi$ la limite d\'ecroissante de la suite
$d_t^{-1}(f^n)_*\widehat\varphi$. On a
$d_t^{-1}f_*\tilde\varphi=\tilde\varphi$ et $\E=\{\tilde\varphi
<c_\varphi\}$.  
Comme $\mu$ est ergodique,
$\tilde\varphi$ est constante $\mu$-presque partout. Si $\tilde
\varphi=c_\varphi$ $\mu$-presque partout, 
on a $\mu(\E)=0$. Sinon, $\mu(\E)=1$. 
\par
On sait que $\limsup \varphi_n\leq c_\varphi$. D'apr\`es le lemme de
Fatou
$$c_\varphi=\int\varphi_n \d\mu \leq \int \limsup \varphi_n \d
\mu \leq \int c_\varphi\d\mu =c_\varphi.$$
Par cons\'equent, l'ensemble $\E'=\{\limsup \varphi_n <c_\varphi\}$ 
est de mesure $\mu$ nulle. 
\par
1. Soit $\nu$ une valeur d'adh\'erence de la suite
$(\mu_n^z)$. D'apr\`es la proposition 3.2.5,
$\int \psi \d\nu\leq \int \psi \d\mu$ pour toute fonction 
$\psi$ p.s.h. au voisinage de $\K$. Comme $\K$ est $B$-r\'egulier, 
cette in\'egalit\'e est valable pour toute fonction
continue $\psi$, en particulier pour $-\psi$. D'o\`u
$\nu=\mu$.
\par
2. Soit $u$ une fonction continue au voisinage de $\K$. Fixons un
$n_0$ suffisament grand tel que $\Cr\cap U_{-n_0}=\emptyset$
et tel que $u$ soit uniform\'ement continue sur $U_{-n_0}$.
Rappelons que l'existence d'une fonction strictement p.s.h. born\'ee
$\Phi$ dans $U$
entra\^{\i}ne que $U$ est Kobayashi hyperbolique \cite{Sibony1}.
Notons $\K_n$ la m\'etrique infinit\'esimale de Kobayashi de
$U_{-n}$.
L'application $f$ est une isom\'etrie de
$(U_{-n_0-1},\K_{n_0+1})$ dans $(U_{-n_0},\K_{n_0})$.  Puisque
$U_{-n_0-1}\subset\subset U_{-n_0}$, il existe $\delta>0$ telle que
$$\K_{n_0}(f(z),f'(z)\zeta)=\K_{n_0+1}(z,\zeta)\geq (1+\delta)
\K_{n_0}(z,\zeta)$$
pour tout $z\in U_{-n_0-1}$ et tout vecteur tangent $\zeta$ de
$V$ en $z$.
Donc $f$ est expansive sur $U_{-n_0-1}$. En particulier, $f$ est
hyperbolique sur $\K$ au sens dynamique. Plus pr\'ecis\'ement, on a
$|(f^n)'(z)\eta|\geq (1+\delta)^n|\eta|$ pour $z\in \K$.
\par
Dans la suite, on note $d$ la distance de Kobayashi sur $U_{-n_0}$.
Il existe une constante $c>0$ telle que
pour tous $z$ et $w$ dans $V$ et tout $n\geq n_0$ on peut
ranger les points $z_i$ de $f^{-n}(z)$ et $w_i$ de $f^{-n}(w)$ de
sorte que $d(z_i,w_i)\leq c(1+\delta)^{-n+n_0}$. 
Comme $u$ est uniform\'ement continue, $|u_n(z)-u_n(w)|$ tend
uniform\'ement vers $0$.
Ceci implique que $|u_n(z)-u_n(w)|$ tend uniform\'ement vers $0$
pour tous $z$ et $w$ dans $V$. 
D'o\`u la convergence uniforme de $\mu_n^z$ sur $V$. 
\par
Pour montrer que $\K$ est $B$-r\'egulier, on observe que
$\Phi\circ f^n$ est born\'ee au voisinage de $\K$ mais que $\ddc (\Phi\circ
f^n)$ tend vers l'infini uniform\'ement sur $\K$. Toute fonction $u$ de
classe ${\cal
  C}^2$ peut \^etre donc perturb\'ee en une fonction
$u+\epsilon (\Phi\circ f^n)$ qui est p.s.h. au voisinage de $\K$.
\par
Montrons que le support de $\mu$ est \'egal \`a $\K$.
On utilise encore la
distance de Kobayashi. L'ensemble $J_k$ est
totalement invariant et $f$ est expansive.
Pour tout $z\in U$, on a $\dist(f(z),J_k)\geq
(1+\delta)\dist(z,J_k)$. Par suite, $\dist(f^n(z),J_k)\geq
(1+\delta)^n \dist(z,J_k)$ lorsque $z\in U_{-n}$. Il en r\'esulte
que si $z\not\in J_k$ on a
$z\not\in\bigcap_{n\geq 0} U_{-n}$. Donc $J_k=\bigcap_{n\geq 0} U_{-n}=\K$.
\end{preuve}
\begin{remarque} \rm Si au point 2 de la proposition pr\'ec\'edente,
  on suppose que $u$ est une fonction
H\"old\'erienne
alors $|u_n(z)-u_n(w)|\leq c_1
(1+\delta)^{-n\alpha}$ avec $c_1>0$ et $0<\alpha<1$. 
Autrement dit, la suite $u_n(z)-\int u
\d\mu$ tend uniform\'ement vers $0$ \`a vitesse exponentielle. Pour ceci, il
suffit d'observer que 
$$u_n(z)-\int u\d\mu = \int \big[u_n(z)-u_n(w) \big]\d
\mu(w).$$ 
En particulier, la vitesse de m\'elange 
est
exponentiellement petite, \cad qu'il existe $c>0$ tel que pour toute fonction
$\varphi$ de classe ${\cal C}^1$ et toute fonction born\'ee $\psi$, on
ait
$$\left|\int \psi(f^n) \varphi \d\mu -\int\psi \d\mu\int\varphi\d\mu 
\right|\leq c\|\varphi\|_{{\cal
    C}^1}\|\psi\|_\infty (1+\delta)^{-n}.$$  
\end{remarque}
\subsection{Degr\'es dynamiques et entropie}
Nous allons donner des majorations des degr\'es dynamiques d'une
application d'allure polynomiale $f:U\longrightarrow V$.
Utilisant ces estimations
et un th\'eor\`eme de Gromov-Yomdin, nous allons montrer, dans le
cas o\`u $V$ est un ouvert d'une vari\'et\'e de Stein, que
l'entropie de $f$ est \'egale \`a $\log d_t$. 
La mesure d'\'equilibre $\mu$ est donc d'entropie maximale.
\begin{proposition} Soit $f:U\longrightarrow V$ une application
d'allure polynomiale comme pr\'ec\'edemment. Alors
$d_{k-1,n}=\o(d_t^n)$ et $\delta_n=\o(1)$. En
particulier, on a $d_{k-1}\leq d_t$.
\end{proposition}
\begin{preuve} Si $\varphi$ est une fonction quasi-p.s.h. au
voisinage de $\K$ et $\mu$-int\'egrable, d'apr\`es la proposition
3.2.5, $\varphi_n$ converge vers une constante $c_\varphi$ dans
$\Ltwoloc(V)$. Par cons\'equent, $\ddc \varphi_n$ tend faiblement
vers $0$. On a
$$\int_U (f^n)_*(\ddc\varphi)\wedge \omega^{k-1}=d_t^n\int_U \ddc
\varphi_n \wedge \omega^{k-1} =\o(d_t^n).$$
En appliquant, cette relation aux fonctions $\Phi$ et $\log J$,
on obtient $d_{k-1,n}=\o(d_t^n)$ et $\delta_n=\o(1)$. Rappelons ici
que $\ddc\Phi=\omega$ et que 
$\log J$ est une fonction quasi-p.s.h., $\mu$-int\'egrable.
\end{preuve}
\par
Le th\'eor\`eme suivant permet de calculer l'entropie.
\begin{theoreme} 
Soit $f:U\longrightarrow V$ une application
d'allure polynomiale de degr\'e topologique $d_t\geq 2$.
Si $V$ est un ouvert d'une vari\'et\'e de Stein, alors
l'entropie topologique de $f$ est \'egale \`a $\log d_t$. Dans ce cas,
$\mu$ est une mesure invariante d'entropie maximale. 
\end{theoreme}
\par
Soient $W_1$, $W_2$ deux vari\'et\'es complexes
de dimensions $k_1$ et $k_2$ munie de m\'etriques hermitiennes. Soient 
$K_1\subset W_1$, $K_2\subset W_2$ des compacts et $m\geq 1$ un entier.
Notons $\pi:W_1^m\times W_2\longrightarrow
W_2$ la projection canonique. 
Soit $\Gamma\subset K_1^m\times W_2$
un sous-ensemble
analytique de dimension pure $k_2$ de $W_1^m\times W_2$. Alors la
restriction de $\pi$ \`a $\Gamma$ r\'ealise un rev\^etement
ramifi\'e au dessus de $W_2$. Notons $d_\Gamma$ le degr\'e de ce
rev\^etement.
Le lemme suivant 
nous sera utile pour estimer les degr\'es dynamiques et pour
d\'emontrer le th\'eor\`eme 3.3.2.
\begin{lemme} Supposons que $W_1$ soit un ouvert d'une vari\'et\'e
de Stein. Alors il existe une constante $c>0$ et un entier $s$,
ind\'ependants de
$\Gamma$ et de $m$, tels que le volume de $\Gamma\cap
(W_1^m\times K_2)$ soit major\'e par $cm^sd_\Gamma$. 
\end{lemme}
\begin{preuve}
Pour simplifier les calculs, la norme consid\'er\'ee dans la
suite pour tout espace complexe $\C^n$ sera la somme des modules
des coordonn\'ees. Ceci est aussi
valable pour d\'efinir les boules. En
revanche, les volumes seront calcul\'es pour la m\'etrique
euclidienne. 
\par  
Puisque toute vari\'et\'e de Stein est plongeable
dans un $\C^N$, on peut supposer que $W_1$ est \'egale \`a
$\C^{k_1}$ et que $K_1$ est la boule unit\'e ferm\'ee de $\C^{k_1}$.
Le probl\`eme est local pour $W_2$. On peut supposer
que $W_2$ est la boule unit\'e de $\C^{k_2}$ et $K_2$ est une
boule ferm\'ee
de rayon $0<r<1$ centr\'ee en $0$. 
\par
Notons $x=(x_1,\ldots , x_{mk_1})$ et $y=(y_1,\ldots,y_{k_2})$ les
coordonn\'ees de $(\C^{k_1})^m$ et de $\C^{k_2}$. Soit $\epsilon$ une
matrice de taille $k_2\times mk_1$ dont les coefficients sont
major\'es par
$(1-r)/2m^2$. Posons $\pi_\epsilon(x,y):=y+\epsilon x$,
$\Gamma_\epsilon:=\Gamma\cap \{\|\pi_\epsilon\|< (1+r)/2\}$ et
$\Gamma^*:=\Gamma\cap (W_1^m\times K_2)$.
\par
Montrons d'abord que $\Gamma^*\subset\Gamma_\epsilon$. Soit
$(x,y)\in \Gamma^*$. On a $\|x\|< 1$ et $\|y\|\leq r$. Par
cons\'equent,
$$\|\pi_{\epsilon}(x,y)\|\leq \|y\| +\|\epsilon x\|< r+ (1-r)/2
= (1+r)/2.$$
D'o\`u $(x,y)\in\Gamma_\epsilon$.
\par
Montrons maintenant que pour tout $a\in\C^{k_2}$ v\'erifiant
$\|a\|<(1+r)/2$, on a $\#\pi_\epsilon^{-1}(a)\cap \Gamma =d_\Gamma$.
Pour ceci, on montre que $\#\pi_{t\epsilon}^{-1}(a)\cap \Gamma$
ne d\'epend pas de $t\in [0,1]$. Il suffit donc de montrer
que la r\'eunion des ensembles $\pi_{t\epsilon}^{-1}(a)\cap \Gamma$
est contenue dans le compact
$\Gamma\cap\{\|\pi\|\leq (1+r)/2\}$ de
$\Gamma$. Soient $(x,y)\in\Gamma$ et
$t\in[0,1]$ tels que $\pi_{t\epsilon}(x,y)=a$. On a
$$(1+r)/2> \|a\| =\|\pi_{t\epsilon}(x,y)\| =\|y\| +
t\|\epsilon x\|.$$
On en d\'eduit que $\|y\|<(1+r)/2$ et donc $(x,y)\in 
\Gamma\cap\{\|\pi\|<(1+r)/2\}$.
\par
Notons $\Omega$ la forme de volume euclidienne de $\C^{k_2}$ et $B$ la
boule de centre $0$ et de rayon $(1+r)/2$ de $\C^{k_2}$. On a
pour une constante $c'>0$
$$\int_{\Gamma^*} (\pi_\epsilon)^* \Omega \leq
\int_{\Gamma_\epsilon} (\pi_\epsilon)^*\Omega = d_\Gamma \int_B
\Omega =c'd_\Gamma.$$
\par
Notons $\Theta:=i\sum \d x_\nu\wedge \d\overline x_\nu +\pi^*\omega$ 
la $(1,1)$-forme euclidienne de $(\C^{k_1})^m\times
\C^{k_2}$ o\`u $\omega$ est celle de $\C^{k_2}$. 
Il suffit majorer la forme $\Theta$ par une combinaison
lin\'eaire finie de $2mk_1+1$ formes du type $(\pi_\epsilon)^*\omega$ dont 
les coefficients sont d'ordre $m^4$. 
Pour ceci, on majore $i\d x_\nu \wedge \d \overline x_\nu$ par une
combinaison de $(1,1)$-formes du type 
$(\pi_\epsilon)^*\omega$. Consid\'erons
$\delta:=(1-r)/2m^2$ et $\pi_\epsilon(x,y):=(y_1+\delta x_\nu,
y_2,\ldots,y_{k_2})$. On a
\begin{eqnarray*}
i\d x_\nu \wedge \d \overline x_\nu & = & \frac{4i}{3\delta^2} 
\Big[3\d y_1\wedge\d \overline y_1 + \d (y_1+\delta x_\nu) \wedge
\d (\overline y_1+\delta \overline x_\nu) \\
& & - \d (2y_1+\delta x_\nu/2) \wedge
\d (2\overline y_1+\delta \overline x_\nu/2)\Big] \\
& \leq &  \frac{4i}{3\delta^2} 
\Big[3\d y_1\wedge\d \overline y_1 + \d (y_1+\delta x_\nu) \wedge
\d (\overline y_1+\delta \overline x_\nu)\Big] 
\end{eqnarray*}
Ceci donne l'estimation du lemme.
\end{preuve}
\begin{corollaire} Soit $f:U\longrightarrow V$ une application
d'allure polynomiale de degr\'e topologique $d_t\geq 2$.
Si $V$ est un ouvert d'une vari\'et\'e de Stein, alors
pour tout $1\leq l\leq k$ on a $d_l\leq d_t$.
\end{corollaire}
\begin{preuve} On applique le lemme 3.3.3 pour $W_1=W_2=V$,
$K_1=K_2=\overline U$, $m=1$ et
$\Gamma$ le graphe de l'application $f^n$
dans $U_{-n}\times V$. Pour $n\geq 1$, la projection de $\Gamma$ sur $W_1$
est contenue dans $K_1$.
Puisque le degr\'e topologique de
$f^n$ est \'egal \`a $d_t^n$, on a $d_\Gamma=d_t^n$.
Or on a la majoration 
$$d_{l,n}\leq \int_{U_{-n-1}}\big(\omega+(f^n)^*\omega\big)^k 
\leq k!\vol(\Gamma\cap K_1\times K_2).$$
D'apr\`es le lemme 3.3.3, on
a $d_{l,n}\leq ck!d_t^n$. D'o\`u $d_l\leq d_t$.
\end{preuve}
{\it D\'emonstration du th\'eor\`eme 3.3.2}-- 
D'apr\`es le th\'eor\`eme 2.3.1, $\mu$ est une
mesure invariante d'entropie au moins $\log d_t$. Puisque l'entropie
topologique de $f$ est minor\'ee par les entropies des mesures
invariantes, on a $\h(f)\geq \log d_t$. Il reste \`a montrer que 
$\h(f)\leq \log d_t$. Consid\'erons $\Gamma_n$ le graphe de
l'application $(f,f^2,\ldots, f^{n-1})$ dans $U_{-n+1}\times
U_{-n+2}\times \cdots \times V$. C'est un sous-ensemble analytique de
dimension $k$ de $V^n$. De plus, $\Gamma_n$ est contenu dans
$U^{n-1}\times V$ et est un rev\^etement de degr\'e $d_t^n$ au dessus
de $V$. D'apr\`es le lemme 3.3.3, on a
$$\lov(f):=\lim_{n\rightarrow\infty} \frac{1}{n}\log \vol
(\Gamma_n\cap U^n)
\leq \log d_t.$$
Montrons que $\h(f)\leq \lov(f)$. Dans le cas des vari\'et\'es
compactes, ceci est le th\'eor\`eme de Yomdin-Gromov 
\cite{Gromov1,Gromov3,Yomdin}. Sa preuve est une application d'une
in\'egalit\'e de Lelong qui reste valide dans notre cas.
\par
Plus pr\'ecis\'ement, 
soit $\epsilon_0$ la distance entre $\K$ et le bord de $U$. Soit
$F\subset \K$ une famille $(n,\epsilon)$-s\'epar\'ee avec
$0<\epsilon<\epsilon_0$. Notons $F^*$ la famille des points
$a^*:=(a,f(a),\ldots, f^{n-1}(a))\in \Gamma_n$ avec $a\in F$ et
$B_{a^*}$ la boule de rayon $\epsilon/2$ centr\'ee en $a^*$. Ces
boules sont disjointes et contenues dans $U^n$. 
Une in\'egalit\'e de Lelong affime que $\vol(\Gamma_n\cap
B(a^*))\geq \pi^k(\epsilon/2)^{2k}$. Par cons\'equent, 
$$(\#F) \pi^k(\epsilon/2)^{2k}\leq \vol(\Gamma_n\cap U^n).$$
D'o\`u $\h(f)\leq \lov(f)$. 
\par \hfill $\square$
\subsection{Ensemble exceptionnel}
Soit $f:U\longrightarrow V$ une application d'allure polynomiale
et soit $X$ un sous-ensemble analytique propre de
$V$. 
Dans ce paragraphe, nous donnons une caract\'erisation du
sous-ensemble analytique maximal $\E_X$ de $X$ qui est
totalement invariant par $f$, \cad telle que $f^{-1}\E_X=\E_X\cap U$. 
Nous utiliserons ce r\'esultat pour
montrer, sous une
hypoth\`ese suppl\'ementaire,
l'analyticit\'e de l'ensemble exceptionnel $\E$. 
\par 
Nous montrons que $z\in\E_X$ si la proportion $\tau_X(z)$, 
de l'orbite n\'egative de $z$ dans $X$, est
strictement positive. 
Introduisons d'abord quelques notations. 
Pour tout $z\in V$, on pose 
\begin{eqnarray*}
\F^n_X(z): & = & \big\{w\in f^{-n}(z),\ f^i(w)\in X
\mbox{ pour tout } i=0,...,n \big\}\\
& = & f^{-1}(\F_X^{n-1}(z))\cap X \mbox{ si } n\geq 1
\end{eqnarray*}
$${\cal N}^n_X(z):=\#\F^n_X(z) \ \ \mbox{ et }\ \   
\tau_X(z):=\lim_{n\rightarrow \infty}\frac{{\cal N}^n_X(z)}{d_t^n}$$
o\`u $\#\F^n_X(z)$ est compt\'e avec multiplicit\'e.
Observons que la suite $d_t^{-n}{\cal N}^n_X(z)$ est d\'ecroissante,
en effet, 
$f(\F^{n+1}_X(z))\subset \F^n_X(z)$.
On pose
$$\E_X:=\big\{z\in V,\  \tau_X(z)>0\big\}.$$
On a bien $\E_X\subset X$.
\begin{theoreme} Soient $V$ une vari\'et\'e complexe,
$U$ un ouvert relativement compact de $V$ et $X$ un sous-ensemble
analytique de $V$.
Soit $f:U\longrightarrow V$ une application
holomorphe propre de degr\'e topologique $d_t\geq 1$. Alors
$\E_X$ est le plus grand sous-ensemble analytique
  de $X$ qui est totalement invariant par $f$, \cad
  $f^{-1}\E_X=\E_X\cap U$. En particulier,
  $$\E_X=\big\{z\in V,\ \tau_X(z)=1\big\}= 
\big\{z\in V,\ {\cal N}_X^n(z)=d_t^n
  \mbox{ pour tout } n\geq 0\big\}.$$
\end{theoreme}
\par
Posons pour tout $\theta$ r\'eel
$$\E_X(\theta):=\big\{z\in V, \ \tau_X(z)\geq\theta\big\}.$$
Soit $T$ un ensemble analytique irr\'eductible immerg\'e dans $V$. On
pose
\begin{eqnarray*}
{\cal N}^n_X(T):=\min_{z\in T}{\cal N}^n_X(z) & \mbox{ et } &
\tau_X(T):=\lim_{n\rightarrow \infty}\frac{{\cal N}^n_X(T)}{d_t^n}
\end{eqnarray*}
Observons que ${\cal N}^n_X(T)\leq {\cal N}^n_X(z)$ pour tout $z\in T$
et qu'on a \'egalit\'e en 
dehors d'un sous-ensemble analytique propre de $T$. Ceci se voit
ais\'ement sur le graphe de l'application $(f,f^2,\ldots,f^n)$.
\begin{lemme} Pour tout $\theta$ strictement positif, 
 $\E_X(\theta)$ est un
  sous-ensemble analytique de $X$.
\end{lemme}
\begin{preuve}
La fonction $z\mapsto {\cal N}_X^n(z)$ est semi-continue
sup\'erieurement au sens o\`u pour tout $\theta$ r\'eel, $\{z,\ 
{\cal N}_X^n(z)
\geq d_t^n\theta\}$ est un sous-ensemble
analytique de $V$. Par suite, pour tout $\theta$,
$\{\tau_X(z)\geq \theta\}$ est un sous-ensemble analytique de $V$
car $\tau_X$ est la limite d\'ecroissante des
$d_t^{-n}{\cal N}_X^n$. Si $\theta>0$, on a $\{\tau_X(z)\geq
\theta \}\subset X$. C'est donc un sous-ensemble analytique 
de $X$.
\end{preuve}
{\it D\'emonstration du th\'eor\`eme 3.4.1}---
L'ensemble $\E'_X$ des points $z$ tels que
  $f^{-n}(z)\subset X$ pour tout $n\geq 0$, est le plus grand
  sous-ensemble analytique totalement invariant de $f$ contenu dans
  $X$. Il est
clair que $\E'_X\subset \E_X$. Montrons
  que $\E_X=\E'_X$. Sinon 
d'apr\`es le lemme 3.4.2, on peut choisir une composante
  irr\'eductible $T$ de $\E_X$ qui v\'erifie
\begin{enumerate}
\item $T\not\subset \E'_X$.
\item $\theta_0:=\tau_X(T)$ est maximal parmi les $T$ qui v\'erifient
1. 
\end{enumerate}
La condition 2 est r\'ealisable. En effet,
si $x\in X$,
$\tau_X(x)$ est la moyenne des valeurs de $\tau_X$
sur les images r\'eciproques de $x$; si $x\not\in X$, $\tau_X(x)=0$. 
Pour tout $x\in V$ on a donc
$f_*\tau_X\geq d_t\tau_X$. 
Par cons\'equent, on peut se limiter \`a 
la famille
des composantes $T$ rencontrant $U$ et qui v\'erifient $\tau_X(T)\geq
\theta_0'$ avec un $\theta_0'>0$. C'est une famille non
vide quand $\theta'$ est assez petit; elle est finie car
$U\subset\subset V$.
\par 
Soit $x$ un point g\'en\'erique de $T$. On a 
$\tau_X(x)=\tau_X(T)$. Le nombre $\theta_0=\tau_X(T)$ \'etant maximal,
on a $\tau(x_i)\leq \tau(x)$ pour tout $x_i\in f^{-1}(x)$. D'autre part, 
on a $f_*\tau_X(x)\geq d_t \tau_X(x)$. Par suite,
$\tau(x_i)=\tau(x)$ et donc $f^{-1}(x)\subset X\cap
\E_X(\theta_0)$. On en d\'eduit que
$f^{-1}(T)\subset X\cap \E_X(\theta_0)$. Par r\'ecurrence, $f^{-n}(T)
\subset X$ pour tout $n\geq 0$. D'o\`u $T\subset \E'_X$.
\par
\hfill $\square$
\begin{remarque} \rm Le th\'eor\`eme 3.4.1 reste valide dans le cas o\`u
$U\subset V$ (sans supposer $U\subset\subset V$), 
$V$ est un ouvert d'une vari\'et\'e alg\'ebrique et o\`u
$f$ et $X$ sont alg\'ebriques. En effet, pour tout $\theta>0$, l'ensemble
$\E_X(\theta)$ d\'efini ci-dessus est alg\'ebrique et poss\`ede donc
un nombre fini de composantes, ce qui permet d'\'eviter l'hypoth\`ese
$U\subset\subset V$. De m\^eme, le
th\'eor\`eme 3.4.1 reste valide pour
les endomorphismes $f:V\longrightarrow V$ d'une vari\'et\'e compacte
$V$. On peut \'egalement \'etendre le
r\'esultat pour les applications
m\'eromorphes dominantes d'une vari\'et\'e
compacte $V$ dans elle-m\^eme.
\end{remarque}
\par
Le th\'eor\`eme suivant g\'en\'eralise un r\'esultat de
Briend-Duval \cite{BriendDuval2} pour les endomorphismes holomorphes
de $\P^k$.
Le fait que les applications $f$ ne
soient pas 
alg\'ebriques et donc qu'on n'a pas de
th\'eor\`eme de B\'ezout, nous oblige
\`a utiliser quelques arguments diff\'erents.
\begin{theoreme} Soit $f:U\longrightarrow V$ une application d'allure
  polynomiale d'ensemble critique $\Cr$.
Supposons que la suite de courants 
$$S_N:=\sum_{n=1}^N
  \frac{1}{d_t^n}(f^n)_*[\Cr\cap U_{-n}]$$ 
converge dans $V$ vers un courant $S$. Alors $\E$ est un
  sous-ensemble analytique de $V$, totalement invariant
  ($f^{-1}\E=\E\cap U$), contenu dans 
$\{a\in V, \nu(S,a)\geq 1\}$ o\`u $\nu(S,a)$ d\'esigne le nombre
de Lelong de $S$ en $a$. De plus, pour toute
  $\varphi$ continue, la suite de fonctions 
$d_t^{-n}(f^n)_*\varphi$ est localement
  \'equicontinue hors de $\E$.  
\end{theoreme}
\begin{remarques} \rm 
1. En g\'en\'eral, on n'a pas $\E\subset \Cr$. On a
  seulement, $\E\subset \PC_n$ pour $n$ assez grand. Il suffit pour le
  voir de consid\'erer l'exemple
$$f(z,w)=(w^d,2z)\mbox{ avec } d\geq 2.$$
L'application v\'erifie la condition du th\'eor\`eme. On a
$\Cr=\{w=0\}$ et $\E=\{zw=0\}=\PC_1$.
\par
2. L'hypoth\`ese du th\'eor\`eme 3.4.4 est satisfaite dans 
le cas de dimension 1 ou lorsque
l'ensemble critique $\Cr$ est pr\'ep\'eriodique. On verra dans 3.9 et
3.10 des grandes familles d'applications d'allure polynomiale
v\'erifiant cette hypoth\`ese.
\end{remarques}
\par
Pour simplifier les notations, on suppose que $V$ est un ouvert de
$\C^k$, $\Phi=|z|^2$ et $\omega=i\sum \d z_j\wedge \d\overline z_j$.
Pour le cas g\'en\'eral, il suffit de
recouvrir $V$ par des cartes biholomorphes \`a des ouverts de
$\C^k$. Pour tout $a\in V$, la famille des droites passant par
$a$ est  param\'etr\'ee par l'espace projective  $\P^{k-1}$ dont
la mesure invariante de masse 1 est not\'ee $\H_{2k-2}$.
Pour tout ensemble
connexe et simplement connexe
$X\subset V$, on appelle {\it bonne composante} 
de
$f^{-n}(X)$ toute composante connexe de $f^{-n}(X)$ qui ne
rencontre pas l'ensemble critique de $f^n$. En particulier, les
bonnes composantes s'envoient
bijectivement par $f^n$ sur $X$. En pratique, $X$ sera une
boule ou un disque holomorphe.
Si $\Delta$ est une droite passant par $a$,
on note $\Delta_r$ le disque de centre $a$ et de rayon $r$ dans
$\Delta$. On pose $\delta_n(\Delta_r):=d_t^{-n}
\#\big[f^n(\Cr\cap U_{-n})\cap \Delta_r\big]$ o\`u les valeurs critiques sont
compt\'ees avec multiplicit\'e.
\begin{lemme} Supposons qu'il existe $0<\nu<1$ tel que  
$\sum_{n\geq 1}\delta_n(\Delta_r)<\nu$.
Alors $f^{-n}(\Delta_{r/2})$ poss\`ede au moins
$(1-\sqrt{\nu})d_t^n$ bonnes
composantes de diam\`etre inf\'erieur \`a
$4\sqrt{\alpha_n(\Delta_r)/(\sqrt{\nu}-\nu)}$
o\`u $\alpha_n(\Delta_r):=d_t^{-n}\int_{\Delta_r}(f^n)_*\omega$.  
\end{lemme}
\begin{preuve} On montre d'abord par r\'ecurrence
que $f^{-n}(\Delta_r)$ poss\`ede
au moins $(1-\sum_1^n\delta_j(\Delta_r))d_t^n$ bonnes composantes.
Supposons le au rang $n-1$. Notons $\Delta^{n-1}$
l'union des $(1-\sum_1^{n-1}\delta_j(\Delta_r))d_t^{n-1}$ bonnes
composantes de $f^{-n+1}(\Delta_r)$. On a
$$\# f(\Cr\cap U)\cap\Delta^{n-1}\leq \# f^{n}(\Cr\cap
U_{-n}) \cap \Delta =\delta_n(\Delta_r)d_t^n.$$
Par cons\'equent, $f^{-1}(\Delta^{n-1})$ contient au moins
$$\left(1-\sum_1^{n-1}\delta_j(\Delta_r)\right)d_t^n-\delta_n(\Delta_r)d_t^n=
\left(1-\sum_1^{n}\delta_j(\Delta_r)\right)d_t^n$$ 
disques qui sont
des bonnes composantes de $f^{-n}(\Delta_r)$ (\ie les composantes ne
rencontrant pas $\Cr$). 
Ceci termine la r\'ecurrence. L'ensemble
$f^{-n}(\Delta_r)$ contient donc au moins $(1-\nu)d_t^n$ bonnes
composantes.
\par
Observons que $d_t^n\alpha_n(\Delta_r)$ est l'aire de $f^{-n}(\Delta_r)$.
Il y a donc au plus $(\sqrt{\nu}-\nu)d_t^n$
composantes de $f^{-n}(\Delta_r)$ dont l'aire est sup\'erieure \`a
$\alpha_n(\Delta_r)/(\sqrt{\nu}-\nu)$.
Par suite, $f^{-n}(\Delta_r)$
poss\`ede au moins $(1-\sqrt{\nu})d_t^n$ bonnes composantes
d'aire inf\'erieure \`a $\alpha_n(\Delta_r)/(\sqrt{\nu}-\nu)$, et qui
tend donc vers z\'ero.
\par
Soient $\Delta^n_i$ une telle composante et
$f_i^{-n}:\Delta_r\longrightarrow \Delta^n_i$ l'inverse de
$f^n$. La formule de Cauchy entra\^{\i}ne que
le diam\`etre de $f_i^{-n}(\Delta_{r/2})$ est inf\'erieur \`a
$4\sqrt{\alpha_n(\Delta_r)/(\sqrt{\nu}-\nu)}$. 
\par
Dans le cas o\`u
$V$ n'est pas un ouvert de $\C^k$, la famille des
applications holomorphes de $\Delta_r$ dans $\overline U$ est
normale.  
Il en r\'esulte que les valeurs adh\'erentes de la famille des $f^{-n}_i$
sont des applications constantes. On peut donc appliquer la formule de
Cauchy sur les cartes de $U$ et on a une majoration du diam\`etre de
$f^{-n}_i(\Delta_{r/2})$.
\end{preuve}
\begin{proposition} Soit $a\in V$. Supposons que
le nombre de Lelong
de $S$ en $a$ v\'erifie $0<\nu:=\nu(S,a)<1$. Alors, pour tout
$n\geq 0$ et toute boule
$B_r$ de centre $a$ et de rayon $r$ suffisamment petit, il existe
au moins $c_1(\nu) d_t^n$
bonnes composantes de $f^{-n}(B_{\delta_\nu r})$ de
diam\`etre inf\'erieur \`a $c_2(\nu)\sqrt{\alpha_n(B_r)}$
o\`u $\alpha_n(B_r):=d_t^{-n}\int_{B_r}
\omega^{k-1}\wedge(f^n)_*\omega$ et
$c_1$, $c_2$, $\delta_\nu$ sont strictement positifs. De plus, $c_1$,
$c_2$ d\'ependent contin\^ument
de $\nu$ et $\lim_{\nu\rightarrow 0}c_1(\nu)=1$. Si $\nu=0$ pour tout
$\epsilon>0$, pour $r$ suffisamment petit, il existe
$(1-\epsilon)d_t^n$ bonnes composantes de $f^{-n}(B_{\delta r})$ de
  diam\`etre inf\'erieur \`a $c_2\sqrt{\alpha_n(B_r)}$ o\`u $\delta$ et
  $c_2$ sont strictement positifs, d\'ependants de $\epsilon$. 
\end{proposition}
\begin{preuve} Rappelons que le nombre de Lelong de $S$ en $a$ est
  d\'efini par
$$\nu(S,a):=\lim_{r\rightarrow 0} \frac{1}{c_{k-1}r^{2k-2}}\int_{B_r}
S\wedge \omega^{k-1}$$
o\`u $c_{k-1}$ d\'esigne le volume de la boule unit\'e de $\C^{k-1}$. 
Posons $\delta_1:=\nu(1-\nu)/3$.
Si la boule $B_r$ est suffisamment petite, par d\'efinition du nombre
de Lelong,
la masse de $S$ dans $B_r$ est
inf\'erieure \`a $(\nu+\delta_1)c_{k-1} r^{2k-2}$. On note $\F'$ la
famille des droites $\Delta$, passant par $a$,
telles que la masse de la mesure $S\cap
\Delta$ dans $B_r$ soit inf\'erieure \`a $\nu':=\nu+2\delta_1$.
Par tranchage,
$$\H_{2k-2}(\F')\geq
1-\frac{\nu+\delta_1}{\nu+2\delta_1}=
\frac{\delta_1}{\nu+2\delta_1}=:2\delta'.$$
On note $\F$ la famille des droites $\Delta\in \F'$ telles que la masse de
$d_t^{-n}(f^n)_*\omega$ sur $\Delta\cap B_r$ soit inf\'erieure
\`a $\alpha_n':=\delta^{'-1}c_{k-1}^{-1}r^{-2k+2}\alpha_n(B_r)$. 
Par tranchage,
$\H_{2k-2}(\F)\geq 1-\delta'$.
\par
Fixons $\Delta$ dans $\F$. D'apr\`es le lemme 3.4.6,
$f^{-n}(\Delta_{r/2})$ poss\`ede au
moins $(1-\sqrt{\nu'})d_t^n$ bonnes composantes de diam\`etre
inf\'erieur \`a $l_n:=4\sqrt{\alpha_n'/(\sqrt{\nu'}-\nu')}$.
Notons $a_1,\ldots,a_m$ les points de $f^{-n}(a)$ et $\F_j$ la
famille des droites $\Delta\in \F$ telles que $f^{-n}(\Delta_{r/2})$
poss\`ede une bonne composante de diam\`etre inf\'erieur \`a $l_n$
passant par $a_j$. On a $m\leq d_t^n$, $\H_{2k-2}(\F_j)\leq
\H_{2k-2}(\F)$ et $\sum \H_{2k-2}(\F_j)\geq
(1-\sqrt{\nu'})d_t^n\H_{2k-2}(\F)$. On en d\'eduit facilement que
la famille suivante est de cardinal au moins
$(1-\sqrt[3]{\nu'})d_t^n$:
$$\tilde A_n:=\big\{a_j,\ \H_{2k-2}(\F_j)>
(1-\sqrt[6]{\nu'}) \H_{2k-2}(\F)\big\}.$$
La preuve de la proposition est compl\'et\'ee gr\^ace au lemme suivant
appliqu\'e au fonctions coordonn\'ees des inverses de $f^n$.
(Le cas $\nu=0$ se traite de fa\c con semblable en prenant $\delta_1>0$
assez petit).
\begin{lemme}[\cite{Alexander}, \cite{SibonyWong}] 
Soit $\F_j$ une famille de droites passant
par $a$. 
Supposons que $\H_{2k-2}(\F_j)\geq 2\delta_\nu>0$.
Alors toute fonction $g$ holomorphe au voisinage de 
$\F_j\cap B_{r/2}$ se prolonge en fonction holomorphe dans
la boule $B_{\delta_\nu r}$. De plus, 
$$\sup_{B_{\delta_\nu r}}|g|\leq \sup_{\F_j\cap B_{r/2}}|g|.$$
\end{lemme}
\end{preuve}
\begin{lemme} 
Soit $0<\nu<1$. 
Supposons $a\in V$ avec $\nu(S,a)\leq \nu$. Soit $\mu^a$ une 
valeur d'adh\'erence de la
  suite $(\mu^a_n)$. Alors $\mu^a$ est \'egale \`a la somme de deux
  mesures positives $\mu^a_r$ et $\mu^a_s$ o\`u $\mu^a_r$ est
  de masse au moins $c_1(\nu)$ et
est absolument continue par rapport \`a $\mu$.
En particulier, on a $\E\subset \PC_\infty$.
\end{lemme}
\begin{preuve} Soit $\mu^a=\mu^a_r+\mu^a_s$ la d\'ecomposition de
  Lebesgue de $\mu^a$ avec $\mu^a_r<<\mu$. Soit
$\psi$ une fonction test de classe ${\cal C}^1$.
Soit $B_r$ une boule v\'erifiant les propri\'et\'es de la proposition 3.4.7.
Posons $B:=B_{\delta_\nu r}$ et soit $B^n$ la
  famille des bonnes composantes de $f^{-n}(B)$.
Notons $B^n_i$ 
ces composantes, $f^{-n}_i:B\longrightarrow B^n_i$ les
inverses de $f^n$ pour $i=1$, $\ldots$, $d_n$ et $\tilde
\mu^z_n:=\sum \delta_{z_i^n}$ o\`u $z_i^n:=f^{-1}(z)\cap B_i^n$. 
Fixons une suite croissante $(n_i)$ telle que $\mu^a_{n_i}$ tende vers
  $\mu^a$ et $\tilde \mu^a_{n_i}$ tende vers
une mesure $\tilde \mu^a$.
Le cardinal $d_n$ de $B^n$
  v\'erifiant $c_1(\nu)d_t^n\leq d_n\leq d_t^n$, la masse de
$\tilde \mu^a_{n_i}$ est plus grande ou \'egale \`a $c_1(\nu)$.  
La famille des applications $f^{-n}_i$ est \'equicontinue car $U$ est
Kobayashi hyperbolique.
Posons pour tout
  $z\in B$ 
$$\tilde \psi_n(z):=\frac{1}{d_t^n}\sum_{i=1}^{d_n}
\psi(f^{-n}_i(z)).$$
La fonction $\psi$ \'etant de classe ${\cal C}^1$ et la famille des
applications $f^{-n}_i$ \'etant \'equicontinue, il existe une constante $c>0$
telle que
$|\psi(f^{-n}_i(a))-\psi(f^{-n}_i(z))|\leq c|a-z|$ pour tous $i$ et
$n$. 
Par cons\'equent, 
$$|\tilde\psi_n(a)-\tilde\psi_n(z)|\leq c|a-z|.$$
Pour toute valeur d'adh\'erence $\tilde \mu^z$ de la suite 
$(\tilde \mu^z_{n_i})$ on a $|\int\psi \d\tilde\mu^a-\int\psi
\d\tilde\mu^z|\leq c|a-z|$. En particulier,
$\tilde\mu^z\rightharpoonup \mu^a$ quand $z\rightarrow a$.
Si $z\not \in \E^{(n_i)}$, on peut d'apr\`es la proposition 3.2.5,
supposer que $\mu_{n_i}^z\rightharpoonup \mu$ et donc 
la mesure $\mu-\tilde \mu^z$
est positive. En prenant $z\rightarrow a$ et
$z\in B\setminus \E^{(n_i)}$, on
a que $\mu-\tilde \mu^a$ est positive. Ceci
implique que la masse de $\mu^a_r$ (qui est minor\'ee par la masse de
$\tilde\mu^a$) est sup\'erieure ou \'egale \`a
$c_1(\nu)$.
\par
Si $a\not \in \PC_\infty$, en utilisant la formule de Poisson-Jenssen,
on v\'erifie facilement que $\nu(S,a)=0$. 
Utilisant le m\^eme raisonnement et la proposition 3.4.7, on montre
que 
la masse de $\mu_r^a$ est \'egale
\`a 1.
Par cons\'equent, 
$\mu^a_r=\mu$ et donc $a\not\in \E$. 
\end{preuve}
\begin{lemme} Soit $X_\nu$ l'ensemble
(analytique) des points $z\in V$
tels que
$\nu(S,z)\geq \nu$ avec $0<\nu<1$.
Alors pour $a\not\in \E_{X_\nu}$, toute valeur
d'adh\'erence $\mu^a$ de la suite $(\mu^a_n)$ est \'egale \`a la
somme de deux mesures positives
$\mu^a_r$ et $\mu^a_s$ o\`u $\mu^a_r$ est
une mesure de masse au moins $c_1(\nu)$ et est absolument continue
par rapport \`a $\mu$.
\end{lemme}
\begin{preuve} Le fait que $X_\nu$ soit analytique r\'esulte du
  th\'eor\`eme de Siu \cite{Siu}. Dans notre cas, on peut montrer
  facilement qu'il existe $n_0$ assez grand tel que
$X_\nu\subset \PC_{n_0}$.
Soit $\mu^a_r$ la partie r\'eguli\`ere
de $\mu^a$ par rapport \`a $\mu$.  
Fixons un $\epsilon>0$. Il suffit de montrer que la masse de
$\mu^a_r$ est plus grande ou \'egale \`a $(1-\epsilon)c_1(\nu)$. 
D'apr\`es le lemme 3.4.9, il reste \`a traiter le cas $a\in X_\nu$.
Puisque $a\not\in
\E_{X_\nu}$, d'apr\`es le th\'eor\`eme 3.4.1, on a
  $\tau_{X_\nu}(a)=0$. Observons que 
\begin{eqnarray*}
1-\tau_{X_\nu}(a)=1-\lim\frac{\#\F_{X_\nu}^n(a)}{d_t^n}
& = & 
\sum_{n=1}^\infty\left[\frac{\#\F_{X_\nu}^{n-1}(a)}{d_t^{n-1}}-
\frac{\#\F_{X_\nu}^n(a)}{d_t^n}\right]\\
& = & \sum_{n=1}^\infty \frac{1}{d_t^n} \#
\left[f^{-1}(\F_{X_\nu}^{n-1}(a))\setminus X_\nu\right].
\end{eqnarray*}
Dans la suite, on consid\`ere une somme partielle de la derni\`ere 
s\'erie de
gauche. Il existe un entier $N_0$, 
des entiers positifs $n_j\leq N_0$,
des points $a_j\in f^{-n_j}(a)\setminus X_\nu$ (les $n_j$ ne sont pas
\`a priori distincts)
 et $b_i\in f^{-N_0}(a)\cap X_\nu$ v\'erifiant
\begin{enumerate}
\item $1-\epsilon\leq \sum d_t^{-n_j}\leq 1$.
\item $\mu^a_{N_0}=d_t^{-N_0}
\big[\sum (f^{N_0-n_j})^*\delta_{a_j}+
\sum \delta_{b_i}\big]$.
\end{enumerate}
On a donc
$$\mu^a_{N}=\frac{1}{d_t^{N}}
\left [\sum d_t^{N-n_j}\mu^{a_j}_{N-n_j}+
\sum d_t^{N-N_0}\mu_{N-N_0}^{b_i}\right ].$$
On choisit une suite croissante $(N_i)$ telle que
$\mu_{N_i-n_j}^{a_j}$ (resp. $\mu_{N_i-n_j}^a$)
converge vers une mesure $\mu^{a_j}$ (resp. vers $\mu^a$) quand
$i\rightarrow\infty$.  
D'apr\`es le lemme 3.4.9, puisque $\nu(S,a_j)<\nu$, on a
$\mu^{a_j}=\mu^{a_j}_r+\mu^{a_j}_s$ avec
$\mu^{a_j}_r$ absolument continue
par rapport \`a $\mu$ et de masse au moins $c_1(\nu)$.
Soit $\mu^a_{r,\epsilon}:=\sum d_j^{-n_j}\mu^{a_j}_r$. Alors
$\mu^a_{r,\epsilon}$
est absolument continue par rapport \`a $\mu$
et sa masse est plus grande ou \'egale \`a
$(1-\epsilon)c_1(\nu)$. La masse de $\mu^a_r$ (qui est minor\'ee
par celle de $\mu^a_{r,\epsilon}$) est donc plus grande ou
\'egale \`a $(1-\epsilon)c_1(\nu)$.
\end{preuve}
{\it Fin de la d\'emonstration du th\'eor\`eme 3.4.4}--- Posons
$\E_n:=\E_{\PC_n}$. Observons que
si $a\in\E_n$ et $\mu^a$
est une valeur adh\'erente \`a la suite $(\mu^a_n)$ alors
$\mu^a(\E_n)=1$ et donc $\mu^a\not=\mu$ car $\mu(\PC_\infty)=0$; si
$a\not\in \bigcup \E_n$, on peut appliquer le lemme 3.4.10 avec
$\nu\rightarrow 0$. Donc
$\E=\bigcup_{n\geq 1} \E_n$. On en d\'eduit que
$f^{-1}(\E)=\E\cap U$. 
\par
Soit $X_\nu:=\{a\in V,\ \nu(S,a)\geq\nu\}$. Alors pour tout
$\nu>0$, $X_\nu$ est un sous-ensemble analytique de $V$ contenu
dans $\PC_\infty$. Donc $\E_{X_\nu}\subset \E=\bigcup_{n\geq
1}\E_n$. Montrons que $\E_{X_1}=\E$.
Sinon, il existe $n\geq 1$ et 
$a\in\E_n\setminus \E_{X_1}$ avec $\nu(S,a)<1$.
Soit $\mu^a$
une valeur d'adh\'erence de la suite $\{\mu_m^a\}$. On a
$\mu^a(\E_n)=1$, car $\E_n$ est un ferm\'e totalement invariant
et donc $\mu^a(\PC_n)=1$. D'apr\`es le lemme 3.4.10, on a
$\mu^a_r(\PC_n)\geq c_1(\nu)>0$.
La mesure $\mu^a_r$ \'etant absolument
continue par rapport \`a $\mu$, on a $\mu(\PC_n)>0$.
C'est la contradiction cherch\'ee car
$\mu$ ne charge pas les ensembles postcritiques, puisque $\log J$ est
$\mu$-int\'egrable.
\par
\hfill $\square$
\par
%
%\begin{corollaire} Soit $f:U\longrightarrow V$ une application
%  d'allure polynomiale. Si $\mu$ est une mesure PLB 
%
%alors l'ensemble
%  exceptionnel $\E$ est un ensemble analytique. Il en est de m\^eme si
%  $\delta:=\lim\sqrt[n]{\delta_n}<1$.
%
%\end{corollaire}
%
%\begin{preuve} Cela r\'esulte des estimations de la proposition 3.3.5
%  et du th\'eor\`eme 3.3.4.
%\end{preuve}
%
\subsection{Points p\'eriodiques r\'epulsifs}
Sous la m\^eme hypoth\`ese qu'au paragraphe pr\'ec\'edent, nous
allons montrer la densit\'e des points p\'eriodiques r\'epulsifs
dans le support de la mesure d'\'equilibre.
\begin{theoreme} Soit $f$ une application v\'erifiant les
  hypoth\`eses du th\'eor\`eme 3.4.4. Notons $\PR_n$ l'ensemble des
  points p\'eriodiques r\'epulsifs d'ordre $n$ de $f$ qui
  appartiennent \`a $\supp(\mu)$. Supposons  
$\limsup d_t^{-n}\#\PR_n\leq 1$. Alors les points p\'eriodiques
r\'epulsifs de $f$ sont \'equidistribu\'es sur le support de $\mu$,
\cad que la suite de mesures
$$\nu_n:=\frac{1}{d_t^n}\sum_{a\in\PR_n}\delta_a$$
converge vers $\mu$. 
\end{theoreme}
\begin{preuve} La d\'emonstration reprend des id\'ees d\'ej\`a
  utilis\'ees par Briend-Duval \cite{BriendDuval1}
pour d\'emontrer le m\^eme r\'esultat
  pour les endomorphismes de $\P^k$. 
On pose $X:=V\setminus \bigcup_{n\geq 0}f^{-n}(\PC_\infty)$. 
Puisque $\mu$ ne charge pas
l'ensemble postcritique, $X$ est de mesure totale. On consid\`ere
$$\widehat{X}:=\big\{\widehat{x}=(x_{-n})_{n\geq 0},\
f(x_{-n})=x_{-n+1}\big\}$$
l'espace des pr\'ehistoires. La mesure $\mu$ et l'application $f$  
se remontent en mesure
$\widehat{\mu}$ sur $\widehat{X}$ et une application 
$\widehat{f}$. La mesure $\widehat\mu$ est m\'elangeante pour $\widehat f$. 
Notons $\pi$ la projection
canonique de $\widehat{X}$ dans $X$, \ie $\pi((x_{-n})):=x_0$ et
$\sigma:=\widehat{f}^{-1}$ le d\'ecalage \`a gauche de $\widehat{X}$,
\ie $\sigma(\widehat{x})=(x_{-n})_{n\geq 1}$. Notons
$f^{-n}_{\widehat{x}}$ la branche inverse de $f^n$ v\'erifiant
$f^{-n}_{\widehat{x}}(x_0)=x_{-n}$. 
Pour tous $\delta>0$ et $c>0$, on pose
\begin{eqnarray*}
E_{\delta,c} & := & \big\{\widehat{x}, \ f^{-n}_{\widehat{x}} \mbox{ est
  d\'efinie sur } B(x_0,2\delta),\\
& &  \diam \big(f^{-n}_{\widehat{x}}(B(x_0,2\delta))\big)\leq 
c\sqrt{d_{k-1,n}d_t^{-n}}\big\}.
\end{eqnarray*}
D'apr\`es la proposition 3.4.7,
  $\bigcup_{\delta,c}E_{\delta,c}=\widehat{X}$. 
\par
Pour tout bor\'elien $B$ de $V$, on note $\widehat{B}:=\pi^{-1}(B)$ et
$\widehat{B}_{\delta,c}:=\widehat{B}\cap E_{\delta,c}$. Posons
$\mu_{\delta,c}:=\pi_*(\widehat{\mu}_{|E_{\delta,c}})$. On a
$\mu_{\delta,c}(B)=\widehat{\mu}(B_{\delta,c})$. 
\par
Par hypoth\`ese, 
toute valeur d'adh\'erence $\nu$ de $(\nu_n)$ est de masse au
plus 1. Pour montrer que $\nu=\mu$, il suffit de comparer $\nu$ et
$\mu_{\delta,c}$. Soient $x\in \pi(E_{\delta,c})$ 
et $r>0$, $\epsilon>0$ assez petits tels
que $r+\epsilon<\delta$. Dans la suite, on consid\`ere uniquement les
boules ferm\'ees.
Il nous suffit de montrer que
$$(1-\epsilon)\mu_{\delta,c}(B(x,r))\leq
\nu(B(x,r+\epsilon)).$$
\par 
La mesure $\widehat\mu$ \'etant m\'elangeante, on
a pour $n$ assez grand:
\begin{eqnarray*}
(1-\epsilon)\mu_{\delta,c}
(B(x,r))\mu(B(x,r)) & = &
(1-\epsilon)\widehat{\mu}
(\widehat{B}_{\delta,c}(x,r))\mu(B(x,r)) \\
 & \leq & \widehat{\mu}(\sigma^n(\widehat{B}_{\delta,c}(x,r))\cap
 \widehat{B}(x,r))\\
& \leq & \mu(\pi(\sigma^n(\widehat{B}_{\delta,c}(x,r)))\cap B(x,r))
\end{eqnarray*}
Observons que pour tout $\widehat{x}\in E_{\delta,c}$,
$f_{\widehat{x}}^{-n}(B(x,r+\epsilon))$ est de diam\`etre
inf\'erieur \`a $\epsilon$ d\`es que $n\geq n_0$ uniform\'ement en
$\widehat{x}$. De plus, si cette composante rencontre $B(x,r)$, 
elle est contenue dans $B(x,r+\epsilon)$. D'apr\`es le th\'eor\`eme du
point fixe,
elle contient exactement un point p\'eriodique r\'epulsif dont
la p\'eriode divise $n$. Notons $\F(x,r)$ la famille de ces
composantes. On a
\begin{eqnarray*}
(1-\epsilon)\mu_{\delta,c}
(B(x,r))\mu(B(x,r))
& \leq & \mu(\bigcup B)\mbox{ avec } B\in\F(x,r)\\
& \leq & \frac{1}{d_t^n}\mu(B(x,r+\epsilon))\#\F(x,r)\\
& \leq & \nu_n(B(x,r+\epsilon))\mu(B(x,r+\epsilon)).
\end{eqnarray*}
Faisant tendre $\epsilon$ vers $0$, on obtient $\nu\geq
\mu_{\delta,c}$ et donc $\lim \nu_n=\mu$. 
\par
Nous terminons la preuve par la remarque
suivante. Il suffit de consid\'erer les boules $B(x,r+\epsilon)$
rencontrant $J_k$. Puisque $J_k$ est totalement invariant, toute
composante de l'image r\'eciproque de $B(x,r+\epsilon)$ rencontre
$J_k$. Par suite, les points p\'eriodiques r\'epulsifs obtenus
ci-dessus sont dans $J_k$. 
\end{preuve}
\begin{remarque} \rm Si on ne suppose pas $\limsup d_t^{-n}\#\PR_n\leq
  1$, on peut montrer qu'il existe des ensembles $\PR_n'$ de points
  p\'eriodiques r\'epulsifs dont l'ordre divise $n$ tels que
la suite de mesures
$$\nu_n':=\frac{1}{d_t^n}\sum_{a\in\PR_n'}\delta_a$$
converge vers $\mu$. 
\end{remarque}
\par
Donnons des exemples o\`u on peut majorer le nombre de points
p\'eriodiques d'ordre $n$. Soit 
$f:\C^k\longrightarrow \C^k$ une application polynomiale propre. 
On dit que {\it l'infini est
  attirant}
si pour tout $R>0$ assez grand, $|f^n(z)|\rightarrow
  \infty$ uniform\'ement sur $\C^k\setminus B(0,R)$ o\`u $B(0,R)$
  d\'esigne la boule 
de rayon $R$ centr\'ee
  en $0$. On d\'efinit {\it l'ensemble de Julia rempli}
par 
$$\K:=\big\{z\in\C^k,\ \big(f^n(z)\big)_{n\geq 0} \mbox{ born\'ee}\big\}.$$
\begin{proposition} Soit $f:U\longrightarrow V$ une application d'allure
  polynomiale, $V$ \'etant un ouvert de $\C^k$. Supposons que le
  diam\`etre $\diam(\K)$ de $\K$ est strictement plus petit que
  $\max_{a\in \K}\dist(a,\partial U)$. Alors pour $n$ assez grand, 
$f$ poss\`ede exactement $d_t^n$ points
p\'eriodiques (compt\'es avec multiplicit\'e) dont la p\'eriode divise
  $n$. En particulier, ceci est vrai pour toute application
  polynomiale propre $f:\C^k\longrightarrow \C^k$ telle que l'infini soit
  attirant et $\K\not=\emptyset$.
\end{proposition}
\begin{preuve} Soit $b\in \K$ tel que $\dist(b,\partial U)=\max_{a\in
  \K} \dist(a,\partial U)$. On a pour $\epsilon>0$ assez petit et 
pour $z\in \partial U_{-n-1}$ avec $n$ assez
  grand
$$|f^n(z)-z-(f^n(z)-b)|=|z-b|\leq\diam(\K)+\epsilon<|f^n(z)-b|$$
car $f^n(z)\in\partial U$. On peut donc appliquer le th\'eor\`eme de
Rouch\'e: le nombre de solutions de $f^n(z)=z$ est \'egal au nombre de
solutions de $f^n(z)=b$.
\par
Si $f:\C^k\longrightarrow \C^k$ est une application avec l'infini
  attirant et $\K\not=\emptyset$, 
on choisit $U$ une boule assez grande et le raisonnement
  ci-dessus est vrai pour tout $n$ assez grand.
\end{preuve}
\begin{proposition} Soit $f:U\longrightarrow V$ une application
  d'allure polynomiale. Supposons que $U$ est
  holomorphiquement contractible dans $V$,
\cad qu'il existe $h:[0,1]\times
  \overline U\longrightarrow V$ continue, $h(t,.)$
holomorphe sur $U$ avec $h(0,z)=p$ et $h(1,z)=z$ pour tout $z\in
\overline U$. 
Alors le nombre de
points p\'eriodiques dont la p\'eriode divise $n$ est \'egale \`a
$d_t^n$ pour tout $n\geq 1$.
\end{proposition}
\begin{preuve}
Observons que le nombre de solutions de l'\'equation
  $f^n(z)=h(t,z)$ avec $z\in U_{-n}$ ne d\'epend pas de $t\in
[0,1]$. Pour $t=0$, le nombre de solutions est \'egal \`a $d_t^n$.   
\end{preuve}
\subsection{Exposants de Lyapounov}
Sous les hypoth\`eses des paragraphes 3.4 et 3.5, 
nous avons une minoration suivante des exposants de Lyapounov de
$f$ qui g\'en\'eralise le r\'esultat de Briend-Duval
\cite{BriendDuval1} pour les endomorphismes de $\P^k$. 
\begin{theoreme} Soit $f$ une application
v\'erifiant les hypoth\`eses
  du th\'eor\`eme 3.4.4. 
%Supposons $\alpha_n:=\int_U
%\omega^{k-1}\wedge \frac{f^n_*\omega}{d_t^n}$ 
%d\'ecroit exponentiellement (\cad que
%$\alpha_n$ est major\'e par $c\lambda^{-n}$ pour un $c>0$ et un
%$\lambda>1$), a
Alors les exposants de Lyapounov de $f$ sont
sup\'erieurs ou \'egaux \`a $\frac{1}{2}\log (d_t/d_{k-1})$.
%$\frac{1}{2}\log \lambda$.
En particulier, ils sont tous positifs
ou nuls et ils sont strictement positifs si $d_{k-1}<d_t$. 
\end{theoreme}
\begin{preuve}
Fixons $\lambda$ tel que $0<\lambda< d_t/d_{k-1}$.
Il existe une constante $c>0$
telle que $d_{k-1,n}\leq c\lambda^{-n} d_t^n$ pour tout $n\geq 1$. 
On peut supposer $V\subset\C^n$. Pour le cas
g\'en\'eral, il suffit de recouvrir $\K$ par une famille finie
d'ouverts biholomorphes \`a la boule unit\'e de $\C^k$. On a vu
dans la proposition 2.3.4 que les exposants de Lyapounov existent
et sont constants. Le plus petit exposant de Lyapounov est
calcul\'e par la formule
$$\lambda_{\min} :=-\lim_{n\rightarrow\infty}\frac{1}{n} \int\log
\|(\Dr f^n)^{-1}\| \d\mu.$$
D'apr\`es la proposition 3.4.7, on peut choisir une famille finie
de boules disjointes $B_1$, $\ldots$, $B_m$ et 
des nombres r\'eels positifs
assez petit $\epsilon$, 
$\delta$ v\'erifiant les
propri\'et\'es suivantes:
\begin{enumerate}
\item $\mu(B_1\cup\ldots\cup B_m)\geq 1-\delta/2$;
\item $\mu(\partial B_i)=0$;
\item Il existe un $c'>0$ tel que pour tout $n\geq 1$,
chaque $B_i$ admette au moins $(1-\epsilon)d_t^n$ branches
inverses $f^{-n}_{i,j}$ de $f^n$ de diam\`etre inf\'erieur \`a
$c'\lambda^{-n/2}$. Notons $B^{-n}_{i,j}$ les images de ces branches.
\end{enumerate}
On choisit les boules $\tilde B_i\subset\subset B_i$ telles que
$\mu(\tilde B_1\cup\ldots\cup \tilde B_m)\geq 1-\delta$.
Posons $\tilde
B^{-n}_{i,j}:= f^{-n}(\tilde B_i)\cap B^{-n}_{i,j}$. 
Puisque $f^n$ envoie injectivement la
composante $B_{i,j}^{-n}$ qui est de diam\`etre inf\'erieur \`a
$c'\lambda^{-n/2}$ dans $B_i$, on a que 
$\|(\Dr f^n)^{-1}\|\leq \tilde c\lambda^{-n/2}$ sur $\tilde
B^{-n}_{i,j}$ avec un $\tilde c>0$. Soit
$M>0$ un majorant des valeurs propres de $\Dr f$ sur $\K$.
On en d\'eduit
une minoration de la plus petite valeur propre de $\Dr f$ sur $\K$:
$$\|(\Dr f^n)^{-1}\|^{-1}\geq \det (\Dr
f^n) M^{(-k+1)n}=\sqrt{J_{f^n}}M^{(-k+1)n}$$
o\`u $J_{f^n}$ est le jacobien r\'eel de $f^n$. Posons
$E_n:=V\setminus\bigcup \tilde B_{i,j}^{-n}$. 
On sait que 
$\mu(\bigcup_{i,j}\tilde B_{i,j}^{-n})
\geq (1-\epsilon)(1-\delta)$. D'o\`u
$\mu(E_n)\leq 1-(1-\epsilon)(1-\delta)$ et donc 
$\mu(f^{-s}(E_n))\leq  1-(1-\epsilon)(1-\delta)$ pour $s\geq 1$. 
Comme $-\int\log J\d\mu<\infty$, pour tout $\epsilon'>0$, on a 
$-\int_E\log J\d\mu<\epsilon'$ lorsque $\delta$, $\epsilon$ sont
assez petits et $\mu(E)\leq 1-(1-\epsilon)(1-\delta)$. On a
les estimations suivantes en utilisant l'in\'egalit\'e sur $\bigcup
\tilde B^{-n}_{i,j}$:
\begin{eqnarray*}
\lefteqn{-\frac{1}{n}\int \log\|(\Dr f^n)^{-1}\|\d \mu \geq
\frac{1}{n}\left[\frac{n}{2}\log \lambda -\log \tilde c\right]
\mu\left(\bigcup_{i,j} \tilde B^{-n}_{i,j}\right)+}\\
& & 
+\int_{E_n} \left[\frac{1}{2n}\log J_{f^n}-(k-1)\log M\right]\d\mu
\\
&\hspace{2cm} \geq & \left[\frac{1}{2}\log \lambda
-\frac{1}{n} \log \tilde
c\right](1-\epsilon)(1-\delta)+ \\
& & 
+\int_{E_n} \left(\frac{1}{2n}\sum_{s=0}^{n-1}\log J\circ f^s\right) \d\mu
-(k-1)\log M\mu(E_n)\\
&\hspace{2cm} = & (1-\epsilon)(1-\delta)\left[\frac{1}{2} \log\lambda
-\frac{1}{n}\log\tilde c\right]+\\
& & 
+\frac{1}{2n}\sum_{s=0}^{n-1}\int_{f^{-s}(E_n)}\log J\d\mu
-(k-1)\log M\mu(E_n)\\
&\hspace{2cm} \geq &
(1-\epsilon)(1-\delta)\left[\frac{1}{2} \log\lambda
-\frac{1}{n}\log\tilde c\right]-\\
& & 
-\frac{\epsilon'}{2}-(k-1)\log M \big[1-(1-\epsilon)(1-\delta)\big]
\end{eqnarray*}
et donc 
$$\lambda_{\min}\geq \frac{(1-\epsilon)(1-\delta)}{2}\log \lambda
-\frac{\epsilon'}{2}-(k-1)\log M \big[1-(1-\epsilon)(1-\delta)\big].$$
Faisant tendre $\epsilon'$, $\epsilon$ et $\delta$ vers $0$, on obtient
$\lambda_{\min}\geq \frac{1}{2}\log \lambda$. On en d\'eduit que
$\lambda_{\min}\geq \frac{1}{2}\log (d_t/d_{k-1})$. D'apr\`es la
proposition 3.3.1, on a $\lambda_{\min}\geq 0$.
\end{preuve}
\subsection{Cas de dimension 1}
Soit
$f:U\longrightarrow V$ une application \`a allure polynomiale de
degr\'e $d_t\geq 2$ o\`u $V$
est une surface de Riemann ouverte et $U\subset\subset V$ est un
ouvert de $V$. Lorsque $U$ et $V$ sont simplement connexes, Douady-Hubbard
\cite{DouadyHubbard} ont montr\'e que $f$ est conjugu\'ee \`a un
polyn\^ome de degr\'e $d_t$ par un hom\'eomorphisme h\"old\'erien. 
Ind\'ependement de cette approche, nous
allons donner quelques r\'esultats de nature m\'etrique sur la mesure
d'\'equilibre $\mu$.
Dans un tr\`es joli article, \cite{Mane}, Ma\~ne a \'etudi\'e les
propri\'et\'es m\'etriques de la mesure $\mu$ pour les applications
rationnelles de $\P^1$. La formule donn\'ee au point 4 du th\'eor\`eme
ci-dessus se trouve dans \cite{Mane} pour les applications
rationnelles. Elle est d\^ue \`a Manning \cite{Manning} pour les
polyn\^omes \`a une variable.
\par
Posons $M_n:=\sqrt[n]{\sup_\K |(f^n)'(z)|}$. 
Les $M_n$ d\'ependent de
la m\'etrique choisie pour $V$. La suite $M_n$ d\'ecroit vers une
constante $M>0$. On v\'erifie facilement que $M$ ne d\'epend pas de la
m\'etrique.
\begin{theoreme} Soit $f:U\longrightarrow V$ 
une application d'allure polynomiale de
degr\'e $d_t\geq 2$ o\`u 
$V$ est une surface de Riemann ouverte et $U$ est un ouvert relativement
compact dans $V$. On note $\mu$ la mesure d'\'equilibre et $\K$
l'ensemble de Julia rempli associ\'es 
\`a $f$. Soit $\alpha$ un nombre r\'eel tel que
$0<\alpha<\log d_t/\log M$. Alors
\begin{enumerate}
\item Tout potentiel $G$ de la mesure $\mu$ est h\"old\'erien d'ordre
  $\alpha$. 
\item Pour tout disque $B(x,r)$, on a $\mu(B(x,r))\leq cr^{\alpha}$
  o\`u $c>0$ est une constante.
\item L'entropie de $f$ est \'egale \`a $\log d_t=\h_\mu(f)$.
\item Si $\HD(\mu)$ d\'esigne la dimension de Hausdorff 
de $\mu$ on a
$$\frac{1}{\HD(\mu)}=\frac{\int\log|f'|\d\mu}{\h_\mu(f)} = 
\frac{\int\log|f'|\d\mu}{\log d_t} \geq 1.$$
\end{enumerate}
En particulier, la mesure $\mu$ v\'erifie les hypoth\`eses des
des th\'eor\`emes 3.4.4, 3.5.1, 3.6.1; elle est
approximable par les points p\'eriodiques r\'epulsifs.
\end{theoreme} 
\begin{preuve} 1. Quitte \`a remplacer $f$ par $f^n$ et $U$, $V$ par
  $U_{-m-n-1}$, $U_{-m-n}$ pour $n$ et $m$ assez grands, on peut
  supposer que $\alpha\leq \alpha_0:=\log d_t/\log M_1$.
Soit $G$ un potentiel de $\mu$ dans $V$. C'est une fonction
  sousharmonique dans $V$ et harmonique en dehors du support de
  $\mu$. En particulier, elle est harmonique sur $V\setminus \K$; elle
  est donc born\'ee sur $U\setminus U_{-2}$. 
Observons
  que $d_t^{-1} G\circ f$ est un potentiel de $d_t^{-1}
  f^*\mu=\mu$ dans $U$. Par cons\'equent, 
il existe une fonction harmonique $u$ dans $U$
telle que pour $w\in U$ on ait
$$G(w)-\frac{G(f(w))}{d_t}=u(w).$$
On en d\'eduit par it\'eration  que sur $U_{-n-1}$ 
\begin{eqnarray}
G(w)-\frac{G(f^n(w))}{d^n_t} & = & \sum_{j=0}^{n-1}
\frac{u(f^j(w))}{d_t^j}.
\end{eqnarray}
Soit $A>0$ une constante telle que $|G|<A$ sur $U \setminus U_2$ et 
$|u|<A$ sur $U_{-2}$.
On a  $|G(f^n(w))|<A$ lorsque
$w\in U_{-n-2}\setminus U_{-n-3}$.
Il r\'esulte de la relation (5)
que pour $w\in U_{-n-2}\setminus U_{-n-3}$ on a
$$|G(w)|\leq \frac{A}{d_t^n}+ A\sum_{j=0}^{n-1} \frac{1}{d_t^j} \leq
3A.$$
La fonction $G$ est donc born\'ee dans $U\setminus \K$. Par
semi-continuit\'e, cela entra\^{\i}ne
que $G$ est born\'ee dans $\overline{U\setminus \K}$.
Pour $w\in \partial\K$, en passant \`a la limite, on obtient 
$$G(w)=\sum_{j=0}^\infty \frac{u(f^j(w))}{d_t^j}.$$
Donc $G_{|\partial\K}$ est continue. D'apr\`es \cite[p.53]{Tsuji}, 
$G$ est continue sur $V$.
\par 
Montrons que $G$ est h\"old\'erienne d'ordre $\alpha_0$. On a
$M_1=d_t^{1/\alpha_0}$. Soient $w$ et
$w'$ deux points suffisamment proches de $\K$.
Supposons que $w\in U_{-n}$ 
et $w'\in U_{-m}$ avec $m\geq n$. 
On a pour tout $0\leq s\leq n$ 
$$G(w)-G(w')= \frac{G(f^s(w))-G(f^s(w'))}{d_t^s}+ 
\sum_{j=0}^{s-1}\frac{u(f^j(w))-u(f^j(w'))}{d_t^j}.$$
Posons $\delta:=|w-w'|$ et
$N:=\log(1/\delta)^{\alpha_0}/\log d_t$ (en principe, nous
devons prendre $N$ la partie enti\`ere de 
$\log(1/\delta)^{\alpha_0}/\log d_t$, cet abus ne change pas
le r\'esultat). On a 
$d_t^N=\delta^{-\alpha_0}$. Nous distinguons trois cas.
\par
Dans le premier cas, on suppose $n\geq N$. 
En prenant $s=N$, on a pour des constantes $c>0$
et $c'>0$
\begin{eqnarray*}
|G(w)-G(w')| & \leq & c\left( \frac{1}{d_t^N}+\sum_{j=0}^{N-1}
  \frac{\delta M_1^j}{d_t^j} \right) 
 =  c\left(\delta^{\alpha_0}
  +\delta\sum_{j=0}^{N-1}\left[\delta^{\frac{\alpha_0-1}{N}}\right]^j
  \right)\\
& = &  c\left(\delta^{\alpha_0}
  +\delta\frac{\delta^{\alpha_0 -1}-1}{\delta^{\frac{\alpha_0-1}{N}}
-1}\right)
 =   c\left(\delta^{\alpha_0}
  +\frac{\delta^{\alpha_0}-\delta}{d_t^{(1-\alpha_0)/\alpha_0}-1}
 \right)\\
& \leq &  c'\delta^{\alpha_0}
\end{eqnarray*}
\par
Pour les autres cas, soit $n$ tel que $w\in U_{-n}\setminus
U_{-n-1}$. On peut encore supposer $m\geq n$. Notons $\rho(w)$ la
distance de $w$ \`a $\K$. On peut suposer $\rho(w)<1$. 
\par
Dans le deuxi\`eme cas, on suppose que $n<N$ et
que de plus $|w-w'|\geq
\rho(w)/2$. Observons que pour $w_0\in \K$ on a
$$c_0\leq |f^n(w)-f^n(w_0)|\leq |w-w_0|M_1^n$$ 
pour un $c_0>0$. Donc
$$\frac{\rho(w)}{c_0}\geq \frac{1}{M_1^n}\ \ \mbox{ et }\ \ 
\frac{1}{d_t^n}=\frac{1}{M_1^{n\alpha_0}}\leq
\left(\frac{\rho(w)}{c_0}\right)^{\alpha_0}.$$
En prenant $s=n$, comme dans
le cas pr\'ec\'edent, on obtient
pour des constantes $c>0$ et $c'>0$
$$|G(w)-G(w')|\leq c\left(\frac{1}{d_t^n}+\frac{\delta
    M_1^n}{d_t^n}\right) \leq
    c\big(\rho(w)^{\alpha_0}+\delta^{\alpha_0}\big)\leq c'
|w-w'|^{\alpha_0}.$$
\par
Avant de traiter le dernier cas, observons que pour tout $w_1\in
    U_{-n_1}\setminus U_{-n_1-1}$ et $w_2\in U_{-n_2}\setminus
    U_{-n_2-1}$  tels que
    $|w_1-w_2|=\rho(w_1)$ on a $|w_1-w_2|\geq \rho(w_2)/2$. 
On utilise les estimations
    pr\'ec\'edentes pour $w:=w_1$ et $w':=w_2$ si
    $n_1\geq n_2$;  pour $w:=w_2$ et $w':=w_1$ sinon. On obtient
$|G(w_1)-G(w_2)|\leq c'|w_1-w_2|^{\alpha_0}$. 
\par
Supposons maintenant que  $n<N$ et $|w-w'|<
\rho(w)/2$. La fonction $G(z)-G(w)$ est
    harmonique dans le disque de centre $w$ et de rayon
    $\rho(w)$. D'apr\`es le principe du maximum,  $|G(z)-G(w)|$ est
    major\'e par $c'\rho(w)^{\alpha_0}$ sur ce disque car c'est le cas
    sur le bord du disque.
Gr\^ace \`a la formule de Poisson, on majore la d\'eriv\'ee de $G$ sur
le disque de rayon $\rho(w)/2$ centr\'e en $w$ par
$c\rho(w)^{\alpha_0-1}$ avec un $c>0$.
On a 
$$|G(w)-G(w')|\leq c|w-w'|\rho(w)^{\alpha_0-1}\leq
c|w-w'|^{\alpha_0}.$$
\par 
2. De fa\c con classique, on consid\`ere une fonction ${\cal
  C}^\infty$, positive 
$\chi$ \'egale \`a 1
sur $B(x,r)$ et \`a support dans $B(x,2r)$ dont le Laplacien est
major\'e par $c'/r^2$ sur $r<|z-x|<2r$ avec $c'>0$. On a
pour une constante $c>0$
\begin{eqnarray*}
\mu(B(x,r))\leq \int\chi\d\mu & = & \int \chi(z) \ddc G(z) = \int
\chi(z) \ddc (G(z)-G(x)) \\
& = & \int\Delta \chi
(G(z)-G(x))\leq cr^{\alpha}.
\end{eqnarray*}
En particulier, en tout point $x$, on a $\limsup_{r\rightarrow 0}
\log\mu(B(x,r))/\log r \geq \alpha$. 
\par
3. C'est une cons\'equence du th\'eor\`eme 3.3.2 car toute surface de
Riemann ouverte est une vari\'et\'e de Stein.
\par
4. Rappelons que {\it la dimension de Hausdorff}
$\HD(\mu)$
d'une mesure de
probabilit\'e
$\mu$ est par d\'efinition la borne inf\'erieure des dimensions de
Hausdorff des bor\'eliens $X$ tels que $\mu(X)=1$. On v\'erifie
facilement que lorsque $\lim_{r\rightarrow 0} \log
  \mu(B(x,r))/\log r$ existe et est constante $\mu$-presque partout
alors elle est \'egale \`a $\HD(\mu)$. Ma\~{n}e \cite{Mane} a montr\'e
que dans le cas des applications polynomiales de $\C$, on a 
$$\lim_{r\rightarrow 0}\frac{\log\mu(B(x,r))}{\log
r}=\frac{\h_\mu(f)}{\int\log |f'|\d\mu}.$$
Sa d\'emonstration n'utilise que le lemme de distorsion de Koebe, elle
est valide dans notre cadre.
\end{preuve}
\subsection{Familles d'applications holomorphes}
Dans ce paragraphe, nous donnons quelques propri\'et\'es des
endomorphismes qui commutent et des endomorphismes d\'ependant
d'un param\`etre.
\par 
Nous avons la proposition suivante qui a \'et\'e d\'emontr\'ee dans
\cite{DinhSibony1} pour les endomorphismes holomorphes de $\P^k$.
\begin{proposition} Soient $f_i:U_i\longrightarrow
V_i$ deux
applications d'allure polynomiale de degr\'e topologique
$d_i\geq 2$ et d'ensemble de Julia rempli $\K_i$ pour $i=1$ ou $2$.
Supposons que $U_1$ (resp. $U_2$) contienne $\K_2\cup
f_2(\K_1)$ (resp.
$\K_1\cup f_1(\K_2)$) et que
$f_1\circ f_2=f_2\circ f_1$ au voisinage de $\K_1\cup \K_2$.
Alors l'ensemble de Julia rempli $\K_2$ (resp. la mesure
d'\'equilibre $\mu_2$ et l'ensemble de Julia $J_k^2$)
de $f_2$ est
\'egal \`a celui de $f_1$.
\end{proposition}
\begin{preuve} On a $f_1\circ f_2(\K_1)=f_2\circ f_1(\K_1)
=f_2(\K_1)$. Donc
$f_2(\K_1)\subset \K_1$ et par suite $\K_1\subset \K_2$ car
$\K_2$ est le plus
grand compact invariant par $f_2$.
De m\^eme, $\K_2\subset \K_1$. Donc $\K_1=\K_2$
\par
Soit $\Omega$ une forme de volume de masse 1
\`a support dans un petit voisinage
de $\K:=\K_1=\K_2$ . D'apr\`es le th\'eor\`eme 3.2.1, 
la mesure $d_1^{-n}(f_1^n)^*\Omega$ tend
vers $\mu_1$ quand $n\rightarrow \infty$. On en d\'eduit que 
$d_2^{-1}d_1^{-n}f_2^*(f_1^n)^*\Omega$
tend vers $d_2^{-1}f_2^*\mu_1$. Comme $f_1$ et $f_2$ commutent, on
a 
$$\frac{f_2^*(f_1^n)^*\Omega}{d_2d_1^n}=
\frac{(f_1^n)^*f_2^*\Omega}{d_2d_1^n}.$$
La derni\`ere mesure tend vers $\mu_1$ car $d_2^{-1}f_2^*\Omega$ est
\'egalement une forme de volume de masse 1. Par cons\'equent,
$d_2^{-1}f_2^*\mu_1=\mu_1$.
D'apr\`es la proposition 3.2.5, pour toute
fonction $\varphi$ p.s.h. au voisinage de $\K$ on a $\int \varphi
\d\mu_1\leq \int \varphi \d\mu_2$. De m\^eme, on a $\int \varphi
\d\mu_2\leq \int \varphi \d\mu_1$. On en d\'eduit que $\mu_1=\mu_2$ et donc
$J_k^1=J_k^2$.
\end{preuve}
\begin{theoreme} Soit $V$ une vari\'et\'e S-convexe. Soit $\Gamma$
  un espace m\'etrique.
Soit $(f_s)_{s\in\Gamma}$ une famille continue
d'applications holomorphes propres
de degr\'e topologique $d_t\geq 2$, $f_s:U_s\longrightarrow V_s$ avec
$V_s\subset V$. On
suppose que pour tout compact $\Gamma_0\subset\Gamma$, $\overline
U_{\Gamma_0} \subset V_{\Gamma_0}$ o\`u 
$$U_{\Gamma_0}:=\big\{(s,z)\in \Gamma_0\times V,\ z\in U_s\big\}$$
et
$$V_{\Gamma_0}:=\big\{(s,z)\in \Gamma_0\times V,\ z\in V_s\big\}.$$
Si l'ensemble
exceptionnel $\E_{s_0}$ de $f_{s_0}$ est contenu dans son ensemble
postcritique d'ordre infini 
alors l'application qui associe \`a $s$ la mesure
d'\'equilibre $\mu_s$ de $f_s$ est continue en $s_0$.
\end{theoreme}
\begin{preuve} Fixons une forme k\"ahl\'erienne
  $\omega=\ddc \Phi$ dans $V$. Notons $J_s$ le
  jacobien r\'eel de $f_s$ pour la m\'etrique consid\'er\'ee.
Notons \'egalement
$M_s$ la famille des mesures de probabilit\'e $\nu_s$ de $U_s$ 
qui v\'erifient $f_s^*\nu_s=d_t\nu_s$ et 
$\int\log J_s \d\nu_s\geq \log d_t$.
\par
V\'erifions que $M_{\Gamma_0}:=\bigcup_{s\in\Gamma_0}M_s$ est
ferm\'e pour la topologie vague.
Soit $s_n\rightarrow s_0$ et soit $\nu$ un point d'adh\'erence
de la suite $(\nu_{s_n})\subset M_{\Gamma_0}$.
On a $f_{s_0}^*\nu=d_t\nu$.
Il suffit de v\'erifier que $\int \log
J_{s_0} \d\nu\geq \log d_t$. Posons pour tout $m\in \R^+$,
$h_{s,m}(z):=\max\big(\log J_s(z),-m\big)$. On a $\int h_{s_n,m}
\d\nu_{s_n}\geq \log
d_t$. Montrons que $\int h_{s_0,m}\d\nu\geq \log d_t$. C'est le cas,
car les applications $(f_s)$ \'etant
holomorphes,
la famille de fonctions $(h_{s,m})_{s\in\Gamma_0}$
est uniform\'ement
continue sur les compacts. 
\par
Soit maintenant $\nu=\lim \mu_{s_n}$ avec $s_n\rightarrow s_0$. 
On sait que $\nu\in M_{s_0}$. Si
$\E_{s_0}$ est contenu dans l'ensemble postcritique de $f_{s_0}$,
alors $\nu(\E_{s_0})=0$ puisque $\log J_{s_0}$ est $\nu$-int\'egrable.
Il en r\'esulte que $d_t^{-n}(f_{s_0}^n)^*\nu\rightharpoonup
\mu_{s_0}$ et donc $\nu=\mu_{s_0}$ car $f_{s_0}^*\nu=d_t\nu$. 
\end{preuve}
\begin{proposition} Soit $\Delta$ une vari\'et\'e complexe connexe. 
Soit $(f_s)_{s\in\Delta}$ une famille
holomorphe d'applications holomorphes de degr\'e topologique
$d_t\geq 2$, satisfaisant aux hypoth\`eses du th\'eor\`eme 3.8.2.
Notons $\Cr_s$ l'ensemble critique de $f_s$. Soit
$(s,z)\mapsto \varphi(s,z)$ une fonction p.s.h. continue \`a l'image
dans $[-\infty,+\infty[$ d\'efinie dans 
$$V_\Delta:=\big\{(s,z)\in \Delta\times V, \  z\in V_s\big\}.$$
Si la fonction
$$\tilde \varphi(s):=\int \varphi(s,z) \d\mu_s(z)$$
n'est pas identiquement \'egale \`a $-\infty$, elle est p.s.h. 
En particulier, lorsque $V$ est un ouvert de $\C^k$, la fonction 
$$h(s):=\int\log J_s \d\mu_s(z)=2\sum_{i=1}^k\lambda_i(s)$$
est p.s.h. 
o\`u les $\lambda_i(s)$ sont les exposants de Lyapounov de $f_s$. Si,
de plus,
$\Cr_s\cap \K_s=\emptyset$ pour tout $s\in \Delta$
alors $h(s)$ est pluriharmonique dans $\Delta$.  
\end{proposition}
\begin{preuve} On peut supposer que $\varphi$ est born\'ee. Soient
  $s_0\in\Delta$ et $s_n\rightarrow s_0$.
Si $\mu_{s_n}\rightharpoonup  \nu$, 
$f_{s_0}^*\nu=d_t\nu$. D'apr\`es le lemme de Hartogs pour $\epsilon>0$, on a
$$\int\varphi(s,z)\d\mu_s \leq \int \big(\varphi(s_0,z)+\epsilon \big)
\d\mu_s.$$ 
Donc
\begin{eqnarray*}
\limsup_{s\rightarrow s_0} \int\varphi(s,z) \d\mu_s & \leq & 
\limsup_{s\rightarrow s_0} \int \big( \varphi(s_0,z) +\epsilon \big)
\d\mu_s \\
& = & \int \big( \varphi(s_0,z) +\epsilon \big)
\d\nu.
\end{eqnarray*}
En utiliant
la proposition 3.2.5, on a
\begin{eqnarray*}
\limsup_{s\rightarrow s_0}\tilde \varphi(s)
& = &
\limsup_{s\rightarrow s_0} \int\varphi(s,z) \d\mu_s
\leq 
\int\varphi(s_0,z) \d\nu \\
& \leq &\int\varphi(s_0,z) \d\mu_{s_0}
= \tilde \varphi(s_0).
\end{eqnarray*}
Donc $\tilde \varphi$ est semi-continue sup\'erieurement.
\par
Notons
$$U_\Delta:=\big\{(s,z)\in \Delta\times V, \  z\in U_s\big\}.$$
Soit $F:U_\Delta\longrightarrow \Delta\times V$ avec
$F(s,z):=(s,f_s(z))$. Posons 
$$\psi(s,z):=\left(\limsup_{n\rightarrow \infty}
\frac{(F^n)_*\varphi}{d_t^n}\right)^*.$$
C'est une fonction p.s.h.
D'apr\`es le lemme 3.2.2, elle
est ind\'ependante de $z$, $\psi(s,z)\geq\tilde\varphi(s)$ et
$\psi(s,z)=\tilde \varphi(s)$ quasi-presque partout. Comme
$\tilde\varphi$ est semi-continue sup\'erieurement, on a
$\tilde\varphi(s)=\psi(s,z)$ partout, elle est donc p.s.h.
\par
Lorsque $V$ est un ouvert de $\C^k$, la fonction $\log J_s$ est
p.s.h. et $\mu_s$-int\'egrable. On peut donc appliquer la propri\'et\'e
prouv\'ee ci-dessus \`a cette fonction.
Si, de plus, $\Cr_s\cap \K_s=\emptyset$ pour tout $s\in\Delta$,
$h(s)$ est pluriharmonique car $\log J_s$ l'est au voisinage de
$\K_s$. 
\end{preuve}
\begin{corollaire} Soit $(f_s)_{s\in\Delta}$ une famille holomorphe
  d'applications \`a allure polynomiale de degr\'e $d_t\geq 2$,
  $f_s:U_s\longrightarrow V_s$, $U_s\subset \subset V_s$ 
et $V_s\subset\subset
  \C$. Alors la fonction $s\mapsto 1/\HD(\mu_s)$ est sous-harmonique.
\end{corollaire}
\begin{preuve} On a vu dans le th\'eor\`eme 3.7.1 que 
$$\frac{1}{\HD(\mu_s)}=\frac{\int\log |f_s'|\d \mu_s}{\log d_t}.$$
La proposition 3.8.3 donne le r\'esultat.
\end{preuve}
%
%\end{document}
\subsection{Mesures PLB}
Nous introduisons une classe de mesures dites mesures PLB. En
dimension 1, ce sont les mesures dont 
le Potentiel est Localement Born\'e. Nous \'etudions aussi les
applications d'allure polynomiale dont la mesure d'\'equilibre est
PLB. On a vu au th\'eor\`eme 3.7.1 qu'en dimension 1 ceci est toujours le cas.
\begin{definition} \rm Soit $U$ une vari\'et\'e complexe.  
Une mesure positive $\nu$ \`a support compact dans $U$ 
est appel\'ee {\it PLB} 
si les
fonctions p.s.h sont $\nu$-int\'egrables, \cad que 
pour toute fonction $\varphi$ p.s.h dans $U$ on a
$\int \varphi \d \nu >-\infty$.
\end{definition} 
\par 
Soient $U$, $V$ des vari\'et\'es complexes et
$f:U\longrightarrow V$ un rev\^etement ramifi\'e fini. 
Alors si $\sigma$ est
une mesure PLB de $V$, $f^*\sigma$ est une mesure PLB de $U$; 
si $\nu$ est
une mesure PLB de $U$, $f_*\nu$ est une mesure PLB de $V$.  
\par
Rappelons qu'un ensemble $E$ d'une vari\'et\'e $V$ est {\it
  pluripolaire} 
\index{ensemble!pluripolaire}
si pour tout point $p$ de $V$, il existe un voisinage
  $W$ de $p$ et une fonction p.s.h $u$ dans $W$ tels que 
$$E\cap W\subset \big\{z\in W,\ u(z)=-\infty \big\}.$$
Il est clair que les mesures PLB ne chargent pas les ensembles
pluripolaires.     
\begin{proposition} Soient $V$ une vari\'et\'e complexe et
  $U\subset\subset V$ un ouvert.
Soit $\nu$ une mesure PLB de $U$ 
\`a support dans un compact $K\subset U$.
Alors pour tout $p\geq 1$
\begin{enumerate}
\item Il existe une constante $c>0$ telle que pour toute fonction
  p.s.h $\psi$ dans $U$ on ait $\|\psi\|_{\Lone(\nu)}\leq
  c\|\psi\|_{\Lp(U)}$. 
\item Il existe une constante $0<c<1$ telle que pour toute fonction
  $\psi$ p.s.h dans $V$ et satisfaisant $\int \psi\d \nu=0$, 
on ait $\sup_K  \psi \leq c\sup_V \psi$.
%\item (????) Le c\^one $\{\psi \mbox{ p.s.h., }\int\psi\d\nu=0,
%  \|\psi\|_{\Lp(U)}\leq 1\}$ est ferm\'e dans $\Lp(U)$.  
\end{enumerate}
\end{proposition}
\begin{preuve} 
1. Si la propri\'et\'e 1 n'est pas v\'erifi\'ee, il
  existe des fonctions $\psi_j$, p.s.h sur $U$,
$\int |\psi_j| \d\nu=1$ et
  $\|\psi_j\|_{\Lp(U)}\leq j^{-2}$. 
L'in\'egalit\'e de sous-moyenne implique que pour tout compact $W$ de
  $U$ on a  $\psi_j\leq c_Wj^{-2}$ sur
$W$ o\`u $c_W>0$ est une constante. 
Posons
  $\psi:=\sum \psi_j$. La fonction $\psi$ est bien d\'efinie,
  p.s.h dans $U$ et $\int\psi\d\nu=-\infty$. 
C'est la contradiction recherch\'ee.
\par
2. La famille $\F$ des fonctions
$\psi$ p.s.h dans $V$
v\'erifiant $\int\psi \d\nu=0$ et $\psi\leq 1$, est
relativement compacte dans $\Lploc(V)$. D'apr\`es la propri\'et\'e
1, la fonction identiquement nulle n'est pas dans
l'adh\'erence de la famille $\{\psi-1,\ \psi\in\F\}$. Par suite,
la fonction identiquement \'egale \`a 1 n'est pas dans
l'adh\'erence de $\F$. Ceci implique la propri\'et\'e 2.
\end{preuve}
\par
Le r\'esultat suivant qui justifie la terminologie
choisie:
\begin{corollaire} Soit $\nu$ une mesure \`a support compact et PLB dans
un ouvert born\'e $W$ de $\C^k$. 
Soit $\psi$ une fonction p.s.h dans $\C^k$. Alors la
  fonction $G_\psi(z):= \int \psi(z-\zeta)\d\nu(\zeta)$ est p.s.h et
  localement born\'ee dans $\C^k$. En particulier, une mesure $\nu$
  d'une surface de Riemann ouverte est 
  PLB si et seulement si elle est \`a potentiel
  localement born\'e. 
\end{corollaire}
\begin{preuve} Il est clair que $G_\psi$ est p.s.h. 
On peut supposer que $W$ est la boule $B(0,R)$ de centre $0$ et de
rayon $R>0$.  
Soient $r>0$ et 
  $W':=B(0,R+r)$.
Posons $W_s:=\{s-z,\ \mbox{ avec } z\in W\}$.
D'apr\`es la Proposition 3.9.2, il existe une constante $c>0$ telle
  que si $|s| <r$ on ait
\begin{eqnarray*}
|G_\psi(s)| & \leq &  \|\psi(s-\zeta)\|_{\Lone(\nu)}
\leq  c\|\psi(s-\zeta)\|_{\Lone(W)}\\
& = & c\|\psi\|_{\Lone(W_s)} \leq c\|\psi\|_{\Lone(W')}.
\end{eqnarray*}
Donc $G_\psi$ est born\'ee dans la boule de centre $0$ et
de rayon $r$.
\par
Consid\'erons maintenant
une mesure $\nu$ d'une surface de Riemann ouverte. 
On peut supposer que $\nu$ est une 
mesure \`a support compact dans
$\C$. En prenant $\psi(z):=\log |z|$, on obtient $G_\psi$ un potentiel
de $\nu$. On a montr\'e dans la premi\`ere partie que si $\nu$ est
PLB, son potentiel $G_\psi$ est localement born\'e.
\par
R\'eciproquement, supposons que $\nu=\ddc G$ avec $G$ une fonction 
sousharmonique localement born\'ee. Il suffit d'appliquer
\cite[proposition 2.10]{Demailly} dans un ouvert de $\C$, qui affirme que
les fonctions sous harmoniques sont int\'egrables par rapport \`a
$\nu=\ddc G$ lorsque $G$ est born\'ee.
\end{preuve}
\par
Fixons $p\geq 1$. Pour tout $c>0$, notons $\M^p_c(K,U)$ 
la famille des mesures positives $\nu$
port\'ees par $K$ v\'erifiant   $\|\psi\|_{\Lone(\nu)}\leq
  c\|\psi\|_{\Lp(U)}$ pour toute fonction p.s.h $\psi$. 
Observons que l'in\'egalit\'e
$|\int\psi\d\nu|\leq c'\|\psi\|_{\Lp(U)}$ pour $c'>0$
implique la condition  $\|\psi\|_{\Lone(\nu)}\leq
  c\|\psi\|_{\Lp(U)}$ pour $c>0$ convenable. 
En effet, l'in\'egalit\'e de sous-moyenne
  implique  $\sup_K \psi^+ \leq
  c''\|\psi^+\|_{\Lp(U)}$; d'o\`u 
$$\int|\psi|\d\nu =\int(-\psi +2
    \psi^+)\d\nu \leq \left|\int \psi \d\nu \right|+2\sup_U\psi^+\leq
    c\|\psi\|_{\Lp(U)}.$$ 
Il est clair que $\M^p_c(K,U)$ est un compact
  convexe et que pour toute fonction $0\leq h\leq 1$ on a
  $h\nu\in\M^p_c(K,U)$. 
\par
Soient $u_1,\ldots,u_k$ des fonctions p.s.h born\'ees dans
  $U$. 
D'apr\`es l'in\'egalit\'e de Chern-Levine-Nirenberg \cite{Demailly}, 
pour toute fonction test $\chi\geq 0$
  \`a support compact, la mesure 
  $\chi \ddc u_k\wedge\ldots\wedge \ddc u_1$ 
est PLB. Elle appartient \`a $\M^p_c(K,U)$ pour un
  $c>0$ convenable. 
\par
Bedford-Taylor \cite{BedfordTaylor} 
ont montr\'e que pour tout compact non pluripolaire 
$K\subset \C^k$, il existe des fonctions
p.s.h born\'ees $u_1,\ldots,u_k$
telles que la mesure $\nu:=\ddc u_k\wedge\ldots\wedge \ddc u_1$
soit \`a support compact et 
v\'erifie $\nu(K)>0$. On a donc la proposition suivante:
\begin{proposition} Soit $(v_i)$ une suite de fonctions
  p.s.h convergeant dans $\Lp(U)$ vers une fonction pluriharmonique
  $v$. Alors on peut en extraire une sous-suite convergeant
  vers $v$ hors d'un  ensemble pluripolaire. 
\end{proposition}
\begin{preuve}
Quitte \`a remplacer $v_j$ par $v_j-v$ on peut supposer $v=0$. Quitte
\`a extraire une sous-suite on peut supposer que $\sum
\|v_j\|_{\Lp(U)}<\infty$. Il en r\'esulte que pour toute mesure
$\nu$ PLB dans un ouvert relativement compact de $U$, 
la s\'erie $\sum |v_i|$ converge $\nu$-presque partout. On a
donc $v_j\rightarrow 0$ $\nu$-presque partout. Le r\'esultat de
Bedford-Taylor rappel\'e ci-dessus entra\^{\i}ne qu'on a la
convergence hors d'un ensemble pluripolaire. Lelong d\'ej\`a avait observ\'e
que dans ce cas l'ensemble $\big\{ \limsup v_j<v \big\}$ est
pluripolaire. 
\end{preuve}
\par
Le th\'eor\`eme suivant donne des crit\`eres pour que
la mesure d'\'equilibre $\mu$ 
d'une application d'allure polynomiale $f$ soit
PLB. Observons que $\mu$ est PLB dans $U$ si et seulement si elle l'est
dans $V$. En effet, $\varphi$ est $\mu$-int\'egrable si et seulement si
$\Lambda\varphi$ l'est. Dans la suite, on suppose que $V$ est de
Stein. 
\begin{theoreme} Soit $f:U\longrightarrow V$ une application d'allure
  polynomiale. Supposons que $V$ est de Stein. 
Alors les conditions suivantes sont \'equivalentes:
\begin{enumerate}
\item La mesure d'\'equilibre $\mu$ est PLB.
\item Pour toute $\varphi$ p.s.h., $(\Lambda^n\varphi)$ ne tend
  pas uniform\'ement vers $-\infty$. 
\item Il existe $0<c_2<1$ telle que pour toute $\varphi$
  p.s.h. et $n\geq 0$, on ait 
$0\leq\sup_U (\Lambda^{n+1}\varphi-c_\varphi)\leq 
c_2\sup_U(\Lambda^n\varphi-c_\varphi)$ o\`u $c_\varphi$ est une constante.
\item Il existe des constantes $A>0$ et $0<c_3<1$ telles que pour
  toute $\varphi$ p.s.h. on ait $\|\Lambda^n(\ddc\varphi)\|_U
\leq Ac_3^n\|\ddc \varphi\|_U$. 
\item Pour toute $\varphi$ p.s.h., il existe une constante $0<c_4<1$
  telle que $\|\Lambda^n(\ddc\varphi)\|_U=\O(c_4^n)$.
\end{enumerate}
\end{theoreme} 
\par
Notons $H$ le sous-espace des fonction pluriharmoniques dans
$\Ltwo(U)$ et $E$ sont orthogonal. Notons $E^*$ l'ensemble des
fonctions p.s.h. appartenant \`a $E$. Soit $\varphi\in \Ltwo(U)$ 
une fonction p.s.h. dans $U$. On a la d\'ecomposition $\varphi=u+v$
avec $u\in H$ et $v\in E^*$. Il existe des applications lin\'eaires
$\Lambda_1:H\longrightarrow H$, $\Lambda_2:E^*\longrightarrow H$ et 
$\Lambda_3:E^*\longrightarrow E^*$ telles que
$\Lambda\varphi=(\Lambda_1u+\Lambda_2v, \Lambda_3v)$. D'apr\`es la
proposition 3.2.5, on a $\|\Lambda_1^n\|\leq Ac_1^n$ et
$\|\Lambda_2\|\leq A$ o\`u $A>0$ est une constante. On a aussi
$$\Lambda^n\varphi=(\Lambda_1^nu+\Lambda_1^{n-1}\Lambda_2v+
\Lambda_1^{n-2}\Lambda_2\Lambda_3 v+\cdots +
\Lambda_2\Lambda_3^{n-1}v, \Lambda_3^nv).$$
Nous avons besoin du lemme suivant:
\begin{lemme}
Il existe une constante $B>0$ telle que pour toute $v\in
  E^*$ on ait 
  $\|\Lambda_3v\|_{\Ltwo(U)}\leq B\|\ddc v\|_U$.
\end{lemme}
\begin{preuve} Soit $W\subset U$ un ouvert contenant $\overline U_{-2}$.
Comme $V$ est de Stein, $U$ est aussi de Stein. 
Il existe une constante $B'>0$
  ind\'ependante de $v$ et un potentiel $\psi$ de $\ddc v$ dans
  $U$ qui v\'erifie $\|\psi\|_{\Ltwo(W)}\leq B'\|\ddc v\|_U$. 
D'apr\`es la proposition 3.2.5,
$\Lambda_3:\PSH(W)\cap\Ltwo(W)\longrightarrow \PSH(U)\cap \Ltwo(U)$
est born\'e. Donc
$\|\Lambda_3\psi\|_{\Ltwo(U)}\leq 
B\|\ddc v\|_U$ pour $B>0$ convenable.
D'autre part, puisque $\Lambda_3v\in E$, on a 
$\|\Lambda_3 v\|_{\Ltwo(U)}\leq  \|\Lambda_3 
\psi\|_{\Ltwo(U)}$. On obtient finalement 
$\|\Lambda_3 v\|_{\Ltwo(U)}\leq B\|\ddc v\|_U$. 
\end{preuve}
\par\noindent
{\it Preuve du th\'eor\`eme 3.9.5---}
1. $\Longleftrightarrow$ 2., 3. $\Longrightarrow$ 1.
et 4. $\Longrightarrow$ 5. sont claires.
\par
1. $\Longrightarrow$ 3. On prend $c_\varphi:=\int\varphi\d\mu$ et on
peut supposer $c_\varphi=0$. D'apr\`es la proposition 3.9.2
(appliqu\'ee \`a des ouverts convenables), on a 
$\sup_U\Lambda\varphi\leq
  \sup_{U_{-2}}\varphi\leq c_2\sup_U\varphi$ o\`u $0<c_2<1$ est une
  constante. On a aussi
  $\sup_U\Lambda\varphi\geq 0$ car $\int\Lambda\varphi\d\mu=0$.  
\par
1. 2. et 
3. $\Longrightarrow$ 4. 
Soient $\varphi^j$ des fonctions p.s.h. v\'erifiant
$\|\ddc\varphi^n\|_U=1$. 
D'apr\`es le lemme 3.9.6, il existe $\psi^n$ 
telle que $\|\psi^n\|_{\Ltwo(U)}\leq B$ et
$\ddc\psi^n=\ddc\Lambda\varphi^n$. Il existe donc une constante $C>0$
telle que $\sup_U\Lambda\psi^n\leq\sup_{U_{-1}}\psi^n\leq C$. 
D'apr\`es le point 1 de la proposition 3.9.2, il existe $C'>0$
telle que $|c_{\psi^n}|=|c_{\Lambda\psi^n}|\leq C'$. 
Posons $\phi^n:=\Lambda\psi^n-c_{\psi^n}$. On a $\int 
\phi^n \d\mu=0$ et
$\sup\phi^n\leq C+C'$. D'apr\`es 3., la
famille $c_2^{-n}\Lambda^{n-1}\phi^n$ est 
born\'ee sup\'erieurement sur $U$. 
On en d\'eduit que la famille $c_2^{-n}\Lambda^{n}\phi^n$ est born\'ee
sur $V$.
Puisque   $\int \phi^n
\d\mu=0$, aucune sous-suite de $\phi^n$ ne tend 
uniform\'ement vers $-\infty$. Par cons\'equent, les courants 
$c_2^{-n}\Lambda^{n+2}\ddc\varphi^n=c_2^{-n}\Lambda^n\ddc
\phi^n$ sont de masse uniform\'ement born\'ee dans $U$. Ceci implique
la propri\'et\'e 4.
\par
5. $\Longrightarrow$ 2. Dans la suite, les constantes $A_1$, $A_2$,
$A_3$ et $A_4$ sont positives et convenablement choisises.
D'apr\`es le lemme 3.9.6 et la propri\'et\'e 5, on a 
$$\|\Lambda_3^n
v\|_{\Ltwo(U)}\leq A_1\|\ddc\Lambda^{n-1}\varphi\|_U\leq
A_2c_4^n.$$ 
Posons $b:=\int u\d\mu$,
$b_n:=\int\Lambda_2\Lambda_3^n v\d\mu$ et 
$s_n:=b+b_1+\cdots+b_{n-1}$.
Puisque $\Lambda_2$ est born\'e et $\Lambda_2\Lambda_3^n v$ est
pluriharmonique, on a $|b_n|\leq A_3c_4^n$ 
et donc la suite $(s_n)$ converge. 
On a aussi
\begin{eqnarray*}
\|\Lambda^n\varphi-s_n\|_{\Ltwo(U)} & \leq &
\|\Lambda_1^nu-b\|_{\Ltwo(U)} + 
\|\Lambda_1^{n-1}\Lambda_2v-b_1\|_{\Ltwo(U)} + \\
& & +\cdots +
 \|\Lambda_2\Lambda_3^{n-1}v-b_{n-1}\|_{\Ltwo(U)}+\|\Lambda_3^nv\|\\
& \leq &  
A_3 (c_1^n+c_1^{n-1}+\cdots+c_4^{n-1}+c_4^n)\\
& \leq &  A_4c^n
\end{eqnarray*}
o\`u la constante $c_1$ est donn\'ee dans 3.2.5, $\tilde c<c<1$ et
$\tilde c:=\max(c_1,c_4)$. On en d\'eduit que $\varphi$ est
$\mu$-int\'egrable. Observons que la meilleure estimation est 
$A_4\tilde c^n$ si $c_1\not=c_4$ et 
$A_4n\tilde c^n$ sinon.
\par
\hfill $\square$
\begin{corollaire} Soient $U\subset\subset V$ les ouverts d'une
  vari\'et\'e complexe $M$, $V$ de Stein 
et $f:U\longrightarrow V$ une application
  d'allure polynomiale de degr\'e topologique $d_t\geq 2$. Supposons
  que sa mesure d'\'equilibre $\mu$ est PLB. Alors pour toute
  pertubation $f_\epsilon:U_\epsilon\longrightarrow V_\epsilon$ 
suffisamment proche 
de $f$ avec $U_\epsilon\subset\subset V_\epsilon\subset W$ 
la mesure  d'\'equilibre $\mu_\epsilon$ de $f_\epsilon$ est PLB.   
\end{corollaire}
\begin{preuve} Quitte \`a modifier l\'eg\`erement les ouverts $V$ et
  $V_\epsilon$, on peut supposer que $V_\epsilon=V$. 
Soit $M>0$ tel que $f^*\omega^{k-1}\leq
M\omega^{k-1}$ sur $U$.
D'apr\`es le th\'eor\`eme 
3.9.5, il existe $n_0\geq 1$ et $0<c<d_t/M$ tel que 
pour toute 
$\varphi$ p.s.h. on a
$\|\Lambda^{n_0-1}\ddc\varphi\|_U\leq c\|\ddc\varphi\|$.
\par
Fixons un $0<\alpha<1-cMd_t^{-1}$. 
Supposons que $f_\epsilon$ 
soit assez proche de $f$ dans le sens o\`u 
\begin{enumerate}
\item $(f_\epsilon^{n_0})^*\omega^{k-1}- (f^{n_0})^*\omega^{k-1} 
\leq d_t^{n_0}\alpha \omega^{k-1}$
sur $W:=f_\epsilon^{-n_0}(U)$.
\item $f^{n_0-1}(W)\subset U$.
\end{enumerate}
On a
\begin{eqnarray*}
\|\Lambda_\epsilon^{n_0}\ddc \varphi\|_{U} & = 
& d_t^{-n_0}\int_W \ddc\varphi \wedge (f^{n_0})^*
\omega^{k-1}+\\
& & +  d_t^{-n_0}
\int_W \ddc\varphi\wedge \big[(f_\epsilon^{n_0})^*\omega^{k-1}-
(f^{n_0})^*\omega^{k-1}\big]\\
& \leq & d_t^{-1}
\int_{U} \Lambda^{n_0-1}\ddc\varphi\wedge f^*\omega^{k-1}
+\alpha\|\ddc\varphi\|_{U}\\
& \leq & 
Md_t^{-1}\|\Lambda^{n_0-1}\ddc\varphi\|_{U}+\alpha\|\ddc\varphi\|_{U}\\
&\leq & (cMd_t^{-1}+\alpha)\|\ddc\varphi\|_{U}.
\end{eqnarray*}
Par cons\'equent, $\|\Lambda_\epsilon^n\ddc\varphi\|_U=\O(\tilde
c_4^n)$ avec $\tilde c_4:=cMd_t^{-1}+\alpha<1$. 
D'apr\`es le th\'eor\`eme 3.9.5 (appliqu\'e \`a l'application
$f_\epsilon^{n_0}$), la mesure $\mu_\epsilon$ est PLB. 
\end{preuve}
\par
Pour estimer la vitesse de m\'elange, nous avons la proposition
suivante:
\begin{proposition} 
Supposons qu'il existe une constante $0<c_5<1$ telle
  que $\|\Lambda^n\omega\|=\O(c_5^n)$ (c'est le cas si $d_{k-1}<c_5d_t$). 
Alors la vitesse de m\'elange de
  $\mu$ est exponentielle d'ordre $c^n$ o\`u $c=\min(c_1,c_5)$ si
  $c_1\not=c_5$ et $c_1<c<1$ si $c_1=c_5$. Plus
  pr\'ecis\'ement, il existe une constante $A>0$ telle que 
pour toute $\varphi$
de classe ${\cal C}^2$ et toute $\psi$ born\'ee, on ait
$$|I_n|:=\left|\int \psi(f^n) \varphi \d\mu -\left(\int \psi
  \d\mu\right)\left(\int \varphi \d\mu\right)\right|\leq
  A\|\psi\|_\infty\|\varphi\|_{{\cal C}^2} c^n.$$ 
\end{proposition}
\begin{preuve} Observons que la valeur de $|I_n|$ 
ne change pas si l'on remplace
  $\psi$ par $-\psi$ ou par $\psi+M$. On peut donc supposer que
  $\psi\geq 0$.
La fonction $\varphi$ s'\'ecrit comme diff\'erence de
  fonctions ${\cal C}^2$ p.s.h. On peut supposer que $\varphi$
  est p.s.h. 
On a $\ddc\varphi\leq A_5\|\varphi\|_{{\cal C}^2}\omega$.
 On en d\'eduit que $\|\ddc\Lambda^n\varphi\|_U
  \leq A_6\|\varphi\|_{{\cal C}^2}c_5^n$.
On utilise les calculs d\'ej\`a faits en 3.9.5 en rempla\c cant
  $c_4$ par $c_5$. On a $c_\varphi=\lim s_n$. Il est facile de voir
  que $|c_\varphi-s_n|\leq A_7\|\varphi\|_{{\cal C}^2}c_5^n$.
Par cons\'equent, 
$\|\Lambda^n\varphi-c_\varphi\|_{\Ltwo(U)}\leq A_8\|\varphi\|_{{\cal
  C}^2}c^n$. 
L'in\'egalit\'e de sous-moyenne implique que  
$\Lambda^n\varphi-c_\varphi\leq A_9\|\varphi\|_{{\cal
  C}^2}c^n$ sur $\K$.
Utilisant
  l'invariance de $\mu$ on obtient
$$I_n=\int(\Lambda^n\varphi-c_\varphi)\psi\d\mu\leq
A_9c^n\|\varphi\|_{{\cal C}^2}\|\psi\|_\infty.$$
Rempla\c cant $\psi$ par $\|\psi\|_\infty-\psi$ on obtient une
in\'egalit\'e analogue pour $-I_n$. 
\end{preuve}
\par
Observons que si $V$ n'est pas de Stein, dans le th\'eor\`eme 3.9.5, 
on a encore 1. $\Longrightarrow$ 2. et 1. $\Longrightarrow$ 5;
la derni\`ere proposition reste aussi valable 
et si la condition $\|\Lambda^n\omega\|=\O(c_5^n)$ est remplac\'ee par 
la condition 2. du th\'eor\`eme 3.9.5.
Les th\'eor\`emes 3.9.5, 3.4.4,
3.5.1, 3.6.1 impliquent le corollaire suivant o\`u on ne suppose pas que
$V$ est de Stein.
\begin{corollaire} Soit $f:U\longrightarrow V$ une 
application d'allure
  polynomiale. Supposons que sa mesure d'\'equilibre $\mu$ est
  PLB. Alors $d_{k-1}<d_t$, $\delta<1$. De plus 
\begin{enumerate}
\item L'ensemble exceptionnel $\E$ est un sous-ensemble 
analytique de $V$ totalement invariant par $f$. 
\item Les points p\'eriodiques r\'epulsifs sont denses dans
  $\supp(\mu)$. 
\item Les exposants de Lyapounov de $\mu$ sont strictement positifs.
\item La mesure $\mu$ est 
  m\'elangeante \`a vitesse exponentielle.
\item Si $V$ est contenue dans une vari\'et\'e de Stein, $\mu$ est
  d'entropie maximale $\log d_t$.
\end{enumerate}
\end{corollaire}
\subsection{Exemples et remarques}
Soit $f$ un endomorphisme polynomial de $\C^k$. Il existe 
$\lambda>0$ et $l>0$ tels que $|f(z)|\geq \lambda |z|^l$ pour $z$
suffisamment grand. La meilleur constante $l$ existe et appel\'ee 
{\it l'exposant de Lojasiewicz} de $f$ \cite{Ploski}. On suppose que $l>1$ ou 
$l=1$ et $\lambda>1$. Si $V$ est une boule assez grande, on a
$U:=f^{-1}(V)\subset\subset V$ et la restriction de $f$ sur $U$ est
une application d'allure polynomiale. 
\begin{proposition} Pour toute fonction $\varphi$, 
$\mu$-int\'egrable, on a $\|\ddc\Lambda^n\varphi\|_U=\o(1/n)$
si $\lambda>1$ et $l=1$. Si $l>1$, on a  
$\|\ddc\Lambda^n\varphi\|_U=\O(l^{-n})$.
En particulier, dans le deuxi\`eme cas, la vitesse de
m\'elange est d'ordre $l^{-n}$.
\end{proposition}  
\begin{preuve} On consid\`ere le cas o\`u $\lambda>1$ et $l=1$. Le
  deuxi\`eme cas se traite de la m\^eme mani\`ere.
Pour simplifier les notations, on peut supposer que $V$
  est la boule unit\'e et que $\varphi$ est une 
fonction p.s.h. sur $V$
  v\'erifiant $\sup_V \varphi=1$, $\int\varphi\d\mu=0$. On peut
  supposer \'egalement que $|f(z)|\geq\lambda|z|$ pour $|z|\geq 1$.
D'apr\`es la proposition 3.2.5, la suite de
  fonctions $\Lambda^n\varphi$ 
tend vers $0$ dans
  $\Ltwoloc(\C^k)$. Puisque $|f(z)|\geq \lambda|z|$ 
pour $|z|\geq 1$, on a 
$f^{-n}(B_{\lambda^n})\subset V$ o\`u $B_{\lambda^n}$ 
d\'esigne la boule de rayon
  $\lambda^n$ centr\'ee en $0$. Par suite, 
$\Lambda^n\varphi \leq 1$ sur $B_{\lambda^n}$. 
\par
Montrons d'abord que la suite
$$s_n:=n\int_V \ddc\Lambda^n\varphi  \wedge \omega^{k-1}$$
est born\'ee. 
Comme $\Lambda^n\varphi $ tend vers $0$ dans $\Ltwo_\loc(\C^k)$, 
on peut supposer
qu'il existe un $r_0\geq 1$ tel que 
$$\lim\int_{|z|=r_0} \Lambda^n\varphi  \d\sigma_{r_0} =0$$
o\`u $\sigma_{r_0}$ est la mesure de Lebesgue normalis\'ee de 
masse $1$ sur la sph\`ere
$\{|z|=r_0\}$.  
D'apr\`es la formule de Posson-Jensen, on a pour tout
$r_0<r\leq \lambda^n$
\begin{eqnarray*}
(\log r -\log r_0) \int_{B_{r_0}}
\ddc\Lambda^n\varphi  \wedge \omega^{k-1} & \leq & 
\int_{|z|=r}\Lambda^n\varphi  \d\sigma_r -  \int_{|z|=r_0}
\Lambda^n\varphi  \d\sigma_{r_0} \\
& \leq & 1-\int_{|z|=r_0}
\Lambda^n\varphi  \d\sigma_{r_0}.
\end{eqnarray*}
Appliquons cette in\'egalit\'e \`a $r=\lambda^n$, on obtient
$$\limsup n\int_{B_{r_0}} \ddc\Lambda^n\varphi  
\wedge \omega^{k-1}\leq 1.$$
Par cons\'equent, il existe $A>0$ telle que 
$\|\ddc\Lambda^n\varphi\|_U\leq An^{-1}$. Par homoth\'etie, pour $A>0$
convenable,  
 $\|\ddc\Lambda^n\varphi\|_U\leq 
A'\|\varphi\|_{\Ltwo(B_2)}n^{-1}$ pour toute
 $\varphi$ v\'erifiant $\int\varphi\d\mu=0$. On obtient en
 particulier, 
$\|\ddc\Lambda^n\varphi\|_U\leq 2A'\|\Lambda^{[n/2]}\varphi
\|_{\Ltwo(B_2)}n^{-1}$. 
Ceci implique
que $\|\ddc\Lambda^n\varphi\|_U=\o(1/n)$ car 
$\Lambda^{[n/2]}\varphi$ tend
vers $0$.
\par
Observons que si $l>1$, on peut supposer que $U$ est suffisamment
petit par rapport \`a $V$ de sorte que la constante $c_1$ soit
strictement inf\'erieure \`a $l$. On peut appliquer la proposition
3.9.8 pour obtenir la vitesse de m\'elange.
\end{preuve}
\begin{remarque} \rm Pour tout
  endomorphisme holomorphe $f$ de degr\'e alg\'ebrique 
$d>1$ de $\P^k$, la
  vitesse de m\'elange est d'ordre $d^{-n}$. En effet, il suffit
  d'appliquer la proposition 3.10.1 \`a un relev\'e polynomial de
$f$ dans $\C^{k+1}$. Ant\'erieurement, Forn\ae ss-Sibony ont montr\'e
que la vitesse est d'ordre $(d-\epsilon)^{-n}$ \cite{FornaessSibony2}.  
\end{remarque}
\begin{exemple} \rm
Consid\'erons l'endomorphisme $f:\C^2\longrightarrow \C^2$
d\'efini par $f(z_1,z_2):=(\lambda
z_1 +P(z_2), Q(z_2))$ o\`u $P$, $Q$ 
sont des polyn\^omes avec $d_t:=\deg Q\geq 2$ et $\lambda\in \C$,
$|\lambda|>1$. 
Cette application est de degr\'e topologique $d_t$. On montre
facilement que pour tout
$1<\lambda'<|\lambda|$ on a $|f(z)|\geq \lambda'|z|$ lorsque $z$ est
suffisament grand.
\par
La dynamique de l'application $f$ est facile \`a \'etudier. 
Notons $\K_Q$, $J_Q$ et $\mu_Q$ l'ensemble de Julia rempli,
l'ensemble de Julia et la mesure d'\'equilibre du polyn\^ome 
$Q$ ($J_Q$ est donc le bord de $\K_Q$).
On a
$$f^n(z)=\left(\lambda^n(z_1+\lambda^{-1}P(z_2)+\cdots+\lambda^{-n}
 P\circ Q^{n-1}(z_2)),Q^n(z_2)\right).$$
Il est clair que la suite $f^n(z)$ est born\'ee si et seulement
si $z_2\in \K_Q$
et $z_1=h(z_2):=
-\sum_{j=1}^\infty \lambda^{-j}P\circ Q^{j-1}(z_2)$. Par
cons\'equent, l'ensemble de Julia rempli
$\K$ de $f$ est le graphe de la
fonction continue $h$ au dessus de $\K_Q$. Observons qu'elle satisfait
 l'\'equation fonctionnelle $h \circ Q=\lambda h+P$. Pour tout $s>0$
 fix\'e, lorsque $\lambda$ est
 suffisamment grand, la fonction $h$ est de 
classe ${\cal C}^s$ sur $\K_Q$.
On v\'erifie sans peine que $\supp(\mu)$
est le graphe de $h$ au dessus de $\partial \K_Q$ et
$\mu=\pi^*\mu_Q$ o\`u $\pi$ est la
projection de $\supp(\mu)$ dans $J_Q$. Il est clair aussi que les
deux exposants de Lyapounov sont $\log|\lambda|$ et celui de $Q$.
\par
Il y a $d_t^n$ points p\'eriodiques de p\'eriode $n$, leurs
multiplicateurs sont $\lambda^n$ et $(Q^n)'(z_2)$. Il est facile de
montrer que les mesures $\nu_n$ d\'efinies par des masses de Dirac 
\'equidistribu\'ees aux points
p\'eriodiques r\'epulsifs d'ordre $n$ convergent vers $\mu$. Si
$Q(z_2)=z_2^{d_t}$ et $P(0)=0$, le point selle $(0,0)$ est un point
p\'eriodique isol\'e dans l'ensemble $\K$ des points d'orbite 
born\'e. Le graphe de $h$
apparait comme sa vari\'et\'e stable.
\par
Si $Q(z_2)=z_2^{d_t}$ et si $P$ est non constant, 
la fonction $h$ admet le
disque unit\'e ferm\'e pour domaine d'existence.
Dans ce cas, l'ensemble exceptionnel $\E$ est
\'egal \`a $\{z_2=0\}$.
Dans le cas o\`u $Q$ n'est pas conjugu\'e \`a $z_2^{d_t}$,
l'ensemble exceptionnel de $Q$ est vide. Celui de $f$ est donc aussi vide.
\par Si $Q$ est un polyn\^ome de Tchebychev, $\supp(\mu)$ et $J_Q$ sont
des courbes r\'eelles. Si $J_Q$ est un Cantor, $\supp(\mu)$ l'est
aussi.
\par
Si $P=0$, la mesure $\mu$ est port\'ee par l'ensemble analytique
$\{z_1=0\}$. 
On peut v\'erifier que pour tout courant positif
ferm\'e $\ddc\varphi$ de bidegr\'e $(1,1)$ de $\C^2$, la s\'erie
$\sum \Lambda^n \ddc\varphi$ converge. En consid\'erant
$\varphi=\varphi(|z|)$ avec $\varphi(0)=0$, on v\'erifie qu'il
n'existe pas d'estimation sur $\|\Lambda^n\ddc\varphi\|$ qui est
uniforme en $\varphi$ et qui implique la convergence de
$\sum\Lambda^n\ddc\varphi$. 
\end{exemple}
\begin{theoreme}
Soit $f:\C^k\longrightarrow \C^k$ 
une application polynomiale propre de degr\'e alg\'ebrique
$d>1$ et de degr\'e topologique $d_t>d^{k-1}$. 
Supposons qu'il existe un ouvert de Stein $V\subset \C^k$ tel que 
$U:=f^{-1}(V)\subset\subset V$. 
Alors il existe une constante
  $A>0$ telle que $\|\Lambda^nT\|_{V}\leq
  A\alpha^{n}\|T\|_{V}$ pour tout
  courant positif, ferm\'e $T$ de bidegr\'e $(1,1)$ dans $V$ o\`u
  $\alpha :=d_t^{-1}d^{k-1}$. En particulier, la mesure $\mu$ est
  PLB. 
\end{theoreme}
\begin{preuve}  
Par d\'efinition, il existe une constante $A'>0$ telle
  que  $\log(1+|f^n(z)|)\leq A'd^n$ sur $V$.
D'apr\`es l'in\'egalit\'e de Chern-Levine-Nirenberg \cite{Demailly}, 
on a pour une constante $A>0$ 
\begin{eqnarray*}
\|\Lambda^nT\|_{V} & \leq &
d_t^{-n}\int_{U_{-n}}T\wedge (f^n)^*\omega^{k-1} \leq 
d_t^{-n}\int_{U} T\wedge (f^n)^*\omega^{k-1}\\
& = & d_t^{-n}\int_{U} T\wedge (\ddc \log(1+|f^n(z)|))^{k-1}
\leq A\alpha^n \|T\|_{V}.
\end{eqnarray*}
\end{preuve}
\begin{remarque} \rm
D'apr\`es le corollaire 3.9.7, en pertubant l'application $f$
ci-dessus, on peut construire des familles
d'applications d'allure polynomiale dont la mesure d'\'equilibre est
PLB. 
\end{remarque}
\begin{exemple} \rm Consid\'erons l'endomorphisme 
$f(z,w)=(z^d,w^d)$ de $\C^2$ et sa restriction sur l'ouvert 
$U:=\{1/2<|z|, |w|<2\}$. On v\'erifie que son degr\'e
dynamique (local) introduit dans la d\'efinition 3.1.3, 
est 1. Il est strictement plus petit que le degr\'e
dynamique global $d$. Dans cet exemple, 
le th\'eor\`eme 3.6.1 donne la meilleure
estimation possible pour les exposants de Lyapounov. 
\end{exemple}
\small

\par\noindent
Tien-Cuong Dinh et Nessim Sibony \\
Math\'ematique - B\^at. 425, UMR 8628,\\
Universit\'e Paris-Sud (Paris 11), 
91405 Orsay, France. \\
E-mails: Tiencuong.Dinh@math.u-psud.fr
et Nessim.Sibony@math.u-psud.fr
\end{document}